\documentclass[aos]{imsart}

\usepackage{amsfonts,amssymb,amsmath,amsthm,amscd,bbm,enumitem}
\usepackage[mathscr]{euscript} 
\usepackage[round]{natbib}
\usepackage{color,graphicx,hyperref}

\usepackage{apptools}
\AtAppendix{\counterwithin{lemma}{section}}
\AtAppendix{\counterwithin{theorem}{section}}
\AtAppendix{\counterwithin{assumption}{section}}
\AtAppendix{\counterwithin{example}{section}}

\startlocaldefs
\newtheorem{lemma}{Lemma}

\newtheorem{theorem}{Theorem}
\newtheorem{corollary}{Corollary}
\newtheorem{example}{Example}
\newtheorem{assumption}{Assumption}

\newtheorem*{assumption-BW}{Assumption BW}
\newtheorem*{assumption-MD}{Assumption MD}
\newtheorem*{assumption-MR}{Assumption MR}
\newtheorem*{assumption-MRC}{Assumption MRC}

\DeclareMathOperator*{\argmin}{argmin}
\DeclareMathOperator*{\argmax}{argmax}

\DeclareMathOperator*{\sign}{sign}

\newcommand{\1}{\mathbbm{1}}
\renewcommand\P{\mathbbm{P}}
\newcommand{\E}{\mathbbm{E}}

\newcommand{\Cov}{\mathbbm{C}\mathrm{ov}}
\newcommand\R{\mathbb{R}}
\newcommand\N{\mathbb{N}}
\newcommand\indep{\protect\mathpalette{\protect\independenT}{\perp}}
\def\independenT#1#2{\mathrel{\rlap{$#1#2$}\mkern2mu{#1#2}}}

\newcommand{\G}{\mathbbm{G}}
\newcommand\x{\mathsf{x}}
\newcommand\q{\mathfrak{q}}
\newcommand{\s}{\mathfrak{s}}

\newcommand{\GCM}{\mathsf{GCM}}
\newcommand{\LSC}{\mathsf{LSC}}

\begin{document}

\begin{frontmatter}

\title{Bootstrap-Assisted Inference for Generalized Grenander-type Estimators} 
\runtitle{Bootstrap-Assisted Inference for Monotone Estimators}

\begin{aug}
\author[A]{\fnms{Matias D.}~\snm{Cattaneo}\ead[label=e1]{cattaneo@princeton.edu}},
\author[B]{\fnms{Michael}~\snm{Jansson}\ead[label=e2]{mjansson@econ.berkeley.edu}}
\and
\author[C]{\fnms{Kenichi}~\snm{Nagasawa}\ead[label=e3]{kenichi.nagasawa@warwick.ac.uk}}
\address[A]{Department of Operations Research and Financial Engineering,
Princeton University\printead[presep={,\ }]{e1}}

\address[B]{Department of Economics,
University of California at Berkeley\printead[presep={,\ }]{e2}}

\address[C]{Department of Economics,
University of Warwick\printead[presep={,\ }]{e3}}
\end{aug}

\begin{abstract}

\citet{Westling-Carone_2020_AoS} proposed a framework for studying the large sample distributional properties of generalized Grenander-type estimators, a versatile class of nonparametric estimators of monotone functions. The limiting distribution of those estimators is representable as the left derivative of the greatest convex minorant of a Gaussian process whose monomial mean can be of unknown order (when the degree of flatness of the function of interest is unknown). The standard nonparametric bootstrap is unable to consistently approximate the large sample distribution of the generalized Grenander-type estimators even if the monomial order of the mean is known, making statistical inference a challenging endeavour in applications. To address this inferential problem, we present a bootstrap-assisted inference procedure for generalized Grenander-type estimators. The procedure relies on a carefully crafted, yet automatic, transformation of the estimator. Moreover, our proposed method can be made ``flatness robust'' in the sense that it can be made adaptive to the (possibly unknown) degree of flatness of the function of interest. The method requires only the consistent estimation of a single scalar quantity, for which we propose an automatic procedure based on numerical derivative estimation and the generalized jackknife. Under random sampling, our inference method can be implemented using a computationally attractive exchangeable bootstrap procedure. We illustrate our methods with examples and we also provide a small simulation study. The development of formal results is made possible by some technical results that may be of independent interest.
    
\end{abstract}



\end{frontmatter}

\section{Introduction}

Monotone function estimators have received renewed attention in statistics, biostatistics, econometrics, machine learning, and other data science disciplines. See \citet{Groeneboom-Jongbloed_2014_Book,Groeneboom-Jongbloed_2018_SS} for a textbook introduction and a review article, respectively, the latter being published in a special issue devoted to nonparametric inference under shape constraints. More recently, \citet{Westling-Carone_2020_AoS} expanded the scope and applicability of monotone function estimators by embedding many such estimators in a unified framework of generalized Grenander-type estimators. Estimation problems covered by \citeauthor{Westling-Carone_2020_AoS}'s \citeyearpar{Westling-Carone_2020_AoS} general theory include many practically relevant examples such as monotone density, regression and hazard estimation, possibly with censoring and/or covariate adjustment.

The large sample theory developed by \citet{Westling-Carone_2020_AoS} offers a general distributional approximation involving the left derivative of the Greatest Convex Minorant (GCM) of a Gaussian process whose mean and covariance kernel depend on unknown functions. Furthermore, both the convergence rate of the estimator and the shape of the mean appearing in the representation of the limiting distribution depend on whether the unknown function of interest exhibits certain degeneracies. For these reasons, the large sample distributional approximation for generalized Grenander-type estimators can be difficult to employ in practice for inference purposes. In their concluding remarks, \citet{Westling-Carone_2020_AoS} recognize these limitations and pose the question of whether it would be possible to employ bootstrap-assisted methods to conduct automatic/robust statistical inference within their framework.

As is well documented, the standard nonparametric bootstrap does not provide a valid distributional approximation for the generalized Grenander-type estimators \citep*{Kosorok_2008_BookCh,Sen-Banerjee-Woodroofe_2010_AoS}. This fact has led scholars to rely on other bootstrap schemes such as subsampling \citep{Politis-Romano_1994_AoS}, the smoothed bootstrap \citep{Kosorok_2008_BookCh,Sen-Banerjee-Woodroofe_2010_AoS}, the $m$-out-of-$n$ bootstrap \citep{Sen-Banerjee-Woodroofe_2010_AoS,Lee-Yang_2020_AoS}, or the numerical bootstrap \citep{Hong-Li_2020_AoS}. See also \citet{Cavaliere-Georgiev_2020_ECMA}, and references therein, for some related recent results. Those approaches could in principle be used to construct bootstrap-based inference methods for (some members of the class of) generalized Grenander-type estimators, but they all would require employing specific regularized multidimensional bootstrap distributions or related quantities involving multiple smoothing and tuning parameters, rendering those approaches potentially difficult to implement in practice. Furthermore, those methods would not be robust to unknown degeneracies determining the convergence rate and shape of the limiting distribution without additional modifications. For example, subsampling methods require knowledge of the precise convergence rate of the statistic, or estimation thereof, as a preliminary step \citep*{Politis-Romano-Wolf_1999_Book}. Another resampling approach was recently proposed by \citet{kuchibhotla-Balakrishnan-Wasserman_2024_JRSSB} and \citet{mallick2024new}, which offers valid confidence intervals under an approximate symmetry condition of the limiting distribution of the estimator.

We complement existing methods by introducing a novel bootstrap-assisted inference approach that restores validity of bootstrap methods by reshaping of the ingredients of the generalized Grenander-type estimator. Our approach is motivated by a constructive interpretation of the source of the bootstrap failure. As a by-product, our interpretation explicitly isolates the role of unknown degeneracies determining the precise form of the limiting distribution, a feature of the interpretation which allow us to develop an automatic inference method that is robust to such degeneracies, ultimately resulting in a more robust bootstrap-assisted inference approach. In the case of random sampling, we show that our method can be implemented using a computationally attractive exchangeable bootstrap procedure. For completeness, we also discuss implementation issues, offering a fully automatic (i.e., data-driven) valid inference method for generalized Grenander-type estimators. Some of the ideas underlying our approach are similar to ideas used in \citet*{Cattaneo-Jansson-Nagasawa_2020_ECMA}, where we introduced a bootstrap-based distributional approximation for $M$-estimators with possibly non-standard \citet{Chernoff_1964_AISM}-type asymptotic distributions \citep{Kim-Pollard_1990_AoS,Seo-Otsu_2018_AoS}. Generalized Grenander-type estimators are not $M$-estimators, however, and as further explained in the next paragraph the analysis of generalized Grenander-type estimators turns out to necessitate the development of technical tools that play no role in the analysis of $M$-estimators.

Although valid distributional approximations for monotone estimators can be obtained in a variety of ways \citep{Groeneboom-Jongbloed_2014_Book}, by far the most common approach is to employ the \textit{switch relation} \citep{Groeneboom_1985_BookCh} to re-express the cumulative distribution function (cdf) of the suitably normalized monotone estimator in terms of a probability statement about the maximizer of a certain stochastic process whose large sample properties can in turn be analyzed by employing standard empirical process methods \citep{vanderVaart-Wellner_1996_Book}. From a technical perspective, this approach requires (at least) two ingredients in order to be successful, namely results establishing (i) validity of the switch relation and (ii) continuity of the limiting cdf of the maximizer of a stochastic process. In the process of developing our main results, we shed new light on both (i) and (ii). First, we show by example that \citeauthor{Westling-Carone_2020_AoS}'s \citeyearpar[Supplement]{Westling-Carone_2020_AoS} generalization of the switch relation is incomplete as stated and then propose a modification. Our modification of \citeauthor{Westling-Carone_2020_AoS}'s \citeyearpar[Supplement]{Westling-Carone_2020_AoS} Lemma 1 shows that the conclusion of that lemma is valid once an additional assumption is made. The additional assumption seems very mild, being satisfied in all examples considered by \citet{Westling-Carone_2020_AoS} and all other examples of which we are aware, including some new examples we consider in this paper. Second, we present a lemma establishing continuity of the cdf of the maximizer of a Gaussian process, a result which can in turn be used to establish continuity of the cdf of the suitably normalized generalized Grenander-type estimator. Interestingly, although these continuity properties are important when deriving limiting distributions with the help of the switch relation and justifying bootstrap-type inference procedures, respectively, it would appear that explicit statements of them are unavailable in the existing literature. (A prominent exception is the one where the limiting distribution is a scaled Chernoff distribution, which is known to be absolutely continuous.) For further details on (i) and (ii), see Appendix \ref{[Section] Appendix - Technical Results}.

In the remainder of this introductory section we outline key notation and definitions used throughout the paper. Section \ref{[Section] Setup} then recalls the setup of \citet{Westling-Carone_2020_AoS} and presents a version of their main distributional approximation result for generalized Grenander-type estimators. Section \ref{[Section] Bootstrap-based Distributional Approximation} contains our main results about bootstrap-assisted distributional approximations, while Section \ref{[Section] Implementation} discusses implementation issues, including tuning parameter selection and a computationally attractive weighted bootstrap procedure. Section \ref{[Section] Examples} illustrates our general theory by means of prominent examples, while Section \ref{[Section] Simulations} reports numerical results from a small-scale simulation experiment. Appendix \ref{[Section] Appendix - Technical Results} reports the two technical results alluded to in the previous paragraph. All proofs and other technical details are given in the supplemental appendix \citep{Cattaneo-Jansson-Nagasawa_2024_supp}.

\subsection{Notation and Definitions}

For an interval $I \subseteq \R$ and a function $f : I \to \R$, $\GCM_I(f)$ denotes its greatest convex minorant (on $I$) and if $f$ is non-decreasing and right-continuous, then $f^-$ denotes its generalized inverse; that is, $f^-(x) = \inf\{u\in I : f(u) \geq x\}$, where the dependence of $f^-$ on $I$ has been suppressed (and where the infimum of the empty set is $\sup I$). Also, if $f$ is non-decreasing and right-continuous, we write $f(x-) = \lim_{\eta \downarrow 0}f(x - \eta)$ and $f^-(y+) = \lim_{\eta \downarrow 0}f^-(y + \eta)$. Assuming the relevant derivatives exist, $\partial^q$ denotes the $q$th partial derivative (operator) and $\partial_-$ denotes the left derivative (operator). In addition, $f\circ g$ denotes the composition of $f$ and $g$; that is, $(f\circ g)(x) = f(g(x))$. Finally, for $a,b \in \R$, we write $a \land b = \min(a,b), a \lor b = \max(a,b)$, and $a^+ = \max(a,0)$, and $\lfloor a \rfloor$ denotes the integer part of $a$.

Limits are taken as $n \to \infty$, unless otherwise stated. For two (possibly) random sequences $\{a_n\}$ and $\{b_n\}$, $a_n = O_\P(b_n)$ is shorthand for $\limsup_{\epsilon \to \infty}\limsup_n\P[|a_n / b_n| \geq \epsilon] = 0$, $a_n = o_\P(b_n)$ is shorthand for $\limsup_{\epsilon \to 0}\limsup_n\P[|a_n / b_n| \geq \epsilon] = 0$, and the subscript ``$\P$'' on ``$O$'' and ``$o$'' is often omitted when $\{a_n\}$ and $\{b_n\}$ are non-random. We use $\to_\P$ to denote convergence in probability and $\rightsquigarrow$ to denote weak convergence, where, for a stochastic process indexed by $\R$, convergence is in the topology of uniform convergence on compacta. When analyzing the bootstrap, $\P^*_n$ denotes the probability measure under the bootstrap distribution conditional on the original data and $\rightsquigarrow_\P$ denotes weak convergence in probability conditionally on the original data. For more details, see \citet{vanderVaart-Wellner_1996_Book}.

\section{Setup}\label{[Section] Setup}

Our setup is that of \citet{Westling-Carone_2020_AoS}. The goal is to conduct inference on $\theta_0(\x)$, where, for some interval $I \subseteq \R$, $\theta_0$ is non-decreasing on $I$ and $\x$ is an interior point of $I$. Assuming it is well defined, the function $\Theta_0$ given by
\[\Theta_0(x) = \int_{\inf I}^x \theta_0(v) dv\]
is convex on $I$ and therefore enjoys the property that if $\theta_0$ is left-continuous at $\x$, then
\[\theta_0(\x) = \partial_- \GCM_I(\Theta_0)(\x).\]
An estimator of $\theta_0(\x)$ obtained by replacing $\Theta_0$ and $I$ in the preceding display with estimators is said to be of the Grenander-type, a canonical example of this class of estimators being the celebrated \citet{Grenander_1956_SAJ} estimator of a non-decreasing density.

\begin{example}[Monotone Density Estimation]\label{[Example] Monotone Density Estimation}
Suppose $X_1,\dots,X_n$ are i.i.d. copies of a continuously distributed non-negative random variable $X$ whose density $f_0$ is non-decreasing on its support $[0,u_0],$ where $u_0$ is (possibly) unknown. For $\x \in (0,u_0)$, the \citet{Grenander_1956_SAJ} estimator of $f_0(\x)$ is
\[\widehat{f}_n(\x) = \partial_- \GCM_{[0,\widehat{u}_n]}(\widehat{F}_n)(\x),\]
where $\widehat{u}_n = \max_{1 \leq i \leq n} X_i \lor \x$ and where $\widehat{F}_n(x) = n^{-1} \sum_{i=1}^n \1(X_i \leq x)$ is the empirical cdf. Section \ref{[Subsection] Examples: Monotone Density Estimation} and the supplemental appendix discuss more general monotone density estimators allowing for censoring and/or covariate adjustment \citep[e.g.,][]{vanderlaan-Robins_2003_Book}.
\end{example}

To define the class of generalized Grenander-type estimators, let $\psi_0 = \theta_0 \circ \Phi_0^-$, where $\Phi_0$ is non-negative, non-decreasing, and continuous on $I$. Defining
\[\Gamma_0 = \Psi_0 \circ \Phi_0, \qquad \Psi_0(x) = \int_0^x \psi_0(v) dv,\]
and assuming that $\Phi_0(x) < \Phi_0(\x) < u_0$ for every $x < \x$, we have
\begin{equation}\label{[Representation] theta0}
\theta_0(\x) = \partial_- \GCM_{[0,u_0]} (\Gamma_0 \circ \Phi_0^-) \circ \Phi_0(\x)
\end{equation}
whenever $\theta_0$ is left-continuous at $\x$. In the terminology of \citet{Westling-Carone_2020_AoS}, an estimator of $\theta_0(\x)$ is of the Generalized Grenander-type if it is obtained by replacing $\Gamma_0$, $\Phi_0$, and $u_0$ in the preceding display with estimators $\widehat{\Gamma}_n,\widehat{\Phi}_n$, and $\widehat{u}_n$ (say); to be specific, an estimator of the Generalized Grenander-type is of the form
\[\widehat{\theta}_n(\x) = \partial_- \GCM_{[0,\widehat{u}_n]}(\widehat{\Gamma}_n \circ \widehat{\Phi}_n^-) \circ \widehat{\Phi}_n(\x),\]
where $\widehat{\Phi}_n$ is non-negative, non-decreasing, and right-continuous. Of course, Grenander-type estimators are of the generalized Grenander-type (with $\widehat{\Phi}_n$ equal to the identity mapping) whenever the associated estimator of $I$ is of the form $[0,\widehat{u}_n]$, but the class of generalized Grenander-type estimators contains many important estimators that are not of Grenander-type, a canonical example being the celebrated isotonic regression estimator of \citet{Brunk_1958_AMS}.

\begin{example}[Monotone Regression Estimation]\label{[Example] Monotone Regression Estimation}
Suppose $(Y_1,X_1),\dots,(Y_n,X_n)$ are i.i.d. copies of $(Y,X)$, where $X$ is a continuously distributed random variable and where the regression function $\mu_0(x) = \E(Y|X = x)$ is non-decreasing. For $\x$ in the interior of the support of $X$, the \citet{Brunk_1958_AMS} estimator of $\mu_0(\x)$ is
\[\widehat{\mu}_n(\x) = \partial_-\GCM_{[0,1]}(\widehat{\Gamma}_n \circ \widehat{F}_n^-) \circ \widehat{F}_n(\x),\]
where $\widehat{\Gamma}_n(x) = n^{-1} \sum_{i=1}^n Y_i \1(X_i \leq x)$ and where $\widehat{F}_n(x) = n^{-1} \sum_{i=1}^n \1(X_i \leq x)$ is an estimator of $F_0$, the cdf of $X$. Section \ref{[Subsection] Examples: Monotone Regression Estimation} discusses more general monotone regression estimators allowing for covariate adjustment \citep*[e.g.,][]{Westling-Gilbert-Carone_2020_JRRSB}.
\end{example}

Under regularity conditions, if $\widehat{\theta}_n(\x)$ is of the generalized Grenander-type, then its rate of convergence is governed by the flatness of $\theta_0$ around $\x$, as measured by the characteristic exponent
\[\q = \min\{j \in \N: \partial^j \theta_0(\x) \neq 0\},\]
where $\N$ is the set of positive integers. (When $\theta_0$ is non-decreasing and suitably smooth, $\q$ is necessarily an odd integer.) To be specific, \citet[Theorem 3]{Westling-Carone_2020_AoS} gave conditions under which
\begin{equation}\label{[Asymptotics] thetahat}
r_n(\widehat{\theta}_n(\x) - \theta_0(\x)) \rightsquigarrow \frac{1}{\partial\Phi_0(\x)} \partial_- \GCM_\R(\mathcal{G}_\x + \mathcal{M}_\x^\q)(0), \qquad r_n = n^{\q / (1 + 2 \q)},
\end{equation}
where $\mathcal{G}_\x$ is a scalar multiple of two-sided Brownian motion and where $\mathcal{M}_\x^\q$ is a monomial given by
\begin{equation}\label{[Definition] M_x,q}
\mathcal{M}_\x^\q(v) = \frac{\partial^\q \theta_0(\x) \partial \Phi_0(\x)}{(\q + 1)!} v^{\q + 1}.
\end{equation}
In addition to governing the rate of convergence, the characteristic exponent $\q$ also governs the shape of $\mathcal{M}_\x^\q$. On the other hand, and as the notation suggests, the covariance kernel of $\mathcal{G}_\x$ does not depend on $\q$. If $\q = 1$, then the distribution of $\partial_- \GCM_\R(\mathcal{G}_\x + \mathcal{M}_\x^\q)(0)$ is a scaled Chernoff distribution. More generally, the distribution of $\partial_- \GCM_\R(\mathcal{G}_\x + \mathcal{M}_\x^\q)(0)$ is a scaled Chernoff-type distribution \citep[in the terminology of][]{Han-Kato_2022_AoAP}.

Among other things, the following assumption guarantees validity of the representation (\ref{[Representation] theta0}) and ensures existence of the right hand side of (\ref{[Definition] M_x,q}).

\begin{assumption}\label{Assumption A}
For some $\delta > 0$ and some $\s \geq \q$, the following are satisfied:
\begin{enumerate}[label=(A\arabic*)]
\item $I \subseteq \R$ is an interval and $I_\x^\delta = \{x \in \R : |x - \x| \leq \delta\} \subseteq I$.

\item \label{Assumption A: theta0} $\theta_0$ is non-decreasing and bounded on $I$, and $\left\lfloor \s \right\rfloor$ times continuously differentiable on $I_\x^\delta$ with
\[\sup_{x,x' \in I_\x^\delta}\frac{|\partial^{\left\lfloor \s \right\rfloor} \theta_0(x) - \partial^{\left\lfloor \s \right\rfloor} \theta_0(x')|}{|x - x'|^{\s - \left\lfloor \s \right\rfloor}} < \infty.\]

\item \label{Assumption A: Phi0}
 $\Phi_0$ is non-negative, non-decreasing, continuous, and bounded on $I$, and $\left\lfloor \s \right\rfloor - \q + 1$ times continuously differentiable on $I_\x^\delta$ with $\partial\Phi_0(\x) \neq 0$ and
\[\sup_{x,x' \in I_\x^\delta} \frac{|\partial^{\left\lfloor \s \right\rfloor - \q + 1} \Phi_0(x) - \partial^{\left\lfloor \s \right\rfloor - \q + 1}\Phi_0(x')|}{|x - x'|^{\s - \left\lfloor \s \right\rfloor}} < \infty.\]

\end{enumerate}
\end{assumption}

\section{Bootstrap-Assisted Distributional Approximation}\label{[Section] Bootstrap-based Distributional Approximation}

The distributional approximation \eqref{[Asymptotics] thetahat} motivates a plug-in procedure for statistical inference: with consistent estimators of  $\mathcal{M}_\x^\q$ and the covariance kernel of $\mathcal{G}_\x$, the limiting law in \eqref{[Asymptotics] thetahat} can be simulated. This approach requires characterizing the covariance kernel explicitly, and then forming a consistent estimator based on preliminary non-parametric smoothing methods (in addition to the other preliminary nonparametric estimators needed). Also, simulating from the limiting law requires specifying or estimating the unknown characteristic exponent $\q$. As a consequence, although conceptually feasible, the plug-in approach is often computationally cumbersome to employ in practice.

Since resampling methods have the potential to provide automatic and robust distributional approximations, it seems natural to consider such methods as an alternative to plug-in methods. Letting $(\widehat{\Gamma}_n^*,\widehat{\Phi}_n^*,\widehat{u}_n^*)$ denote a generic (not necessarily nonparametric) bootstrap analog of $(\widehat{\Gamma}_n,\widehat{\Phi}_n,\widehat{u}_n)$ and assuming that $\widehat{\Phi}_n^*$ is non-decreasing and right-continuous, the associated bootstrap analog of $\widehat{\theta}_n(\x)$ is
\[\widehat{\theta}_n^*(\x) = \partial_- \GCM_{[0,\widehat{u}_n^*]}(\widehat{\Gamma}_n^* \circ \widehat{\Phi}_n^{*-}) \circ \widehat{\Phi}_n^*(\x).\]

As is well documented for the classical Grenander estimator \citep[e.g.][]{Kosorok_2008_BookCh,Sen-Banerjee-Woodroofe_2010_AoS}, the bootstrap analog of (\ref{[Asymptotics] thetahat}) does not (necessarily) hold when $(\widehat{\Gamma}_n^*,\widehat{\Phi}_n^*,\widehat{u}_n^*)$ is obtained by means of the nonparametric bootstrap. In other words,
\[r_n (\widehat{\theta}_n^*(\x) - \widehat{\theta}_n(\x)) \not \rightsquigarrow_\P \frac{1}{\partial\Phi_0(\x)} \partial_- \GCM_\R(\mathcal{G}_\x + \mathcal{M}_\x^\q)(0)\]
in general; that is, the bootstrap is inconsistent in general. (A precise statement is provided in Theorem \ref{[Theorem]: Bootstrap inconsistency} in Appendix \ref{[Subsection] Appendix - Bootstrap Inconsistency}.) It turns out, however, that under plausible conditions on $(\widehat{\Gamma}_n^*,\widehat{\Phi}_n^*,\widehat{u}_n^*)$, a valid bootstrap-based distributional approximation can be obtained by employing
\[\widetilde{\theta}_n^*(\x) = \partial_- \GCM_{[0,\widehat{u}_n^*]}(\widetilde{\Gamma}_n^* \circ \widehat{\Phi}_n^{*-}) \circ \widehat{\Phi}_n^*(\x),\]
where, for some judiciously chosen $\widetilde{M}_{\x,n}$,
\[\widetilde{\Gamma}_n^*(x) = \widehat{\Gamma}_n^*(x) - \widehat{\Gamma}_n(x) + \widehat{\theta}_n(\x) \widehat{\Phi}_n(x) + \widetilde{M}_{\x,n}(x - \x).\]
(As the notation suggests, and as further discussed below, a defining property of $\widetilde{M}_{\x,n}$ is that a suitably re-scaled version of it can be interpreted as an estimator of $\mathcal{M}_\x^\q$.)

\subsection{Heuristics}\label{[Subsection] Heuristics}

To explain the source of the bootstrap failure and motivate the functional form of $\widetilde{\Gamma}_n^*$, it is useful to begin by sketching the derivation of (\ref{[Asymptotics] thetahat}). For any $t \in \R$ and under mild conditions on $\widehat{\Gamma}_n$ and $\widehat{\Phi}_n$, it follows from the so-called (generalized) switch relation (e.g., Appendix \ref{[Subsection] Appendix - Generalized Switch Relation}) that
\[\P \left[r_n (\widehat{\theta}_n(\x) - \theta_0(\x)) \leq t\right] = \P \left[\argmin_{v \in \widehat{V}_{\x,n}^\q} \{\widehat{G}_{\x,n}^\q(v) + M_{\x,n}^\q(v) - t \widehat{L}_{\x,n}^\q(v)\} \geq \widehat{Z}_{\x,n}^\q \right],\]
where $\widehat{V}_{\x,n}^\q$ is a ``large'' subset of $\R$, $\widehat{G}_{\x,n}^\q$ is a stochastic process given by

\begin{align*}
\widehat{G}_{\x,n}^\q(v) &= \sqrt{n a_n} [\widehat{\Gamma}_n(\x + v a_n^{-1}) - \widehat{\Gamma}_n(\x) - \Gamma_0(\x + v a_n^{-1}) + \Gamma_0(\x)]\\
&- \theta_0(\x) \sqrt{n a_n} [\widehat{\Phi}_n(\x + v a_n^{-1}) - \widehat{\Phi}_n(\x) - \Phi_0(\x + v a_n^{-1}) + \Phi_0(\x)], \quad a_n = n^{1 / (1 + 2 \q)},
\end{align*}
$M_{\x,n}^\q$ is a smooth non-random function, $\widehat{L}_{\x,n}^\q$ can be interpreted as an estimator of the linear function $\mathcal{L}_\x$ given by $\mathcal{L}_\x(v) = v \partial \Phi_0(\x)$,
and where $\widehat{Z}_{\x,n}^\q$ is zero when $\widehat{\Phi}_n$ is the identity mapping and can be interpreted as an ``estimator'' of zero more generally. (The exact definitions of $\widehat{V}_{\x,n}^\q,M_{\x,n}^\q,\widehat{L}_{\x,n}^\q$, and $\widehat{Z}_{\x,n}^\q$ are given in Appendix \ref{[Subsection] Appendix - Omitted Formulas}.) In the above and the sequel, the superscript $\q$ indicates that the superscripted object is defined on a ``localized'' domain. 

Under mild conditions, $\widehat{G}_{\x,n}^\q$ converges weakly to a scalar multiple of two-sided Brownian motion $\mathcal{G}_\x$ (say), while $M_{\x,n}^\q,\widehat{L}_{\x,n}^\q$, and $\widehat{Z}_{\x,n}^\q$ converge in probability to $\mathcal{M}_\x^\q,\mathcal{L}_\x$, and zero, respectively. In other words, we would expect that
\begin{equation}\label{[Asymptotics] Ghat,Lhat,Zhat}        
(\widehat{G}_{\x,n}^\q,\widehat{L}_{\x,n}^\q,\widehat{Z}_{\x,n}^\q) \rightsquigarrow (\mathcal{G}_\x,\mathcal{L}_\x,0),
\end{equation}
and that
\begin{equation}\label{[Asymptotics] M}
M_{\x,n}^\q \rightsquigarrow \mathcal{M}_\x^\q.
\end{equation}
It therefore stands to reason that under mild additional conditions (including suitable convergence of $\widehat{V}_{\x,n}^\q$ to $\R$), we have
\begin{align*}
\lim_{n \to \infty} \P \left[r_n (\widehat{\theta}_n(\x) - \theta_0(\x)) \leq t\right] &= \P \left[\argmin_{v \in \R}\{\mathcal{G}_\x(v) + \mathcal{M}_\x^\q(v) - t \mathcal{L}_\x(v)\} \geq 0 \right]\\
&= \P \left[\frac{1}{\partial \Phi_0(\x)} \partial_- \GCM_\R(\mathcal{G}_\x + \mathcal{M}_\x^\q)(0) \leq t\right],
\end{align*}
where the first equality follows from an $\argmax$ continuous mapping theorem, and the second equality is obtained by another application of the switch relation.

Next, consider the bootstrap analog of $\widehat{\theta}_n(\x)$. For any $t \in \R$, the switch relation gives
\[\P_n^*\left[r_n (\widehat{\theta}_n^*(\x) - \widehat{\theta}_n(\x)) \leq t\right] = \P_n^*\left[\argmin_{v \in \widehat{V}_{\x,n}^{\q,*}}\{\widehat{G}_{\x,n}^{\q,*}(v) + \widehat{M}_{\x,n}^\q(v) - t \widehat{L}_{\x,n}^{\q,*}(v)\} \geq \widehat{Z}_{\x,n}^{\q,*} \right],\]
where
\begin{align*}
\widehat{G}_{\x,n}^{\q,*}(v) &= \sqrt{n a_n} [\widehat{\Gamma}_n^*(\x + v a_n^{-1}) - \widehat{\Gamma}_n^*(\x) - \widehat{\Gamma}_n(\x + v a_n^{-1}) + \widehat{\Gamma}_n(\x)] \\
&-\widehat{\theta}_n(\x) \sqrt{n a_n} [\widehat{\Phi}_n^*(\x + v a_n^{-1}) - \widehat{\Phi}_n^*(\x) -\widehat{\Phi}_n(\x + v a_n^{-1}) + \widehat{\Phi}_n(\x)],
\end{align*}
and where $\widehat{V}_{\x,n}^{\q,*},\widehat{M}_{\x,n}^\q,\widehat{L}_{\x,n}^{\q,*}$, and $\widehat{Z}_{\x,n}^{\q,*}$ are bootstrap analogs of $\widehat{V}_{\x,n}^\q,M_{\x,n}^\q,\widehat{L}_{\x,n}^\q$, and $\widehat{Z}_{\x,n}^\q$, respectively. (The exact definitions of $\widehat{V}_{\x,n}^{\q,*},\widehat{M}_{\x,n}^\q,\widehat{L}_{\x,n}^{\q,*}$, and $\widehat{Z}_{\x,n}^{\q,*}$ are given in Appendix \ref{[Subsection] Appendix - Omitted Formulas}.) If (\ref{[Asymptotics] Ghat,Lhat,Zhat}) holds, then the following bootstrap counterpart thereof can also be expected to hold:
\begin{equation}\label{[Asymptotics] Ghatstar,Lhatstar,Zhatstar}
(\widehat{G}_{\x,n}^{\q,*},\widehat{L}_{\x,n}^{\q,*},\widehat{Z}_{\x,n}^{\q,*}) \rightsquigarrow_\P (\mathcal{G}_\x,\mathcal{L}_\x,0).
\end{equation}
On the other hand, unlike $M_{\x,n}^\q$ the process $\widehat{M}_{\x,n}^\q$ is (random and) non-smooth, so the bootstrap counterpart of (\ref{[Asymptotics] M}) typically fails, implying in turn that the bootstrap is inconsistent; for details, see Theorem \ref{[Theorem]: Bootstrap inconsistency} (stated in Appendix \ref{[Subsection] Appendix - Bootstrap Inconsistency}) and its proof.

In other words, the sole source of the bootstrap inconsistency is the failure of $\widehat{M}_{\x,n}^\q$ to replicate the salient properties of $M_{\x,n}^\q$. Recognizing this, our proposed estimator $\widetilde{\theta}_n^*(\x)$ has been carefully constructed so that it differs from $\widehat{\theta}_n^*(\x)$ only in terms of the implied estimator of $\mathcal{M}_\x^\q$. To be specific, $\widetilde{\theta}_n^*(\x)$ is similar to $\widehat{\theta}_n^*(\x)$ insofar as it satisfies the following switch relation: For any $t\in\R$,
\[\P_n^* \left[r_n (\widetilde{\theta}_n^*(\x) -\widehat{\theta}_n(\x)) \leq t \right] = \P_n^* \left[\argmin_{v \in \widehat{V}_{\x,n}^{\q,*}} \{\widehat{G}_{\x,n}^{\q,*}(v) + \widetilde{M}_{\x,n}^\q(v) - t \widehat{L}_{\x,n}^{\q,*}(v)\} \geq \widehat{Z}_{\x,n}^{\q,*}\right],\]
where $\widetilde{M}_{\x,n}^\q$ is the following transformation of $\widetilde{M}_{\x,n}$:
\[\widetilde{M}_{\x,n}^\q(v) = \sqrt{n a_n} \widetilde{M}_{\x,n}(v a_n^{-1}).\]
As a consequence, if (\ref{[Asymptotics] Ghatstar,Lhatstar,Zhatstar}) holds and if $\widetilde{M}_{\x,n}^\q \rightsquigarrow_\P \mathcal{M}_\x^\q$, then it stands to reason that under mild additional conditions $\widetilde{\theta}_n^*(\x)$ satisfies the following bootstrap counterpart of (\ref{[Asymptotics] thetahat}):
\begin{equation}\label{[Asymptotics] thetatildestar}
r_n (\widetilde{\theta}_n^*(\x) - \widehat{\theta}_n(\x)) \rightsquigarrow_\P \frac{1}{\partial\Phi_0(\x)} \partial_- \GCM_\R (\mathcal{G}_\x + \mathcal{M}_\x^\q)(0).
\end{equation}

In words, our proposed bootstrap-assisted inference procedure can be explained as follows: we employ the bootstrap to approximate the ``fluctuations'' of the process $\mathcal{G}_\x$, while we use the estimator $\widetilde{M}_{\x,n}^\q$ to approximate $\mathcal{M}_\x^\q$. It is an interesting question whether our approach offers provably better distributional approximations than the plug-in procedure mentioned previously, but it is beyond the scope of this paper to answer that question.

\subsection{Main result}\label{[Subsection] Main result}

Our heuristic derivation of (\ref{[Asymptotics] thetahat}) can be made rigorous by providing conditions under which four properties hold. First, the switch relation(s) must be valid. Second, the convergence properties (\ref{[Asymptotics] Ghat,Lhat,Zhat}) and (\ref{[Asymptotics] M}) must hold. Third, to use (\ref{[Asymptotics] Ghat,Lhat,Zhat}) and (\ref{[Asymptotics] M}) to obtain the result
\[\argmin_{v \in \widehat{V}_{\x,n}^\q} \{\widehat{G}_{\x,n}^\q(v) + M_{\x,n}^\q(v) - t \widehat{L}_{\x,n}^\q(v)\} \rightsquigarrow \argmin_{v \in \R} \{\mathcal{G}_\x(v) + \mathcal{M}_\x^\q(v) - t \mathcal{L}_\x(v)\}\]
with the help of a suitable $\argmax$ continuous mapping theorem, tightness of the left hand side in the previous display must hold. Finally, to furthermore obtain the conclusion
\begin{align*}
&\lim_{n \to \infty} \P \left[\argmin_{v \in \widehat{V}_{\x,n}^\q} \{\widehat{G}_{\x,n}^\q(v) + M_{\x,n}^\q(v) - t \widehat{L}_{\x,n}^\q(v)\} \geq \widehat{Z}_{\x,n}^\q \right] \\
&= \P \left[\argmin_{v \in \R}\{\mathcal{G}_\x(v) + \mathcal{M}_\x^\q(v) - t \mathcal{L}_\x(v)\} \geq 0 \right],
\end{align*}
the cdf of $\argmin_{v \in \R}\{\mathcal{G}_\x(v) + \mathcal{M}_\x^\q(v) - t \mathcal{L}_\x(v)\}$ must be continuous at zero.

Conditions under which the second and third properties hold can be formulated with the help of well known empirical process results. We base our formulations on \citet{vanderVaart-Wellner_1996_Book} and employ the generalized $\argmax$ continuous mapping theorem of \citet{Cox_2022}, which applies also to settings where the feasibility set varies with $n$. The first and fourth properties, on the other hand, seem more difficult to verify. Regarding the first property, it turns out that the generalization of the switch relation employed by \citet{Westling-Carone_2020_AoS} requires additional conditions in order to be valid. To address this concern, we present a modified generalized switch relation whose assumptions include a condition not present in Lemma 1 of \citet[Supplement]{Westling-Carone_2020_AoS}. Thankfully, the condition in question seems very mild and having imposed it we are able to preserve the main implication of Lemma 1 of \citet[Supplement]{Westling-Carone_2020_AoS}; for details, see Lemma \ref{[Lemma] Generalized Switch Relation} in Appendix \ref{[Subsection] Appendix - Generalized Switch Relation}. In the special case where $\q = 1$, the fourth property follows from well known properties of the Chernoff distribution. More generally, however, we are unaware of existing results guaranteeing the requisite continuity property when $\q\neq1$, but fortunately it turns out that the continuity property of interest can be established; for details, see Lemma \ref{[Lemma] Continuity of argmax} in Appendix \ref{[Subsection] Appendix - Continuity of Limiting Distribution}.

The following assumption collects the conditions under which our verification of the four above-mentioned properties will proceed.

\begin{assumption}\label{Assumption B}
For some $\delta>0$, the following are satisfied:
\begin{enumerate}[label=(B\arabic*)]

\item \label{Assumption B: weak convergence}
$\widehat{G}_{\x,n}^\q \rightsquigarrow \mathcal{G}_\x$ and $\widehat{G}_{\x,n}^{\q,*} \rightsquigarrow_\P \mathcal{G}_\x$, where $\mathcal{G}_\x$ is a centered Gaussian process whose covariance function is $\mathcal{C}_\x(s,t) = C(|s| \land |t|) \1(\sign(s)=\sign(t))$ for some $C>0$. 
    
\item \label{Assumption B: moment bound} 
There exist $\beta < \q+1$ and events $\mathcal{A}_n$ with $\lim_{n \to \infty}\P [\mathcal{A}_n] = 1$,
\[\sup_{V \in [1,a_n \delta]} \E \left[V^{-\beta} \sup_{|v| \in [V,2V]} |\widehat{G}_{\x,n}^\q(v)|\1_{\mathcal{A}_n}\right] = O(1)\]
and
\[\sup_{V \in [1,a_n \delta]} \E \left[V^{-\beta} \sup_{|v| \in [V,2V]} |\widehat{G}_{\x,n}^{\q,*}(v)|\1_{\mathcal{A}_n}\right] = O(1).\]

\item \label{Assumption B: uniform consistency Gamma}
$\sup_{x \in I}|\widehat{\Gamma}_n(x) - \Gamma_0(x)| = o_\P(1)$ and $\sup_{x \in I}|\widehat{\Gamma}_n^*(x) - \widehat{\Gamma}_n(x)| = o_\P(1)$.

\item \label{Assumption B: uniform convergence Phi}
$\sup_{x \in I} |\widehat{\Phi}_n(x) - \Phi_0(x)| = o_\P(1)$ and $\sup_{x \in I} |\widehat{\Phi}_n^*(x) - \widehat{\Phi}_n(x)| = o_\P(1).$
In addition, 
\[a_n \sup_{x \in I_\x^\delta} |\widehat{\Phi}_n(x) - \Phi_0(x)| = o_\P(1) \quad \text{and} \quad a_n \sup_{x \in I_\x^\delta} |\widehat{\Phi}_n^*(x) - \widehat{\Phi}_n(x)| = o_\P(1).\]

\item \label{Assumption B: uhat}
For some $u_0 > \Phi_0(\x)$, $\widehat{u}_n \geq u_0 + o_\P(1)$ and $\widehat{u}_n^* \geq \widehat{u}_n + o_\P(1)$.

\item \label{Assumption B: Phi uhat}
$\widehat{\Phi}_n$ and $\widehat{\Phi}_n^*$ are non-negative, non-decreasing, and right-continuous on $I$. In addition, $\{0,\widehat{u}_n\} \subseteq \widehat{\Phi}_n(I)$ and $\{0,\widehat{u}_n^*\} \subseteq \widehat{\Phi}_n^*(I)$. Also, $\widehat{\Phi}_n(I) \cap [0,\widehat{u}_n]$ and $\widehat{\Phi}_n^*(I) \cap [0,\widehat{u}_n^*]$ are closed. 

\item \label{Assumption B: jump of Phi hat}
\[\sqrt{n a_n} \sup_{x \in I_\x^\delta} |\widehat{\Phi}_n (x) - \widehat{\Phi}_n (x-)| = o_\P(1) \quad \text{and} \quad \sqrt{n a_n} \sup_{x \in I_\x^\delta} |\widehat{\Phi}_n^* (x) - \widehat{\Phi}_n^* (x-)| = o_\P(1).\]

\end{enumerate}    
\end{assumption}

Verification of the bootstrap parts of Assumption \ref{Assumption B} will be discussed in Section \ref{[Subsection] Bootstrapping} below. When combined with Assumption \ref{Assumption A}, Assumption \ref{Assumption B} suffices in order to establish (\ref{[Asymptotics] thetahat}) and (\ref{[Asymptotics] Ghatstar,Lhatstar,Zhatstar}). In addition, (\ref{[Asymptotics] thetatildestar}) can be shown to hold if $\widetilde{M}_{\x,n}$ satisfies the following

\begin{assumption}\label{Assumption C}
$\widetilde{M}_{\x,n}^\q \rightsquigarrow_\P \mathcal{M}_\x^\q$ and, for some $c>0$ and every $K>0$,
\[\liminf_{n \to \infty} \P \left[\inf_{|v| > K^{-1}} \widetilde{M}_{\x,n}(v) \geq c K^{-(\q + 1)} \right] = 1.\]
\end{assumption}

Moreover, continuity of the cdf of $\partial_- \GCM_\R (\mathcal{G}_\x + \mathcal{M}_\x^\q)(0)$ can be shown with the help of the arguments used to show that the cdf of $\argmin_{v \in \R} \{\mathcal{G}_\x(v) + \mathcal{M}_\x^\q(v) - t \mathcal{L}_\x(v)\}$ is continuous at zero. As a consequence, we obtain the following result.

\begin{theorem}\label{[Theorem] Main result}
Suppose Assumptions \ref{Assumption A}, \ref{Assumption B}, and \ref{Assumption C} are satisfied. Then (\ref{[Asymptotics] thetahat}) and (\ref{[Asymptotics] thetatildestar}) hold, and \begin{equation}\label{Bootstrap consistency}
\sup_{t \in \R} \left|\P_n^* \left[\widetilde{\theta}_n^*(\x) - \widehat{\theta}_n(\x) \leq t \right] - \P \left[\widehat{\theta}_n(\x) - \theta_0(\x) \leq t \right]\right| = o_\P(1).
\end{equation}
\end{theorem}

In an attempt to emphasize the rate-adaptive nature of the consistency property enjoyed by the bootstrap-based distributional approximation based on $\widetilde{\theta}_n^*(\x)$, the formulation (\ref{Bootstrap consistency}) deliberately omits the rate term $r_n$ present in (\ref{[Asymptotics] thetahat}) and (\ref{[Asymptotics] thetatildestar}). Theorem \ref{[Theorem] Main result} has immediate implications for inference. For instance, it follows from (\ref{[Asymptotics] thetahat}) and (\ref{Bootstrap consistency}) that for any $\alpha \in (0,1)$, we have
\[\lim_{n \to \infty} \P[\theta_0(\x) \in \mathsf{CI}_{1 - \alpha,n}(\x)] = 1 - \alpha,\]
where, defining $Q_{a,n}^*(\x) = \inf\{Q \in \R : \P_n^* [\widetilde{\theta}_n^*(\x) - \widehat{\theta}_n(\x) \leq Q] \geq a\}$,
\[\mathsf{CI}_{1 - \alpha,n}(\x) = \left[\widehat{\theta}_n(\x) - Q_{1-\alpha/2,n}^*(\x)~,~\widehat{\theta}_n(\x) - Q_{\alpha/2,n}^*(\x)\right]\]
is the (nominal) level $1 - \alpha$ bootstrap confidence interval for $\theta_0$ based on the ``percentile method'' \citep[in the terminology of][]{vanderVaart_1998_Book}.

\section{Implementation}\label{[Section] Implementation}

Suppose Assumption \ref{Assumption A} is satisfied and suppose $(\widehat{\Gamma}_n,\widehat{\Phi}_n,\widehat{u}_n)$ satisfies the non-bootstrap parts of Assumption \ref{Assumption B}. Then, in order to compute the estimator $\widetilde{\theta}_n^*(\x)$ upon which our proposed bootstrap-assisted distributional approximation is based, two implementation issues must be addressed, namely the choice/construction of $\widetilde{M}_{\x,n}$ and $(\widehat{\Gamma}_n^*,\widehat{\Phi}_n^*,\widehat{u}_n^*)$, respectively. Section \ref{[Subsection] Mean function estimation} demonstrates the plausibility of Assumption \ref{Assumption C} by exhibiting estimators $\widetilde{M}_{\x,n}$ satisfying it under Assumptions \ref{Assumption A} and \ref{Assumption D}, the latter being a high-level condition that typically holds whenever the non-bootstrap parts of Assumption \ref{Assumption B} hold. Then, Section \ref{[Subsection] Bootstrapping} exhibits easy-to-compute $(\widehat{\Gamma}_n^*,\widehat{\Phi}_n^*,\widehat{u}_n^*)$ satisfying the bootstrap parts of Assumption \ref{Assumption B} under a random sampling assumption and other mild conditions.

\subsection{Mean Function Estimation}\label{[Subsection] Mean function estimation}

The ease with which an $\widetilde{M}_{\x,n}$ satisfying Assumption \ref{Assumption C} can be constructed depends in particular on whether $\q$ is known. To facilitate the subsequent discussion, let
$\mathcal{D}_0(\x) = 0$ and for $j = 1,\dots,\lfloor \s \rfloor$, define (recursively)
\[\mathcal{D}_{j}(\x) = \lim_{x \to \x} \frac{\Upsilon_0(x) - \Upsilon_0(\x) - \sum_{k = 0}^{j-1} [(k + 1)!]^{-1} \mathcal{D}_k(\x) (x - \x)^{k + 1}}{(x - \x)^{j + 1}}, \quad \Upsilon_0 = \Gamma_0 - \theta_0(\x) \Phi_0.\]
Noting that $\partial \Gamma_0(x) = \theta_0(x) \partial \Phi_0(x)$ for $x$ near $\x$, we have $\mathcal{D}_0(\x) = \dots = \mathcal{D}_{\q - 1}(\x) = 0$ and
\[\mathcal{D}_j(\x) = \frac{1}{(j + 1)!} \sum_{k = \q}^{j}\binom{j}{k} \partial^{k} \theta_0(\x) \partial^{j + 1 - k} \Phi_0(\x), \qquad j = \q,\dots,\left\lfloor \s \right\rfloor.\]
In particular, $\mathcal{D}_\q(\x) = \partial^\q \theta_0(\x) \partial \Phi_0(\x) / (\q + 1)!$ and therefore $\mathcal{M}_\x^\q(v) = \mathcal{D}_\q(\x)v^{\q + 1}$.

\subsubsection{Mean Function Estimation with Known \texorpdfstring{$\q$}{TEXT}}
First, consider the (simpler) case where $\q$ is known. In this case, if
\begin{equation}\label{Consistency condition}
\mathcal{\widetilde{D}}_{\q,n}(\x) \to_\P \mathcal{D}_\q(\x),
\end{equation}
then Assumption \ref{Assumption C} holds when
\begin{equation}\label{eq: simple mean function estimator}
\widetilde{M}_{\x,n}(x) = \mathcal{\widetilde{D}}_{\q,n}(\x) x^{\q + 1}.
\end{equation}

Example-specific estimators $\mathcal{\widetilde{D}}_{\q,n}(\x)$ satisfying the consistency requirement (\ref{Consistency condition}) may be readily available. For instance, in the case of the \citet{Grenander_1956_SAJ} estimator described in Example \ref{[Example] Monotone Density Estimation}, we have $\mathcal{D}_\q(\x) = \partial^\q f_0(\x) / (\q + 1)!$, a consistent estimator of which can be based on any consistent estimator of $\partial^\q f_0(\x)$, such as a standard kernel estimator or, if the evaluation point $\x$ is near the boundary of the support of $X$, boundary adaptive versions thereof.

More generic estimators are also available. For specificity, the remainder of this section focuses on estimators of $\mathcal{D}_{j}(\x)$ obtained by applying numerical derivative-type operators to the following (possibly) non-smooth estimator of $\Upsilon_0$:
\[\widehat{\Upsilon}_n = \widehat{\Gamma}_n - \widehat{\theta}_n(\x) \widehat{\Phi}_n.\]
The fact that $\mathcal{D}_j(\x) = 0$ for $j < \q$ implies that $\mathcal{D}_\q(\x)$ admits the ``monomial approximation'' representation
\[\mathcal{D}_\q(\x) = \lim_{\epsilon \to 0} \left\{\epsilon^{-(\q + 1)}[\Upsilon_0(\x + \epsilon) - \Upsilon_0(\x)]\right\},\]
motivating the estimator
\[\mathcal{\widetilde{D}}_{\q,n}^{\mathtt{MA}}(\x) = \epsilon_n^{-(\q + 1)}[\widehat{\Upsilon}_n(\x + \epsilon_n) - \widehat{\Upsilon}_n(\x)],\]
where $\epsilon_n > 0$ is a tuning parameter. Similarly, the generic ``forward difference'' representation
\[\mathcal{D}_\q(\x) = \lim_{\epsilon \to 0} \left\{\epsilon^{-(\q + 1)} \sum_{k = 1}^{\q + 1}(-1)^{k + \q + 1}\binom{\q + 1}{k}[\Upsilon_0(\x + k \epsilon) - \Upsilon_0(\x)]\right\}\]
motivates the estimator
\[\mathcal{\widetilde{D}}_{\q,n}^{\mathtt{FD}}(\x) = \epsilon_n^{-(\q + 1)} \sum_{k = 1}^{\q + 1} (-1)^{k + \q + 1}\binom{\q + 1}{k}[\widehat{\Upsilon}_n(\x + k \epsilon_n) - \widehat{\Upsilon}_n(\x)].\]

For $\eta > 0$, let
\begin{align*}
\widehat{G}_{\x,n}(v;\eta)  &= \sqrt{n \eta^{-1}} [\widehat{\Gamma}_n(\x + v \eta) - \widehat{\Gamma}_n(\x) - \Gamma_0(\x + v \eta) + \Gamma_0(\x)]\\
&- \theta_0(\x) \sqrt{n \eta^{-1}} [\widehat{\Phi}_n(\x + v \eta) - \widehat{\Phi}_n(\x) - \Phi_0(\x + v \eta) + \Phi_0(\x)],
\end{align*}
and
\[\widehat{R}_{\x,n}(v;\eta) = \eta^{-1}[\widehat{\Phi}_n(\x + v\eta) - \widehat{\Phi}_n(\x) - \Phi_0(\x + v \eta ) +\Phi_0(\x)].\]
Using this notation, the first part of \ref{Assumption B: weak convergence} can be restated as
\[\widehat{G}_{\x,n}(\cdot;a_n^{-1}) \rightsquigarrow \mathcal{G}_\x\]
and the second displayed condition of \ref{Assumption B: uniform convergence Phi} implies that
\[\widehat{R}_{\x,n}(\cdot;a_n^{-1}) = o_\P(1).\]
In the preceding displays, one can typically replace $a_n^{-1} = n^{-1 / (1 + 2 \q)}$ by any $\eta_n > 0$ with $\eta_n = o(1)$ and $a_n^{-1} \eta_n^{-1} = O(1)$. As a consequence, validity of the following assumption usually follows as a by-product of the arguments used to justify \ref{Assumption B: weak convergence} and \ref{Assumption B: uniform convergence Phi}.

\begin{assumption}\label{Assumption D}
For every $\eta_n > 0$ with $\eta_n = o(1)$ and $a_n^{-1} \eta_n^{-1} = O(1)$,
\[\widehat{G}_{\x,n}(1;\eta_n) = O_\P(1) \qquad \text{and} \qquad \widehat{R}_{\x,n}(1;\eta_n) = O_\P(1).\]
\end{assumption}

In turn, Assumption \ref{Assumption D} is useful for the purposes of analyzing $\mathcal{\widetilde{D}}_{j,n}^{\mathtt{MA}}$, and $\mathcal{\widetilde{D}}_{j,n}^{\mathtt{FD}}$.

\begin{lemma}\label{[Lemma] Mean function estimation with known q}
Suppose Assumptions \ref{Assumption A} and \ref{Assumption D} are satisfied, that $r_n (\widehat{\theta}_n(\x) - \theta_0(\x)) = O_\P(1)$, and that $\epsilon_n \to 0$ and $n \epsilon_n^{1 + 2 \q} \to \infty$. If ${\mathcal{\widetilde{D}}_{\q,n}} \in \{\mathcal{\widetilde{D}}_{\q,n}^{\mathtt{MA}} ,\mathcal{\widetilde{D}}_{\q,n}^{\mathtt{FD}}\}$, then
(\ref{Consistency condition}) is satisfied and Assumption \ref{Assumption C} holds with $\widetilde{M}_{\x,n}$ in \eqref{eq: simple mean function estimator}. 
\end{lemma}

\subsubsection{Mean Function Estimation: Adaptation to Unknown Degeneracy}

Next, consider the somewhat more complicated case where $\q$ is unknown. In this case, we are able to obtain a positive result under the (additional) assumption that there are known integers $\bar\q$ and $\underline{\s}$ for which $\q \leq \bar \q \leq \underline{\s} < \s$. As the notation suggests, $\bar\q$ and $\underline{\s}$ are upper and lower bounds on $\q$ and $\s$, respectively. The stated condition furthermore implies that $\s > \q$, a stronger smoothness condition than the condition $\s \geq \q$ of Assumption \ref{Assumption A}.

For any $\bar\q \geq \q$, noting that $\mathcal{D}_{j}(\x) = 0$ for $j < \q$ and that $\mathcal{D}_\q(\x) > 0$ with $\q$ an odd integer, $\mathcal{M}_\x^\q$ can be written as
\[\mathcal{M}_\x^\q(v) = \sum_{\ell = 1}^{\left\lfloor(\bar\q + 1) / 2 \right\rfloor} \1(2 \ell - 1 \leq \q) \mathcal{D}_{2 \ell - 1}(\x)^+ v^{2 \ell}.\]
Dropping the indicator function term from each summand, we obtain
\[\mathcal{\bar{M}}_\x(v) = \sum_{\ell = 1}^{\left\lfloor (\bar\q + 1) / 2 \right\rfloor} \mathcal{D}_{2 \ell - 1}(\x)^+ v^{2 \ell}.\]
The majorant $\mathcal{\bar{M}}_\x$ is an ``adaptive'' approximation to $\mathcal{M}_\x^\q$ in the sense that it does not depend on $\q$, yet satisfies the local approximation property
\[\mathcal{\bar{M}^\q}_\x \rightsquigarrow \mathcal{M}_\x^\q, \qquad \mathcal{\bar{M}^\q}_\x(v) = \sqrt{n a_n} \mathcal{\bar{M}}_\x(v a_n^{-1}).\]
Moreover, the following (``global'') positivity property automatically holds: for some $c > 0$,
\[\inf_{|v| > K^{-1}} \mathcal{\bar{M}}_\x(v) > c K^{-(\q + 1)} \qquad \text{for every } K > 0.\]
As a consequence, it seems plausible that a ``plug-in'' estimator of $\mathcal{\bar{M}}_\x$ would satisfy Assumption \ref{Assumption C} under reasonable conditions. Indeed, if
\begin{equation}\label{Rate of convergence condition}
a_n^{\q - (2 \ell - 1)} (\mathcal{\widetilde{D}}_{2 \ell - 1,n}(\x) - \mathcal{D}_{2 \ell - 1}(\x)) = o_\P(1), \qquad \ell = 1,\dots,\left\lfloor (\bar\q + 1) / 2 \right\rfloor,
\end{equation}
then Assumption \ref{Assumption C} is satisfied by
\begin{equation}\label{eq: adaptive mean function estimator}
\widetilde{M}_{\x,n}(v) = \sum_{\ell = 1}^{\left\lfloor (\bar\q + 1) / 2 \right\rfloor } \mathcal{\widetilde{D}}_{2 \ell - 1,n}(\x)^+ v^{2\ell}.
\end{equation}

For ``small'' $\ell$ (namely, for $2 \ell - 1 < \q$), the precision requirement (\ref{Rate of convergence condition}) is stronger than consistency, but fortunately it turns out that the requirement can be met as long as $\s$ admits an integer $\underline{\s}$ satisfying the condition $\bar\q \leq \underline{\s} < \s$. To describe estimators satisfying \eqref{Rate of convergence condition} under this condition, let the constants $c_1 ,\dots, c_{\underline{\s} + 1}$ be such that invertibility holds on the part of the $(\underline{\s} + 1) \times (\underline{\s}+1)$ matrix with element $(k,p)$ given by $c_{k}^p$ and let the defining property of $\{\lambda_{j}^{\mathtt{BR}}(k):1 \leq 
 k \leq \underline{\s} + 1\}$ be
\[\sum_{k = 1}^{\underline{\s} + 1}\lambda_{j}^{\mathtt{BR}}(k) c_k^{p} = \1(p = j + 1), \qquad p = 1,\dots,\underline{\s} + 1.\]
Then, for any $j = 1,\dots,\underline{\s}$, the ``bias-reduced'' estimator
\[\mathcal{\widetilde{D}}_{j,n}^{\mathtt{BR}}(\x) = \epsilon_n^{-(j + 1)} \sum_{k = 1}^{\underline{\s} + 1} \lambda_{j}^{\mathtt{BR}}(k)[\widehat{\Upsilon}_n(\x + c_k \epsilon_n) - \widehat{\Upsilon}_n(\x)]\]
is motivated by the fact that as $\epsilon\to0$,
\[\mathcal{D}_{j}(\x) = \epsilon^{-(j + 1)}\sum_{k = 1}^{\underline{\s} + 1} \lambda_{j}^{\mathtt{BR}}(k)[\Upsilon_0(\x + c_k \epsilon) - \Upsilon_0(\x)] + O(\epsilon^{\min(\underline{\s} + 1,\s) - j})\]
when $\s > j$. Relative to $\mathcal{\widetilde{D}}_{j,n}^{\mathtt{MA}}(\x)$ and $\mathcal{\widetilde{D}}_{j,n}^{\mathtt{FD}}(\x)$, this is a distinguishing feature of $\mathcal{\widetilde{D}}_{j,n}^{\mathtt{BR}}(\x)$ and as it turns out this feature will enable us to formulate sufficient conditions for (\ref{Rate of convergence condition}).

\begin{lemma}\label{[Lemma] Mean function estimation unknown q}
Suppose Assumptions \ref{Assumption A} and \ref{Assumption D} are satisfied, that $r_n (\widehat{\theta}_n(\x) - \theta_0(\x)) = O_\P(1)$, and that $\epsilon_n \to 0$. If $\q \leq \bar\q \leq \underline{\s} < \s$, then
\[a_n^{\q - j} (\mathcal{\widetilde{D}}_{j,n}^{\mathtt{BR}}(\x) -\mathcal{D}_{j}(\x)) = O(a_n^{\q - j} \epsilon_n^{\min(\underline{\s} + 1,\s) - j}) + O_{\P}[(a_n \epsilon_n)^{-(j + 1/2)} + (a_n \epsilon_n)^{-j}]\]
for $j = 1,\dots,\underline{\s}$, implying in particular that if also
\[n \epsilon_n^{(1 + 2 \bar \q) \min(\underline{\s},\s - 1)/(\bar \q - 1)} \to 0 \quad \text{and} \quad n \epsilon_n^{1 + 2 \bar \q} \to \infty,\]
then (\ref{Rate of convergence condition}) is satisfied and Assumption \ref{Assumption C} holds with $\widetilde{M}_{\x,n}$ in \eqref{eq: adaptive mean function estimator} if $\mathcal{\widetilde{D}}_{2\ell-1,n}=\mathcal{\widetilde{D}}_{2\ell-1,n}^{\mathtt{BR}}$.
\end{lemma}

As alluded to previously, the ability of $\mathcal{\widetilde{D}}^{\mathtt{BR}}$ to satisfy (\ref{Rate of convergence condition}) is attributable in large part to its bias properties. In an attempt to highlight this, the first display of the lemma gives a stochastic expansion wherein the $O(a_n^{\q-j}\epsilon_n^{\min(\underline{\s}+1,\s)-j})$ term is a (possibly) negligible bias term. For $\mathcal{\widetilde{D}}\in\{\mathcal{\widetilde{D}}^{\mathtt{MA}},\mathcal{\widetilde{D}}^{\mathtt{FD}}\}$, the analogous stochastic expansions are of the form
\[a_n^{\q - j} (\mathcal{\widetilde{D}}_{j,n}(\x) - \mathcal{D}_{j}(\x)) = O(a_n^{\q - j}\epsilon_n^{\q - j}) + O_{\P}[(a_n \epsilon_n)^{-(j + 1/2)} + (a_n \epsilon_n)^{-j}],\]
the $O(a_n^{\q - j} \epsilon_n^{\q - j})$ term also being a bias term. When $a_n \epsilon_n \to \infty$ (as is required for the ``noise'' term in the stochastic expansion to be $o_\P(1)$), this bias term is non-negligible and the estimators $\mathcal{\widetilde{D}}^{\mathtt{MA}}$ and $\mathcal{\widetilde{D}}^{\mathtt{FD}}$ therefore do not satisfy (\ref{Rate of convergence condition}).

Under additional assumptions (including $\s \geq \underline{\s} + 1$ and additional smoothness on $\Phi_0$), $\mathcal{\widetilde{D}}^{\mathtt{BR}}$ admits a Nagar-type mean squared error (MSE) expansion that can be used to select $\epsilon_n. $ The resulting approximate MSE formula is of the form
\[\epsilon_n^{2(\underline{\s} + 1 - j)} \mathsf{B}_{j}^{\mathtt{BR}}(\x)^{2} + \frac{1}{n\epsilon_n^{1 + 2 j}} \mathsf{V}_{j}^{\mathtt{BR}}(\x),\]
where the bias constant is
\[\mathsf{B}_{j}^{\mathtt{BR}}(\x) = \mathcal{D}_{\underline{\s} + 1}(\x)\sum_{k = 1}^{\underline{\s} + 1} \lambda_{j}^{\mathtt{BR}}(k) c_k^{\underline{\s} + 2}\]
and the variance constant is
\[\mathsf{V}_{j}^{\mathtt{BR}}(\x) = \sum_{k = 1}^{\underline{\s} + 1}\sum_{l = 1}^{\underline{\s} + 1} \lambda_{j}^{\mathtt{BR}}(k) \lambda_{j}^{\mathtt{BR}}(l) \mathcal{C}_\x(c_k,c_l);\]
for details, see Section SA.1.7 in the Supplemental Appendix \citep{Cattaneo-Jansson-Nagasawa_2024_supp}.
Assuming $\mathsf{B}_{j}^{\mathtt{BR}}(\x)\neq0$, the approximate MSE is minimized by
\[\epsilon_{j,n}^{\mathtt{BR}}(\x) = \left(\frac{1 + 2 j}{2(\underline{\s} + 1 - j)} \frac{\mathsf{V}_{j}^{\mathtt{BR}}(\x)}{\mathsf{B}_{j}^{\mathtt{BR}}(\x)^{2}}\right)^{1/(3 + 2 \underline{\s})} n^{-1/(3 + 2 \underline{\s})},\]
a feasible version of which can be constructed by replacing $\mathcal{D}_{\underline{\s}+1}(\x)$ and $\mathcal{C}_\x$ with estimators in the expressions for $\mathsf{B}_{j}^{\mathtt{BR}}(\x)$ and $\mathsf{V}_{j}^{\mathtt{BR}}(\x)$, respectively.

\subsection{Bootstrapping}\label{[Subsection] Bootstrapping}

Suppose inference is to be based on a random sample $\mathbf{Z}_1,\dots,\mathbf{Z}_n$ from the distribution of some $\mathbf{Z.}$ In all examples of which we are aware, the bootstrap parts of Assumption \ref{Assumption B} are satisfied when $(\widehat{\Gamma}_n^*,\widehat{\Phi}_n^*,\widehat{u}_n^*)$ is given by the nonparametric bootstrap analog of $(\widehat{\Gamma}_n,\widehat{\Phi}_n,\widehat{u}_n)$. Nevertheless, computationally simpler alternatives are often available and in what follows we will present one such alternative. To motivate our proposal, it is instructive to begin by revisiting Example \ref{[Example] Monotone Density Estimation}.

\setcounter{example}{0}
\begin{example}[Monotone Density Estimation, continued]
In this example, $\mathbf{Z}=X$. Moreover, defining $\gamma_0(x;\mathbf{z})=\1(\mathbf{z}\leq x)$ and $\phi_0(x;\mathbf{z})=x$, we have the representations
\[\Gamma_0(x)=\E[\gamma_0(x;\mathbf{Z})]\qquad\text{and} \qquad \Phi_0(x)=\E[\phi_0(x;\mathbf{Z})],\]
and the estimators $\widehat{\Gamma}_n$ and $\widehat{\Phi}_n$ are linear in the sense that they are of the form
\[\widehat{\Gamma}_n(x)=\frac{1}{n}\sum_{i=1}^n\gamma_0(x;\mathbf{Z}_i)\qquad\text{and}\qquad\widehat{\Phi}_n(x)=\frac{1}{n}\sum_{i=1}^n\phi_0(x;\mathbf{Z}_i),\]
respectively. Finally, $\widehat{u}_n = \max_{1\leq i\leq n}\mathbf{Z}_i \lor \x$.

Letting $\mathbf{Z}_{1,n}^*,\dots,\mathbf{Z}_{n,n}^*$ be a random sample from the empirical distribution of $\mathbf{Z}_1,\dots,\mathbf{Z}_n$, $\widehat{u}_n^* = \max_{1 \leq i \leq n}\mathbf{Z}_{i,n}^* \lor \x$ is the nonparametric bootstrap analog of $\widehat{u}_n$, while
\[\widehat{\Gamma}_n^*(x)=\frac{1}{n}\sum_{i = 1}^n\gamma_0(x;\mathbf{Z}_{i,n}^*) \qquad \text{and} \qquad \widehat{\Phi}_n^*(x) = \frac{1}{n} \sum_{i = 1}^n\phi_0(x;\mathbf{Z}_{i,n}^*)\]
are the nonparametric bootstrap analogs of $\widehat{\Gamma}_n$ and $\widehat{\Phi}_n$, respectively.

In the case of $\widehat{u}_n$, the alternative bootstrap analog $\widehat{u}_n^*=\widehat{u}_n$ is computationally trivial and automatically satisfies the bootstrap parts of \ref{Assumption B: uhat} and \ref{Assumption B: Phi uhat}. As for $\widehat{\Gamma}_n$ and $\widehat{\Phi}_n$, their nonparametric bootstrap analogs admit the weighted bootstrap representations
\[\widehat{\Gamma}_n^*(x)=\frac{1}{n}\sum_{i=1}^nW_{i,n}\gamma_0(x;\mathbf{Z}_i)\qquad\text{and}\qquad\widehat{\Phi}_n^*(x)=\frac{1}{n}\sum_{i=1}^nW_{i,n}\phi_0(x;\mathbf{Z}_i),\]
where, conditionally on $\mathbf{Z}_1,\dots,\mathbf{Z}_n$, $(W_{1,n},\dots,W_{n,n})$ is multinomial distributed with probabilities $(n^{-1},\dots,n^{-1})$ and number of trials $n$. (In turn, the weighted bootstrap interpretation of the nonparametric bootstrap analog of $(\widehat{\Gamma}_n,\widehat{\Phi}_n)$ in this example is interesting partly because it can be used to embed the nonparametric bootstrap in a class of bootstraps also containing the Bayesian bootstrap and the wild bootstrap.)
\end{example}
\setcounter{example}{2}

Looking beyond Example \ref{[Example] Monotone Density Estimation}, finding a computationally trivial $\widehat{u}_n^*$ satisfying the bootstrap part of \ref{Assumption B: uhat} and \ref{Assumption B: Phi uhat} is usually straightforward. On the other hand, although a weighted bootstrap interpretation of the nonparametric bootstrap version of the estimators $\widehat{\Gamma}_n$ and $\widehat{\Phi}_n$ is available whenever they are linear (e.g., in Example \ref{[Example] Monotone Regression Estimation}), there is no shortage of examples for which linearity does not hold. Nevertheless, the weighted bootstrap representation of the nonparametric bootstrap in Example \ref{[Example] Monotone Density Estimation} turns out be useful for our purposes, as it is suggestive of computationally attractive alternatives to the nonparametric bootstrap in more complicated examples.

When the non-bootstrap part of \ref{Assumption B: weak convergence} holds, the estimators $\widehat{\Gamma}_n$ and $\widehat{\Phi}_n$ are typically asymptotically linear in the sense that they admit (possibly) unknown functions $\gamma_0$ and $\phi_0$ (satisfying $\Gamma_0(x)=\E[\gamma_0(x;\mathbf{Z})]$ and $\Phi_0(x)=\E[\phi_0(x;\mathbf{Z})]$, respectively) for which the approximations
\[\widehat{\Gamma}_n(x)\approx\bar{\Gamma}_n(x)=\frac{1}{n}\sum_{i=1}^n \gamma_0(x;\mathbf{Z}_i)\qquad\text{and}\qquad\widehat{\Phi}_n(x)\approx\bar{\Phi}_n(x)=\frac{1}{n}\sum_{i=1}^n\phi_0(x;\mathbf{Z}_i)\]
are suitably accurate. Assuming also that $\gamma_0$ and $\phi_0$ admit sufficiently well behaved estimators $\widehat{\gamma}_n$ and $\widehat{\phi}_n$ (say), it then stands to reason that the salient properties of $\widehat{\Gamma}_n$ and $\widehat{\Phi}_n$ are well approximated by those of the easy-to-compute exchangeable bootstrap-type pair
\[\widehat{\Gamma}_n^*(x) = \frac{1}{n} \sum_{i = 1}^nW_{i,n} \widehat{\gamma}_n(x;\mathbf{Z}_i) \qquad \text{and} \qquad \widehat{\Phi}_n^*(x) = \frac{1}{n} \sum_{i = 1}^n W_{i,n} \widehat{\phi}_n(x;\mathbf{Z}_i),\]
where $W_{1,n},\dots,W_{n,n}$ denote exchangeable random variable (independent of $\mathbf{Z}_1,\dots,\mathbf{Z}_n$).

To give a precise statement, let
\[\psi_x(v;\mathbf{z}) = \gamma_0(x + v;\mathbf{z}) - \gamma_0(x;\mathbf{z}) - \theta_0(x) [\phi_0(x + v;\mathbf{z}) -\phi_0(x;\mathbf{z})]\]
and define
\[\check{\Gamma}_n(x) = \frac{1}{n} \sum_{i = 1}^n \widehat{\gamma}_n(x;\mathbf{Z}_i) \qquad \text{and} \qquad  \check{\Phi}_n(x) = \frac{1}{n} \sum_{i = 1}^n \widehat{\phi}_n(x;\mathbf{Z}_i).\]
In addition, for any function class $\mathfrak{F}$, let $N_{U}(\varepsilon,\mathfrak{F})$ denote the associated uniform covering numbers relative to $L_{2};$ that is, for any $\varepsilon>0$, let
\[N_{U}(\varepsilon,\mathfrak{F}) = \sup_{Q} N(\varepsilon \left\| \bar{F}\right\|_{Q,2},\mathfrak{F},L_{2}(Q)),\]
where $\bar{F}$ is the minimal envelope function of $\mathfrak{F}$, $\left\|\cdot\right\|_{Q,2}$ is the $L_{2}(Q)$ norm, $N(\cdot)$ is the covering number, and the supremum is over all discrete probability measure $Q$ with $\left\|\bar{F}\right\|_{Q,2} > 0$.

\begin{assumption}\label{Assumption E}
For some $\delta > 0$, the following are satisfied:
\begin{enumerate}[label=(E\arabic*)]
    \item \label{Assumption E: iid} $\mathbf{Z}_1,\dots,\mathbf{Z}_n,$ are independent and identically distributed.
    
    \item \label{Assumption E: bootstrap weights}$W_{1,n},\dots,W_{n,n}$ are exchangeable random variables independent of $\mathbf{Z}_1,\dots ,\mathbf{Z}_n$, $\widehat{\gamma}_n$, and $\widehat{\phi}_n$.
    In addition, for some $\mathfrak{r}>(4\q+2)/(2\q-1)$,
    \[\frac{1}{n}\sum_{i=1}^nW_{i,n}=1,\quad\frac{1}{n}\sum_{i=1}^n    (W_{i,n}-1)^{2}\to_\P1,\quad\text{and}\quad\E[|W_{1,n}|^{\mathfrak{r}}]=O(1).\]
    
    \item \label{Assumption E: Gamma hat}$\sup_{x\in I}|\widehat{\Gamma}_n(x)-\bar{\Gamma}_n(x)|=o_{\P}(1)$ and
    \[\frac{1}{n}\sum_{i=1}^n\sup_{x\in I}|\widehat{\gamma}_n(x;\mathbf{Z}_i)-\gamma_0(x;\mathbf{Z}_i)|^2=o_{\P}(1).\]
    Also,
    \begin{equation*}
            \frac{a_n}{n}\sum_{i=1}^n\sup_{x \in I_{\x}^{\delta}}\left|\widehat{\gamma}_n(x;\mathbf{Z}_i) - \gamma_0(x;\mathbf{Z}_i) \right|^2 =o_{\P}(1),
        \end{equation*} 
        and there exist $\beta_\gamma \in [1/2,\q+1)$ and random variables $A_{\gamma,n}=o_{\P}(1)$ and $B_{\gamma,n}=o_{\P}(a_n^{\beta_\gamma})$ such that, for $V \in (0,2\delta]$,
        \begin{equation*}
            \sqrt{na_n}\sup_{|v|\leq V } \left|\widehat{\Gamma}_n(\x+v ) -\widehat{\Gamma}_n(\x) - \bar{\Gamma}_n(\x+v ) + \bar{\Gamma}_n(\x)\right|  \leq  A_{\gamma,n}+  V^{\beta_\gamma}B_{\gamma,n} ,
        \end{equation*}
        \begin{equation*}
            \sqrt{na_n}\sup_{|v|\leq V} \left|\check{\Gamma}_n(\x+v) -\check{\Gamma}_n(\x) - \bar{\Gamma}_n(\x+v ) + \bar{\Gamma}_n(\x)\right| \leq A_{\gamma,n}+  V^{\beta_\gamma}B_{\gamma,n}.
        \end{equation*}

        In addition, for some $\rho_\gamma\in (0,2)$,
        \begin{equation*}
            \limsup_{\varepsilon\downarrow 0}\frac{\log N_{U}(\varepsilon,\mathfrak{F}_{\gamma})}{\varepsilon^{-\rho_\gamma}} < \infty,\quad\E[\bar{F}_{\gamma}(\mathbf{Z})^2]<\infty, \quad \limsup_{\varepsilon\downarrow 0} \frac{\log N_{U}(\varepsilon,\widehat{\mathfrak{F}}_{\gamma,n})}{\varepsilon^{-\rho_\gamma}} = O_{\P}(1),
        \end{equation*}
        where $\mathfrak{F}_{\gamma}=\{\gamma_0(x;\cdot):x\in I\}$, $\bar{F}_{\gamma}$ is its minimal envelope, and $\widehat{\mathfrak{F}}_{\gamma,n}=\{\widehat{\gamma}_n(x;\cdot):x\in I\}$.
        
        Also,
        \begin{equation*}
            \limsup_{\eta\downarrow 0} \frac{\E[\bar{D}_{\gamma}^{\eta}(\mathbf{Z})^2+\bar{D}_{\gamma}^{\eta}(\mathbf{Z})^4]}{\eta} <\infty,
        \end{equation*}
        where $\bar{D}_{\gamma}^{\eta}$ is the minimal envelope of $\{\gamma_0(x;\cdot)-\gamma_0(\x;\cdot):x\in I_\x^\eta\}$.

    \item \label{Assumption E: Phi hat}
    $\sup_{x\in I}|\widehat{\Phi}_n(x)-\bar{\Phi}_n(x)|=o_{\P}(1)$ and
    \[\frac{1}{n}\sum_{i=1}^n\sup_{x\in I}|\widehat{\phi}_n(x;\mathbf{Z}_i)-\phi_0(x;\mathbf{Z}_i)|^2=o_{\P}(1).\] 
    Also, $a_n |\widehat{\Phi}_n(\x)-\bar{\Phi}_n(\x)|=o_{\P}(1),a_n |\check{\Phi}_n(\x)-\bar{\Phi}_n(\x)|=o_{\P}(1)$,
    \begin{equation*}
        \frac{a_n}{n}\sum_{i=1}^n\sup_{x \in I_{\x}^{\delta}}\left|\widehat{\phi}_n(x;\mathbf{Z}_i) -\phi_0(x;\mathbf{Z}_i) \right|^2 =o_{\P}(1),
    \end{equation*}
    and there exist $\beta_\phi \in [1/2,\q]$ and random variables $A_{\phi,n}=o_{\P}(1)$ and $B_{\phi,n}=o_{\P}(a_n^{\beta_\phi})$ such that, for $V \in (0,2\delta]$,
    \begin{equation*}
        \sqrt{na_n} \sup_{|v|\leq V}\left|\widehat{\Phi}_n(\x+v) -\widehat{\Phi}_n(\x) - \bar{\Phi}_n(\x+v) + \bar{\Phi}_n(\x)\right|\leq  A_{\phi,n}+  V^{\beta_\phi}B_{\phi,n} ,
        \end{equation*}
    \begin{equation*}
        \sqrt{na_n} \sup_{|v|\leq V}\left|\check{\Phi}_n(\x+v) -\check{\Phi}_n(\x) - \bar{\Phi}_n(\x+v) + \bar{\Phi}_n(\x)\right|\leq  A_{\phi,n}+  V^{\beta_\phi}B_{\phi,n}.
    \end{equation*}
     
        In addition, for some $\rho_\phi\in (0,2)$,
        \begin{equation*}
            \limsup_{\varepsilon\downarrow 0}\frac{\log N_{U}(\varepsilon,\mathfrak{F}_{\phi})}{\varepsilon^{-\rho_\phi}} < \infty,\quad\E[\bar{F}_{\phi}(\mathbf{Z})^2]<\infty, \quad \limsup_{\varepsilon\downarrow 0} \frac{\log N_{U}(\eta,\widehat{\mathfrak{F}}_{\phi,n})}{\varepsilon^{-\rho_\phi}} = O_{\P}(1),
        \end{equation*}
        where $\mathfrak{F}_{\phi}=\{\phi_0(x;\cdot):x\in I\}$, $\bar{F}_{\phi}$ is its minimal envelope, and $\widehat{\mathfrak{F}}_{\phi,n}=\{\widehat{\phi}_n(x;\cdot):x\in I\}$.
        
        Also,
        \begin{equation*}
            \limsup_{\eta\downarrow 0} \frac{\E[\bar{D}_{\phi}^{\eta}(\mathbf{Z})^2+\bar{D}_{\phi}^{\eta}(\mathbf{Z})^4]}{\eta} <\infty,
        \end{equation*}
        where $\bar{D}_{\phi}^{\eta}$ is the minimal envelope of $\{\phi_0(x;\cdot)-\phi_0(\x;\cdot):x\in I_\x^\eta\}$.
    
    \item \label{Assumption E: covariance kernel} For every $\eta_n>0$ with $a_n\eta_n=O(1)$,
    \[\sup_{v,v'\in[-\eta_n,\eta_n]}\frac{\E [|\psi_\x(v;\mathbf{Z})-\psi_\x(v';\mathbf{Z})|]}{|v-v'|}=O(1)\]
    and, for all $v,u\in\R$, and for some $\mathcal{C}_\x$ (as defined in \ref{Assumption B: weak convergence}),
    \[\frac{\E[\psi_\x(\eta_n v;\mathbf{Z})\psi_{\x}(\eta_n u;\mathbf{Z})]}{\eta_n}\to\mathcal{C}_{\x}(v,u).\]
\end{enumerate}

\end{assumption}

\begin{lemma}
\label{[Lemma] Bootstrapping}Suppose Assumptions \ref{Assumption A} and \ref{Assumption E} are satisfied. Then \ref{Assumption B: weak convergence}-\ref{Assumption B: uniform convergence Phi} are satisfied. 
\end{lemma}

\section{Examples}\label{[Section] Examples}

We apply our main results to two distinct sets of examples, both previously analyzed in \citet{Westling-Carone_2020_AoS} and \citet{Westling-Gilbert-Carone_2020_JRRSB}, and references therein. In the supplemental appendix, we also consider two other set of examples, namely monotone hazard estimation \citep{Huang-Wellner_1995_SJS} and monotone distribution estimation \citep{vanderVaart-vanderLaan_2006_IJB}. All the proofs are found in the supplemental appendix. 

For future reference, we collect the conditions on the bootstrap weights in an assumption.
\begin{assumption-BW}\label{Assumption BW}
    The bootstrap weights $W_{1,n},\dots,W_{n,n}$ are exchangeable , independent of the random sample $\mathbf{Z}_1,\dots,\mathbf{Z}_n$, and satisfy, for some $\mathfrak{r}>6$,
    \begin{equation*}
        \frac{1}{n}\sum_{i=1}^nW_{i,n}=1,\qquad\frac{1}{n}\sum_{i=1}^n(W_{i,n}-1)^{2}\to_\P1,\qquad\text{and}\qquad\E[|W_{1,n}|^{\mathfrak{r}}]=O(1).
    \end{equation*}
\end{assumption-BW}

\subsection{Monotone Density Estimation}\label{[Subsection] Examples: Monotone Density Estimation}

Consider the problem of estimating the density of a non-negative, continuously distributed random variable, possibly with censoring. Let $\mathbf{Z}_1,\dots,\mathbf{Z}_n$ be $i.i.d.$ copies of $\mathbf{Z}=(\check{X},\Delta)'$, with $\check{X} = X \land C$ and $\Delta=\1(X\leq C)$.
Assuming that $f_0$, the density of $X$, (exists and) is non-decreasing on $I=[0,u_0]$ where $u_0$ is some point in the support of $X$, the parameter of interest is $\theta_0(\x)$ for some $\x\in(0,u_0)$, where $\theta_0=f_0$. This setup generalizes Example \ref{[Example] Monotone Density Estimation}. We take $\Phi_0(x)=x$ and therefore have $\Gamma_0=F_0$, where $F_0$ is the cdf of $X$.

The canonical case of no censoring (i.e.,\ $\P[C \geq X]=1$) has been considered in Example \ref{[Example] Monotone Density Estimation}. Assumptions under which this example is covered by our general theory is given in the following result.
\begin{corollary}\label{[Corollary] monotone density without censoring}
    Let $u_0$ be the supremum of the support of $X$. Suppose $\x$ is in the interior of $I$, $\theta_0$ satisfies \ref{Assumption A: theta0}, and Assumption BW holds. Then, Assumptions \ref{Assumption A}, \ref{Assumption B}, and \ref{Assumption D} hold with 
    \begin{equation*}
        \widehat{\Gamma}_n(x) = \frac{1}{n}\sum_{i=1}^n \1(X_i\leq x),\quad \widehat{\Gamma}_n^*(x) = \frac{1}{n}\sum_{i=1}^n W_{i,n}\1(X_i\leq x),
    \end{equation*}
    \begin{equation*}
        \widehat{\Phi}_n(x) = \widehat{\Phi}_n^*(x) = x, \quad \widehat{u}_n = \widehat{u}_n^* = \max_{1\leq i\leq n}X_i \lor \x,
    \end{equation*}
\begin{equation*}
    \mathcal{C}_\x(s,t) = f_0(\x) (|s| \land |t|) \1(\sign(s)=\sign(t)), \quad \mathcal{D}_\q(\x) = \frac{\partial^\q f_0(\x)}{(\q+1)!}.
\end{equation*}
\end{corollary}

Under the assumptions of Corollary \ref{[Corollary] monotone density without censoring}, it follows from Theorem \ref{[Theorem] Main result} that the bootstrap consistency result \eqref{Bootstrap consistency} holds for any $\widetilde{M}_{\x,n}$ satisfying Assumption \ref{Assumption C}. In particular, any pointwise consistent estimator of $\partial^\q f_0(\x)$ could be used to estimate $\mathcal{D}_\q(\x)$. Alternatively, since Assumption \ref{Assumption D} also holds, $\widetilde{M}_{\x,n}$ in \eqref{eq: simple mean function estimator} with $\widetilde{\mathcal{D}}_{\q,n} \in \{ \widetilde{\mathcal{D}}_{\q,n}^{\mathtt{MA}},\widetilde{\mathcal{D}}_{\q,n}^{\mathtt{FD}},\widetilde{\mathcal{D}}_{\q,n}^{\mathtt{BR}}\}$ can also be used, provided that $\epsilon_n \to 0$ and $n\epsilon_n^{1+2\q}\to \infty$.

Suppose now that censoring occurs completely at random; that is, suppose $X\indep C$ \citep[e.g.,][and references therein]{Huang-Wellner_1995_SJS}.
Let $\widehat{S}_n$ and $\widehat{G}_n$ be the Kaplan-Meier estimators of the survival functions $S_0(\cdot)=\P[X>\cdot]$ and $G_0(\cdot)=\P[C>\cdot]$, and define $\widehat{F}_n=1-\widehat{S}_n$ and $\widehat{\Lambda}_n(x) = \int_0^x \widehat{S}_n(u)^{-1}d\widehat{F}_n(u)$. In anticipation of the following result, we note that $\Gamma_0(x)=F_0(x)=\E[\gamma_0(x;\mathbf{Z})]$, where
\[\gamma_0(x;\mathbf{Z})= F_0(x)+S_0(x)\left[\frac{\Delta\1(\check{X}\leq x) }{S_0(\check{X})G_0(\check{X})} - \int_0^{\check{X} \land x}\frac{\Lambda_0(du)}{S_0(u)G_0(u)}\right],\]
with $\Lambda_0(x) = \int_0^x S_0(u)^{-1}dF_0(u)$.
\begin{corollary}\label{[Corollary] monotone density independent right-censoring}
    Let $u_0$ be an interior point in the support of $X$. Suppose $\x$ is in the interior of $I$, $\theta_0$ satisfies \ref{Assumption A: theta0}, $C\indep X$, $C$ is absolutely continuous on $I$ with bounded density, and $S_0(u_0)G_0(u_0) > 0$. Then, Assumptions \ref{Assumption A}, \ref{Assumption B}, and \ref{Assumption D} hold with 
    \begin{equation*}
        \widehat{\Gamma}_n(x)= \widehat{F}_n(x),\quad \widehat{\Gamma}_n^*(x)= \frac{1}{n}\sum_{i=1}^n W_{i,n}\widehat{\gamma}_n(x;\mathbf{Z}_i)
    \end{equation*}
    \begin{equation*}
        \widehat{\gamma}_n(x;\mathbf{Z}) = \widehat{F}_n(x) + \widehat{S}_n(x)\left[\frac{\Delta\1(\check{X}\leq x)}{\widehat{S}_n(\check{X})\widehat{G}_n(\check{X})} - \int_0^{\check{X} \land x} \frac{d\widehat{\Lambda}_n(u) }{\widehat{S}_n(u)\widehat{G}_n(u)} \right],
\end{equation*}
    \begin{equation*}
        \widehat{\Phi}_n(x)=\widehat{\Phi}_n^*(x)=x,\quad \widehat{u}_n=\widehat{u}_n^*=u_0,
    \end{equation*}
\begin{equation*}
    \mathcal{C}_\x(s,t) = \frac{f_0(\x)}{G_0(\x)}(|s| \land |t|) \1(\sign(s)=\sign(t)), \quad \mathcal{D}_\q(\x) = \frac{\partial^\q f_0(\x)}{(\q+1)!}.
\end{equation*}
\end{corollary}
\noindent Theorem \ref{[Theorem] Main result} implies that inference based on our proposed bootstrap-assisted method is asymptotically valid for this example as well. Additionally, Assumption \ref{Assumption D} also holds under the above conditions, and we can use one of the numerical derivative estimators discussed in Section \ref{[Subsection] Mean function estimation}. 

In the supplemental appendix, we also consider estimation of monotone density under conditionally independent right-censoring, where covariates are present. We do not discuss this example here because it involves additional notation and technicalities.

\subsection{Monotone Regression Estimation}\label{[Subsection] Examples: Monotone Regression Estimation}

Consider now the problem of regression estimation, possibly with covariate adjustment. We assume that $\mathbf{Z}_1,\dots,\mathbf{Z}_n$ are $i.i.d.$ copies of $\mathbf{Z}=(Y,X,\mathbf{A}')'$, where $X$ is continuous on its support $I$ with density $f_0$ and cdf $F_0$. Defining $\mu_0(x,\mathbf{a})=\E[Y|X=x,\mathbf{A}=\mathbf{a}]$, the parameter of interest is $\theta_0(\x)$ for some $\x$ in the interior of $I$, where $\theta_0(\cdot)=\E[\mu_0(\cdot,\mathbf{A})]$ is assumed to be non-decreasing on $I$. (If there are no covariates $\mathbf{A}$, then $\theta_0(\cdot)=\mu_0(\cdot)=\E[Y|X=\cdot]$.) We take $\Phi_0=F_0$ and set $u_0=1$. 

The classical monotone regression estimator has been considered in Example \ref{[Example] Monotone Regression Estimation}. Assumptions under which this example is covered by our general theory is given in the following result. To state it, let $\varepsilon=Y-\E[Y|X]$ and define $\sigma_0^2(x) = \E[\varepsilon^2|X=x]$.
\begin{corollary}\label{[Corollary] classical monotone regression}
    Suppose that $\x$ is in the interior of $I$, $\theta_0$ satisfies \ref{Assumption A: theta0}, $\Phi_0$ satisfies \ref{Assumption A: Phi0}, and Assumption BW holds. If $\sigma_0^2(x)$ is continuous and positive at $\x$, $\E[\varepsilon^2]<\infty$, and if $\sup_{x \in I_\x^\delta}\E[\varepsilon^4|X=x]<\infty$ for some $\delta >0$, then Assumptions \ref{Assumption A}, \ref{Assumption B}, and \ref{Assumption D} hold with 
    \begin{equation*}
        \widehat{\Gamma}_n(x)= \frac{1}{n}\sum_{i=1}^n Y_i\1(X_i\leq x),\quad \widehat{\Gamma}_n^*(x)= \frac{1}{n}\sum_{i=1}^n W_{i,n}Y_i\1(X_i\leq x),
    \end{equation*}
    \begin{equation*}
        \widehat{\Phi}_n(x)=\frac{1}{n}\sum_{i=1}^n \1(X_i\leq x),\quad \widehat{\Phi}_n^*(x)=\frac{1}{n}\sum_{i=1}^n W_{i,n}\1(X_i\leq x),\quad \widehat{u}_n=\widehat{u}_n^*=1,
    \end{equation*}
\begin{equation*}
    \mathcal{C}_\x(s,t) = f_0(\x) \sigma_0^2(\x) (|s| \land |t|) \1(\sign(s)=\sign(t)),\qquad \mathcal{D}_\q(\x) = \frac{f_0(\x)\partial^\q\mu_0(\x)}{(\q+1)!}.
\end{equation*}
\end{corollary}%
\noindent Under the assumptions of Corollary \ref{[Corollary] classical monotone regression}, it follows from Theorem \ref{[Theorem] Main result} that the bootstrap consistency result \eqref{Bootstrap consistency} holds for any $\widetilde{M}_{\x,n}$ satisfying Assumption \ref{Assumption C}. In this example, any pointwise consistent estimators of $f_0(\x)$ and $\partial^\q \mu_0(\x)$ could be used to estimate $\mathcal{D}_\q(\x)$. Alternatively, instead of using two distinct estimators, we can use the numerical derivative-type estimators since Assumption \ref{Assumption D} also holds.

Next, consider the case of monotone regression estimation with covariate-adjustment.
Introducing covariates $\mathbf{A}$, let $f_{X|\mathbf{A}}(\cdot|\mathbf{a})$ denote the conditional Lebesgue density of $X$ given $\mathbf{A}=\mathbf{a}$ and define $g_0(x,\mathbf{a})=f_{X|\mathbf{A}}(x|\mathbf{a})/f_0(x)$. Noting that $\Gamma_0(x)=\E[\gamma_0(x;\mathbf{Z})]$, where
\[\gamma_0(x;\mathbf{Z})= \1(X\leq x) \left[\frac{Y-\mu_0(X,\mathbf{A})}{g_0(X,\mathbf{A})} - \theta_0(X)\right],\]
we assume that $\widehat{\Gamma}_n$ is of the form $\widehat{\Gamma}_n(x)= n^{-1}\sum_{i=1}^n \widehat{\gamma}_n(x;\mathbf{Z}_i)$, where
\[\widehat{\gamma}_n(x;\mathbf{Z})= \1(X\leq x)\left[\frac{Y-\widehat{\mu}_n(X,\mathbf{A})}{\widehat{g}_n(X,\mathbf{A})} + \frac{1}{n}\sum_{j=1}^n \widehat{\mu}_n(X,\mathbf{A}_j)\right],\]
with $\widehat{\mu}_n$ and $\widehat{g}_n$ being preliminary estimators of $\mu_0$ and $g_0$ based on $\mathbf{Z}_1,\dots,\mathbf{Z}_n$. As discussed by \cite{Westling-Gilbert-Carone_2020_JRRSB}, there are a number of candidates for $\widehat{\mu}_n$ and $\widehat{g}_n$. To accommodate a relatively wide class of estimators, the following assumption employs high-level conditions on $\widehat{\mu}_n$ and $\widehat{g}_n$. To state the assumption, let $\varepsilon=Y-\E[Y|X,\mathbf{A}]$ and define $\sigma_0^2(x,\mathbf{a})=\E[\varepsilon^2|X=x,\mathbf{A}=\mathbf{a}]$.

\begin{assumption-MRC}\label{Assumption MRC}
    For some $\delta>0$, the following are satisfied:
    \begin{enumerate}[label=\normalfont(\roman*),noitemsep,itemindent=*]
        \item[(i)] $\E[g_0(\x,\mathbf{A})^{-1}\sigma_0^2(\x,\mathbf{A})]>0$, $\E[\varepsilon^2]<\infty$, and $\sup_{x \in I_\x^\delta}\E[\varepsilon^4|X=x]<\infty$.
        
        \item[(ii)] $f_{X|\mathbf{A}}$ is bounded and $g_0$ is bounded away from zero.

        \item[(iii)] 
        There exist random variables $A_n=o_{\P}(1)$ and $B_n=O_\P(a_n^{1/2})$ such that
        \[\sqrt{na_n}\sup_{\vert v\vert\leq  V} \vert\widehat{\Gamma}_n(\x+v)-\widehat{\Gamma}_n(\x) - \bar{\Gamma}_n(\x+v)+\bar{\Gamma}_n(\x)\vert\leq A_n + V B_n, \qquad V\in (0,2\delta].\]
        In addition,
        \[\frac{a_n}{n}\sum_{i=1}^n|\widehat{\mu}_n(X_i,\mathbf{A}_i)-\mu_0(X_i,\mathbf{A}_i)|^2=o_{\P}(1),\]
        \[\frac{a_n}{n^2}\sum_{i=1}^n\sum_{j=1}^n|\widehat{\mu}_n(X_i,\mathbf{A}_j)-\mu_0(X_i,\mathbf{A}_j)|^2=o_{\P}(1),\]
        and
        \[\frac{a_n}{n}\sum_{i=1}^n\varepsilon_i^2|\widehat{g}_n(X_i,\mathbf{A}_i)-g_0(X_i,\mathbf{A}_i)|^2=o_{\P}(1).\]
        
        \item[(iv)] $\E[\bar{\mu}(\mathbf{A})^2]<\infty$, where
        \[\bar{\mu}(\mathbf{A})=\sup_{|x-x'|\leq\delta}\frac{|\mu_0(x,\mathbf{A})-\mu_0(x',\mathbf{A})|}{|x-x'|}.\]
        
        \item[(v)] $\E[\bar{\sigma}^2(\mathbf{A})]<\infty$, where, for some function $\omega$ with $\lim_{\eta\downarrow0}\omega(\eta)=0$,
        \[\left|\frac{\sigma_0^2(x,\mathbf{A})f_0(x)}{g_0(x,\mathbf{A})} -\frac{\sigma_0^2(\x,\mathbf{A})f_0(\x)}{g_0(\x,\mathbf{A})}\right| \leq\omega(|x-\x|) \bar{\sigma}^2(\mathbf{A}), \qquad x \in I_\x^\delta.\]
    \end{enumerate} 
\end{assumption-MRC}
Part (iii) of the assumption contains high-level conditions. In the supplemental appendix, we provide primitive sufficient conditions using nonparametric estimators based on sample splitting.

\begin{corollary}\label{[Corollary] monotone regression with covariates}
Suppose that $\x$ is in the interior of $I$, $\theta_0$ satisfies \ref{Assumption A: theta0}, $\Phi_0$ satisfies \ref{Assumption A: Phi0}, and Assumption BW holds. If Assumption MRC holds, then Assumptions \ref{Assumption A} and \ref{Assumption B} hold with 
    \begin{equation*}
        \widehat{\Gamma}_n^*(x)= \frac{1}{n}\sum_{i=1}^n W_{i,n}\widehat{\gamma}_n(x;\mathbf{Z}_i),
    \end{equation*}
    \begin{equation*}
        \widehat{\Phi}_n(x)=\frac{1}{n}\sum_{i=1}^n \1(X_i\leq x),\quad \widehat{\Phi}_n^*(x)=\frac{1}{n}\sum_{i=1}^n W_{i,n}\1(X_i\leq x),\quad \widehat{u}_n=\widehat{u}_n^*=1,
    \end{equation*}
    \begin{equation*}
    \mathcal{C}_\x(s,t) = f_0(\x) \E\left[\frac{\sigma_0^2(\x,\mathbf{A})}{g_0(\x,\mathbf{A})}\right] (|s| \land |t|) \1(\sign(s)=\sign(t)),\quad \mathcal{D}_\q(\x) = \frac{f_0(\x)\partial^\q\theta_0(\x)}{(\q+1)!}.
\end{equation*}
\end{corollary}
\noindent As a consequence, if $\widetilde{M}_{\x,n}$ satisfies Assumption \ref{Assumption C}, then the ``reshaped'' bootstrap estimator $\widetilde{\theta}_n^*(\x)=\partial_-\GCM_{[0,1]}(\widetilde{\Gamma}_n^*\circ\widehat{\Phi}_n^{*-})\circ\widehat{\Phi}_n^*(\x)$ gives a bootstrap-assisted distributional approximation satisfying \eqref{Bootstrap consistency}. In this example, the covariance kernel contains $\E[g_0(\x,\mathbf{A})^{-1}\sigma_0^2(\x,\mathbf{A})]$ and the plug-in procedure discussed in Section \ref{[Section] Bootstrap-based Distributional Approximation} would involve estimating this nuisance parameter, which in turn requires estimation of $\sigma_0^2(\x,\cdot)$. When $\mathbf{A}$ is of moderate dimension, this would introduce a challenging estimation problem on top of estimating $\mu_0$ and $g_0$. Our bootstrap-assisted procedure circumvents estimation of this additional nuisance parameter.

\section{Simulations}\label{[Section] Simulations}

We consider the canonical case of monotone regression (i.e., Example \ref{[Example] Monotone Regression Estimation}). We estimate the non-decreasing regression function $\theta_0(\cdot)=\E[Y|X=\cdot]$ at an interior point $\x$ using a random sample of observations. Three distinct data generating processes (DGPs) are considered. For each DGP, we let $X$ be a uniform $(0,1)$ random variable and we let $Y=\theta_0(X)+\sigma_0(X)\widetilde{\varepsilon}$, where $\widetilde{\varepsilon}$ is a standard normal random variable independent of $X$. DGP 1 sets $\theta_0(x)=2\exp(x-0.5)$ and $\sigma_0=1$, DGP 2 sets $\theta_0(x)=2\exp(x-0.5)$ and $\sigma_0=\exp(x)$, and DGP 3 sets $\theta_0(x)=24\exp(X-0.5) -24(X-0.5) - 12(X-0.5)^2$ and $\sigma_0=0.1$. In all three cases, we set $\x=0.5$. The second DGP exhibits heteroskedastic regression errors, and the third DGP features a regression function whose first derivative equals zero at the evaluation point $\x$; that is, DGP 3 exhibits a degree of degeneracy and if inference is conducted without knowledge of the fact that $\q=3$, then this feature makes the inference problem more challenging. 

The Monte Carlo experiment employs a sample size $n=1,000$ with $2,000$ bootstrap replications and $4,000$ simulations, and compares four types of bootstrap-based inference procedures: the standard non-parametric bootstrap, $m$-out-of-$n$ bootstrap, the ``plug-in method'' which simulates the limit law with estimated nuisance parameters, and our proposed bootstrap-assisted inference method implemented using the bias-reduced numerical derivative estimator discussed in Section \ref{[Subsection] Mean function estimation}. For the plug-in method, we estimate the covariance kernel and the nuisance parameter $\partial \Phi_0(\x)$ using local linear kernel estimators, and we estimate the mean function by the same numerical derivative estimator as for our bootstrap-assisted procedure. For $m$-out-of-$n$ bootstrap and the plug-in method, researchers need to specify the characteristic exponent $\q$, and we use the true (infeasible) value. For the numerical derivative estimator, we developed a rule of thumb for the step size $\epsilon_n$ to operationalize the procedure. For further details, see Section SA.6 of \cite{Cattaneo-Jansson-Nagasawa_2024_supp}. For our proposed method, we report results for three implementations: (i) (infeasible) procedure using the true value of $\mathcal{D}_0(\x)$ (referred to as ``oracle''), (ii) implementation using the numerical derivative estimator with a correct specification of $\q$ (referred to as ``known $\q$''), and (iii) ``robust'' implementation only assuming $\q\in\{1, 3\}$ (referred to as ``robust''). 

Table \ref{TABLE: isoreg} presents the numerical results. We report empirical coverage for nominal $95\%$ confidence intervals and their average interval length. For all of the DGPs considered, our proposed bootstrap-assisted inference method leads to confidence intervals with excellent empirical coverage and average interval length. The infeasible ``oracle'' procedure attains empirical coverage very close to the nominal $95\%$, which aligns with our theoretical results, and feasible procedures using the numerical derivative estimators perform almost identical to the infeasible oracle version. 

In DGPs 1 and 2, our procedures outperform both the standard non-parametric bootstrap (which is inconsistent) and the $m$-out-of-$n$ bootstrap (which is consistent) in empirical coverage and average length. In DGP 3, the $m$-out-of-$n$ with the subsample size $m=\lceil n^{1/2}\rceil$ performs comparable to our procedure, but two caveats should be noted. First, the $m$-out-of-$n$ bootstrap performance is sensitive to the choice of the subsample size. Therefore, to operationalize the $m$-out-of-$n$ bootstrap procedure, one needs to develop a reliable procedure to choose the subsample size. Another caveat, arguably more important in this context, is that the $m$-out-of-$n$ bootstrap procedure requires the knowledge of the convergence rate of the estimator. In our simulation, we use the true convergence rate for $m$-out-of-$n$ bootstrap, but in practice, one needs to assume or estimate from data the convergence rate. Since the convergence rate and the limit distribution of the generalized Grenander-type estimator crucially hinge on unknown $\q$, this feature of the $m$-out-of-$n$ bootstrap may be unappealing.
In contrast, our proposed bootstrap-based procedure only requires specifying an upper bound on $\q$ and it automatically adapts to the unknown convergence rate.

The plug-in method performs on par with our method, which is not unexpected because (i) the covariance kernel in this setting takes a relatively simple form and estimating it is not challenging with sample size $n=1,000$ and (ii) we use the true $\q$ to simulate the limit law. In other examples where covariates are present, estimating covariance kernels typically involves nonparametric estimation of additional, possibly high-dimensional, nuisance functions. Our proposed method reduces the number of nuisance parameters that need to be estimated. Moreover, the good performance of the plug-in method for DGP 3 crucially depends on correctly specifying $\q$, whereas our ``robust'' procedure only requires specifying an upper bound on $\q$, and this is another desirable feature of our proposed method.

\begin{table}\renewcommand{\arraystretch}{1.2}
	\caption{Simulations, Monotone Regression Estimator, 95\% Confidence Intervals.}
	{\resizebox{\columnwidth}{!}{ \label{TABLE: isoreg}
\begin{tabular}{lrrrrcrrrrcrrrr}
\hline\hline
\multicolumn{1}{l}{\bfseries }&\multicolumn{4}{c}{\bfseries DGP 1}&\multicolumn{1}{c}{\bfseries }&\multicolumn{4}{c}{\bfseries DGP 2}&\multicolumn{1}{c}{\bfseries }&\multicolumn{4}{c}{\bfseries DGP 3}\tabularnewline
\cline{2-5} \cline{7-10} \cline{12-15}
\multicolumn{1}{l}{}&\multicolumn{1}{c}{$\tilde{\mathcal{D}}_{1,n}$}&\multicolumn{1}{c}{$\tilde{\mathcal{D}}_{3,n}$}&\multicolumn{1}{c}{Coverage}&\multicolumn{1}{c}{Length}&\multicolumn{1}{c}{}&\multicolumn{1}{c}{$\tilde{\mathcal{D}}_{1,n}$}&\multicolumn{1}{c}{$\tilde{\mathcal{D}}_{3,n}$}&\multicolumn{1}{c}{Coverage}&\multicolumn{1}{c}{Length}&\multicolumn{1}{c}{}&\multicolumn{1}{c}{$\tilde{\mathcal{D}}_{1,n}$}&\multicolumn{1}{c}{$\tilde{\mathcal{D}}_{3,n}$}&\multicolumn{1}{c}{Coverage}&\multicolumn{1}{c}{Length}\tabularnewline
\hline
{\bfseries Standard}&&&&&&&&&&&&&&\tabularnewline
~~&$$&$$&$0.828$&$0.373$&&$$&$$&$0.838$&$0.519$&&$$&$$&$0.912$&$0.029$\tabularnewline
\hline
{\bfseries m-out-of-n}&&&&&&&&&&&&&&\tabularnewline
~~$m = \lceil n^{1/2} \rceil$&$$&$$&$0.896$&$0.413$&&$$&$$&$0.915$&$0.583$&&$$&$$&$0.939$&$0.031$\tabularnewline
~~$m = \lceil n^{2/3} \rceil$&$$&$$&$0.868$&$0.400$&&$$&$$&$0.884$&$0.557$&&$$&$$&$0.928$&$0.029$\tabularnewline
~~$m = \lceil n^{4/5} \rceil$&$$&$$&$0.851$&$0.393$&&$$&$$&$0.864$&$0.545$&&$$&$$&$0.918$&$0.029$\tabularnewline
\hline
{\bfseries Plug-in}&&&&&&&&&&&&&&\tabularnewline
~~&$1.053$&$0.000$&$0.950$&$0.399$&&$1.052$&$0.000$&$0.946$&$0.553$&&$0.000$&$1.003$&$0.936$&$0.028$\tabularnewline
\hline
{\bfseries Reshaped}&&&&&&&&&&&&&&\tabularnewline
~~Oracle&$1.000$&$0.000$&$0.941$&$0.395$&&$1.000$&$0.000$&$0.951$&$0.550$&&$0.000$&$1.000$&$0.946$&$0.029$\tabularnewline
~~ND known q&$1.053$&$0.000$&$0.949$&$0.398$&&$1.052$&$0.000$&$0.945$&$0.547$&&$0.000$&$1.003$&$0.936$&$0.028$\tabularnewline
~~ND robust&$1.053$&$0.599$&$0.950$&$0.401$&&$1.052$&$0.938$&$0.954$&$0.559$&&$0.014$&$1.003$&$0.960$&$0.030$\tabularnewline
\hline
\end{tabular}
 } }
	\raggedright Notes:
	(i) Panel \textbf{Standard} refers to standard non-parametric bootstrap, Panel \textbf{m-out-of-n} refers to $m$-out-of-$n$ non-parametric bootstrap with subsample $m$, Panel \textbf{Plug-in} refers to the plug-in method, Panel \textbf{Reshaped} refers to our proposed bootstrap-assisted procedure. \newline
	(ii) Columns ``$\tilde{\mathcal{D}}_{1,n}$'' and ``$\tilde{\mathcal{D}}_{3,n}$'' report the averages of the estimated $\mathcal{D}_{1},\mathcal{D}_3$ across simulations, and Columns ``Coverage'' and ``Length'' report empirical coverage and average length of bootstrap-based $95\%$ percentile confidence intervals, respectively.\newline
	(iii) ``Oracle'' corresponds to the infeasible version of our proposed procedure using the true value of $\mathcal{D}_\q$, ``ND known $\q$'' corresponds to our proposed procedure using the bias-reduced numerical derivative estimator with a correct specification of $\q$, and ``ND robust'' corresponds to our proposed procedure only assuming $\q\in \{1,3\}$. The step size choice for the numerical derivative estimator is described in the supplemental appendix.\newline 
	(iv) The sample size is $1,000$, the number of bootstrap iterations is $2,000$, and the number of Monte Carlo simulations is $4,000$.
\end{table}

\begin{appendix}

\section{Technical Results and Omitted Details}\label{[Section] Appendix - Technical Results}

This appendix collects several technical results and details omitted from the main paper to improve exposition. First, in Section \ref{[Subsection] Appendix - Generalized Switch Relation} we present a corrected version of the generalization of the switch relation stated by \citet[Supplement]{Westling-Carone_2020_AoS}. Second, in Section \ref{[Subsection] Appendix - Continuity of Limiting Distribution} we present a lemma that can be used to establish continuity of the cdf of the maximizer of a Gaussian process whose covariance kernel is that of two-sided Brownian motion. Both lemmas are used in the proof of Theorem \ref{[Theorem] Main result} and may be of independent interest as well. Third, in Section \ref{[Subsection] Appendix - Bootstrap Inconsistency} we establish inconsistency of the plug-in nonparametric bootstrap for a large class of Generalized Grenander-type Estimators, following early work by \cite{Kosorok_2008_BookCh} and \cite{Sen-Banerjee-Woodroofe_2010_AoS}. Finally, in Section \ref{[Subsection] Appendix - Omitted Formulas} we report omitted formulas from Section \ref{[Subsection] Heuristics}.  

\subsection{Generalized Switch Relation}\label{[Subsection] Appendix - Generalized Switch Relation}
For a real-valued function $f$ defined on a set $X \subset \R$, $\argmax_{x \in X} \{f(x)\}$ denotes the (possibly empty) set of maximizers of $f$ over $X$. If $\argmax_{x \in X} \{f(x)\}$ is non-empty and contains a largest element, we denote that element by $\max \argmax_{x \in X} \{f(x)\}$.

In their analysis of generalized Grenander-type estimators, \cite{Westling-Carone_2020_AoS} relied on a generalization of the \textit{switch relation} \citep{Groeneboom_1985_BookCh}. Their generalized switch relation is given in Lemma 1 of \citet[Supplement]{Westling-Carone_2020_AoS}. For the purposes of comparing it with Lemma \ref{[Lemma] Generalized Switch Relation} below, it is convenient to restate that lemma as follows:

\begin{description}
\item[Statement GSR \citep{Westling-Carone_2020_AoS}.] \emph{Let $\Phi$ and $\Gamma$ be real-valued functions defined on an interval $I \subseteq \R$ and suppose that $\Phi$ is non-decreasing and right-continuous. Fix $l,u$ in $\Phi(I)$ with $l < u$ and let $\theta = \partial_- \GCM_{[l,u]}(\Gamma \circ \Phi^-) \circ \Phi$. If $I$ is closed, $\Phi(I) \subseteq [l,u]$, and if $\Gamma$ and $\Gamma \circ \Phi^-$ are lower semi-continuous, then}
\begin{equation}\label{Westling and Carone Lemma 1}
\theta(\x) > t \quad \iff \quad \sup \argmax_{x \in \Phi^-([l,u])} \{t \Phi(x) - \Gamma(x)\} < \Phi^-(\Phi(\x))
\end{equation}
\emph{for any $t \in \R$ and any $\x \in I$ with $\Phi(\x) \in (l,u)$}.
\end{description}

The main purpose of the following example is to show that without further restrictions, the $\argmax$ in \eqref{Westling and Carone Lemma 1} can be empty and the relation \eqref{Westling and Carone Lemma 1} can be violated (if we interpret $\sup \argmax_{x \in X} \{f(x)\}$ as $\inf X$ when $\argmax_{x \in X} \{f(x)\}$ is empty). In other words, Statement GSR can fail to hold.

\begin{example}
Let $I = [l,u] = [0,1], \Gamma(x) = \gamma x$ (for some $\gamma \in \R$), and let
\begin{equation*}
\Phi(x)=
\begin{cases}
x & \text{if } 0 \leq x < 1/2\\
1 & \text{if } 1/2 \leq x \leq 1.
\end{cases}
\end{equation*}
Then
\begin{equation*}
t \Phi(x) - \Gamma(x) =
\begin{cases}
(t - \gamma) x & \text{if } 0 \leq x < 1/2\\
t - \gamma/2 & \text{if } 1/2 \leq x \leq 1,
\end{cases}
\end{equation*}
so $\argmax_{x \in \Phi^-([l,u])}\{t \Phi(x) - \Gamma(x)\}$ is empty when $0 > t > \gamma$.
	
In particular, if $0 > t > \gamma$ and if $\x \in (0,1/2)$, then
\[\sup\argmax_{x \in \Phi^-([l,u])}\{t \Phi(x) - \Gamma(x)\} = 0 < \x = \Phi^-(\Phi(\x)),\]
whereas $\theta(\x) = \partial_- \GCM_{[l,u]} (\Gamma \circ \Phi^-) \circ \Phi(\x) = \gamma < t,$ so \eqref{Westling and Carone Lemma 1} is violated.
	
For future reference, we note that $\Phi(I) = \Phi(I) \cap [l,u] = [0,1/2) \cup \{1\}$ in this example.
\end{example}

The problems highlighted by the example are attributable to the fact that $\Phi(I) \cap [l,u]$ is not closed. Fortunately, it turns out that one can obtain a result in the spirit of \eqref{Westling and Carone Lemma 1} as long as $\Phi(I) \cap [l,u]$ is closed. The following lemma gives a precise statement.

\begin{lemma}[Generalized Switch Relation]\label{[Lemma] Generalized Switch Relation} Let $\Phi$ and $\Gamma$ be real-valued functions defined on an interval $I \subseteq \R$ and suppose that $\Phi$ is non-decreasing and right-continuous. Fix $l,u \in \Phi(I)$ with $l < u$ and let $\theta = \partial_- \GCM_{[l,u]} (\Gamma \circ \Phi^-) \circ \Phi$. If $\Phi(I) \cap [l,u]$ is closed, then
\begin{equation}\label{Generalized Switch Relation}
\theta(\x) >t \quad \iff \quad \max \argmax_{x \in \Phi^-([l,u])}\{t\Phi(x) - \LSC_{\Phi}({\Gamma})(x)\} < \Phi^-(\Phi(\x))
\end{equation}
for any $t \in \R$ and any $\x \in I$ with $\Phi(\x) \in (l,u)$, where
\[\LSC_{\Phi} ({\Gamma})(x) = \LSC (\Gamma \circ \Phi^-) \circ \Phi(x),\]
and where $\LSC (\cdot)$ denotes the greatest lower semi-continuous minorant.
\end{lemma}

The assumptions of the lemma seem mild. In particular, the assumption that $\Phi(I) \cap [l,u]$ is closed is satisfied not only in the examples of Section \ref{[Section] Examples}, but in all examples of which we are aware. 

The lemma does not assume lower semi-continuity of $\Gamma$ or $\Gamma \circ \Phi^-$. Instead, the conclusion involves the function $\LSC_{\Phi} ({\Gamma})$. As the notation suggests, $\LSC_{\Phi} (\cdot) = \LSC (\cdot)$ when $\Phi$ is the identity. More generally, $\LSC_{\Phi}(\cdot) = \LSC (\cdot)$ when $\Phi$ is continuous and strictly increasing and whether or not these properties hold it follows from Lemma 4.3 of \cite{vanderVaart-vanderLaan_2006_IJB} that $\LSC_{\Phi}({\Gamma})$ admits the following characterization:
\[\LSC_{\Phi} ({\Gamma})(x) = \liminf_{y \to \Phi(x)} (\Gamma \circ \Phi^-) (y).\]

\subsection{Continuity of argmax}\label{[Subsection] Appendix - Continuity of Limiting Distribution}

Let $\{\G(v) : v \in \R\}$ be a Gaussian process with $\E[\G(v)]=\mu(v)$ and $\Cov[\G(v),\G(u)]=C(|v|\wedge |u|)\1(\sign(v)=\sign(u))$ for some $C>0$, i.e.\ $\G$ is a scalar multiple of two-sided Brownian motion with mean shift. 
Under conditions on $\mu$ stated below, there exists a unique maximizer of $\G(v)$ over $v \in \R$ with probability one \citep{Kim-Pollard_1990_AoS}. We complement this known fact with a result showing that the cdf of $\argmax_{v \in \R} \{\G(v)\}$ is continuous.

\begin{lemma}\label{[Lemma] Continuity of argmax}
Suppose that for each $v \in \R$, $\lim_{\delta \downarrow 0}[\mu(v + \delta) -\mu(v)] / \sqrt{\delta} = 0$ and that $\limsup_{|v| \to \infty} \mu(v) / |v|^{c} = - \infty $ for some $c > 1/2$. Then $x \mapsto \P[\argmax_{v \in \R} \{\G(v)\} \leq x]$ is continuous.
\end{lemma}

Under Assumptions \ref{Assumption A} and \ref{Assumption B} of the main text, $-\mathcal{M}_\x^\q + t \mathcal{L}_\x$ satisfies the hypothesis of Lemma \ref{[Lemma] Continuity of argmax} for any $t \in \R$. It therefore follows from the lemma that the function
\[ x \mapsto \P \left[\argmin_{v \in \R} \{\mathcal{G}_\x(v) + \mathcal{M}_\x^\q(v) - t\mathcal{L}_\x(v)\} \geq x \right]\]
is continuous at $x = 0$ for any $t \in \R$. We utilize that fact in our proof of \eqref{[Asymptotics] thetahat} and note in passing that most of the existing literature on monotone function estimators seems to implicitly utilize a similar continuity result.

\subsection{Bootstrap Inconsistency}\label{[Subsection] Appendix - Bootstrap Inconsistency}

Employing the exchangeable bootstrap setup introduced in Section \ref{[Subsection] Bootstrapping}, we show the inconsistency of the plug-in bootstrap-based distributional approximation. More precisely, we consider the ``na\"ive'' (plug-in) bootstrap estimator
\[\widehat{\theta}_n^* (\x) = \partial_- \GCM_{[0,\widehat{u}_n^*]} (\widehat{\Gamma}_n^* \circ \widehat{\Phi}_n^{*-}) \circ \widehat{\Phi}_n^* (\x)\]
with
\[\widehat{\Gamma}_n^*(x) = \frac{1}{n} \sum_{i = 1}^n W_{i,n} \widehat{\gamma}_n(x;\mathbf{Z}_i), \qquad \widehat{\Phi}_n^*(x) = \frac{1}{n}\sum_{i=1}^n W_{i,n}\widehat{\phi}_n(x;\mathbf{Z}_i),\]
and establish the following negative result.

\begin{theorem}\label{[Theorem]: Bootstrap inconsistency}
Suppose Assumptions \ref{Assumption A}, \ref{Assumption B: uhat}-\ref{Assumption B: jump of Phi hat}, and \ref{Assumption E} hold. Then
\[\sup_{t \in \R}\left| \P_n^* \left[\widehat{\theta}_n^*(\x) - \widehat{\theta}_n(\x) \leq t\right] - \P \left[\widehat{\theta}_n(\x) - \theta_0(\x) \leq t\right] \right| \not = o_\P(1).\]
\end{theorem}

Theorem \ref{[Theorem]: Bootstrap inconsistency} accommodates exchangeable bootstrap schemes and a wide class of generalized Grenander-type estimators and contains as a special case the the well-known fact that the nonparametric bootstrap approximation to the distribution of the classical Grenander estimator is inconsistent \cite[e.g.,][]{Kosorok_2008_BookCh,Sen-Banerjee-Woodroofe_2010_AoS}.

\subsection{Omitted Formulas}\label{[Subsection] Appendix - Omitted Formulas}

The objects $\widehat{V}_{\x,n}^\q,M_{\x,n}^\q,\widehat{L}_{\x,n}^\q$, and $\widehat{Z}_{\x,n}^\q$ appearing in Section \ref{[Subsection] Heuristics} are defined as follows:
\[\widehat{V}_{\x,n}^\q = \{a_n (x-\x): x \in \widehat{\Phi}_n^-([0,\widehat{u}^*_n])\},\]
\[M_{\x,n}^\q(v) = \sqrt{n a_n} [\Gamma_0(\x + v a_n^{-1}) - \Gamma_0(\x)]-\theta_0(\x)\sqrt{na_n}[\Phi_0(\x+va_n^{-1})-\Phi_0(\x)],\]
and
\[\widehat{L}_{\x,n}^\q(v) = a_n [\widehat{\Phi}_n(\x + v a_n^{-1}) -\widehat{\Phi}_n(\x)], \qquad \widehat{Z}_{\x,n}^\q = a_n [(\widehat{\Phi}_n^- \circ \widehat{\Phi}_n)(\x) - \x].\]

Similarly, the bootstrap analogs of $\widehat{V}_{\x,n}^\q,M_{\x,n}^\q,\widehat{L}_{\x,n}^\q$, and $\widehat{Z}_{\x,n}^\q$ are given by
\[\widehat{V}_{\x,n}^{\q,*} = \{a_n (x-\x): x \in \widehat{\Phi}_n^{*-}([0,\widehat{u}^*_n])\},\]
\[\widehat{M}_{\x,n}^\q(v) = \sqrt{n a_n} [\widehat{\Gamma}_n(\x + v a_n^{-1}) - \widehat{\Gamma}_n(\x)] - \widehat{\theta}_n(\x) \sqrt{n a_n} [\widehat{\Phi}_n(\x + v a_n^{-1}) - \widehat{\Phi}_n(\x)],\]
and
\[\widehat{L}_{\x,n}^{\q,*}(v) = a_n [\widehat{\Phi}_n^*(\x + v a_n^{-1}) - \widehat{\Phi}_n^*(\x)], \qquad \widehat{Z}_{\x,n}^{\q,*} = a_n [(\widehat{\Phi}_n^{*-} \circ \widehat{\Phi}_n^*)(\x) - \x],\]
respectively.
\end{appendix}

\begin{acks}[Acknowledgments]
We are grateful to Boris Hanin, Jason Klusowski, Whitney Newey, and seminar participants at various institutions for their comments. We specially thank Marco Carone and Ted Westling for insightful discussions. The co-Editor, an Associate Editor, and two reviewers offered detailed suggestions that improved our paper.
\end{acks}
\begin{funding}
Cattaneo gratefully acknowledges financial support from the National Science Foundation through grants SES-1947805 and DMS-2210561, and Jansson gratefully acknowledges financial support from the National Science Foundation through grant SES-1947662.
\end{funding}

\begin{supplement}
\stitle{Supplement to ``Bootstrap-Assisted Infernce for Generalized Grenander-Type Estimators.''}
\sdescription{The supplementary material contains proofs of Theorems \ref{[Theorem] Main result} and \ref{[Theorem]: Bootstrap inconsistency}, Corollaries \ref{[Corollary] monotone density without censoring}-\ref{[Corollary] monotone regression with covariates}, and Lemmas \ref{[Lemma] Mean function estimation with known q}-\ref{[Lemma] Bootstrapping} and \ref{[Lemma] Generalized Switch Relation}-\ref{[Lemma] Continuity of argmax}. It discusses additional examples that our main results cover.}
\end{supplement}


\bibliographystyle{imsart-nameyear.bst} 
\bibliography{CJN_2024_AOS--bib}


\end{document}


\begin{frontmatter}

\title{Supplement to ``Bootstrap-Assisted Inference for Generalized Grenander-type Estimators''
}
\runtitle{Bootstrap-Assisted Inference for Monotone Estimators}

\begin{aug}
\author[A]{\fnms{Matias D.}~\snm{Cattaneo}\ead[label=e1]{cattaneo@princeton.edu}},
\author[B]{\fnms{Michael}~\snm{Jansson}\ead[label=e2]{mjansson@econ.berkeley.edu}}
\and
\author[C]{\fnms{Kenichi}~\snm{Nagasawa}\ead[label=e3]{kenichi.nagasawa@warwick.ac.uk}}
\address[A]{Department of Operations Research and Financial Engineering,
Princeton University\printead[presep={,\ }]{e1}}

\address[B]{Department of Economics,
University of California at Berkeley\printead[presep={,\ }]{e2}}

\address[C]{Department of Economics,
University of Warwick\printead[presep={,\ }]{e3}}
\end{aug}

\begin{keyword}[class=MSC]
\kwd[Primary ]{62G09}
\kwd{62G20}
\kwd[; secondary ]{62G07}
\kwd{62G08}
\end{keyword}

\begin{keyword}
\kwd{Monotone estimation}
\kwd{bootstrapping}
\kwd{robust inference}
\end{keyword}

\end{frontmatter}	


\tableofcontents

\clearpage


\section{Proofs}

\subsection{Proof of Lemma \ref{M-[Lemma] Generalized Switch Relation}}

Lemma 4.2 of \cite{vanderVaart-vanderLaan_2006_IJB} implies that, on $(l,u)$,
\begin{equation*}
    \GCM_{[l,u]} (\Gamma\circ\Phi^{-}) = \GCM_{[l,u]} (\LSC (\Gamma\circ\Phi^{-})).
\end{equation*}
Since $\Phi(\x) \in (l,u)$, we therefore have
\begin{equation*}
    \theta(\x) = \partial_{-} \GCM_{[l,u]} (\Gamma\circ\Phi^{-}) \circ \Phi(\x) = \partial_{-} \GCM_{[l,u]} (\LSC (\Gamma\circ\Phi^{-})) \circ \Phi(\x),
\end{equation*}
so Lemma 4.1 of \cite{vanderVaart-vanderLaan_2006_IJB} implies that
\begin{equation*}
    \theta(\x) > t \quad \iff \quad y^{\star} < \Phi(\x),
\end{equation*}
where
\begin{equation*}
    y^{\star} = \max\argmax_{y \in [l,u]} \left\{ ty -\LSC(\Gamma\circ\Phi^{-})(y) \right\}.
\end{equation*}

Suppose $y^{\star} \in \Phi(I)$. Then $y^{\star} = \Phi(x^{\star})$, where
\begin{equation*}
    x^{\star} = \Phi^{-}(y^{\star}) \in \argmax_{x \in \Phi^{-}([l,u])} \left\{ t\Phi(x) - \LSC_{\Phi}(\Gamma)(x) \right\}.
\end{equation*}
In particular,
\begin{equation*}
    \theta(\x) > t \quad \iff \quad y^{\star} = \Phi(x^{\star}) < \Phi(\x) \quad \iff \quad x^{\star} < \Phi^{-}(\Phi(\x)),
\end{equation*}
where, in fact,
\begin{equation*}
    x^{\star} = \max\argmax_{x \in \Phi^{-}([l,u])} \left\{ t\Phi(x) -\LSC_{\Phi}(\Gamma)(x) \right\},
\end{equation*}
because if
\begin{equation*}
    x^{\star} < x' \in \argmax_{x \in \Phi^{-}([l,u])} \left\{ t\Phi(x) -\LSC_{\Phi}(\Gamma)(x) \right\},
\end{equation*}
then, contradicting the definition of $y^{\star}$, we have
\begin{equation*}
    y' = \Phi(x') \in \argmax_{y\in [l,u]} \left\{ty -\LSC(\Gamma\circ\Phi^{-})(y) \right\} \qquad \text { and } \qquad y' > y^{\star}.
\end{equation*}
The proof can therefore be completed by showing that $y^{\star}\in \Phi(I)$.

If $\Phi(I) \supseteq [l,u]$, then there is nothing to show, so suppose $\Phi(I) \not\supseteq [l,u]$. If $y \in \Phi(I)^c \cap [l,u]$, then, since $\Phi(I) \cap [l,u]$ is closed and $l,u \in \Phi(I)$, we have $[y-\eta,y+\eta] \cap \Phi(I) = \emptyset$ for some $\eta > 0$ with $[y-\eta,y+\eta] \subset [l,u]$. Therefore, the function $\LSC (\Gamma\circ\Phi^{-})$ is constant on the interval $[y-\eta/2,y+\eta/2]$, implying in particular that
\begin{equation*}
    y^{\star}  = \max \argmax_{y' \in [l,u]} \left\{ty' -\LSC (\Gamma\circ\Phi^{-})(y') \right\} \neq y. \qed
\end{equation*}

\subsection{Proof of Lemma \ref{M-[Lemma] Continuity of argmax}}

$\G(v)$ is a Gaussian process with $\E[\G(v)]=\mu(v)$ for all $v\in \R$ and $\Cov[\G(v),\G(v')] = C \min\{|v|,|v'|\}\1(\sign(v)=\sign(v')) =: \mathcal{K}(v,v')$ for some positive constant $C$ and for all $v,v' \in \R$. The form of the covariance kernel $\mathcal{K}$ implies that for every $\tau> 0$ and every $v,v' \in \R$, 
\begin{equation*}
    \mathcal{K}(v \tau,v' \tau) = \mathcal{K}(v,v') \tau
\end{equation*}
and
\[\mathcal{K}(v + v',v + v') - 2 \mathcal{K}(v + v',v') + \mathcal{K}(v',v') = \mathcal{K}(v,v).\]
In addition, $\mathcal{K}(1,1) > 0$ and $\lim_{v \downarrow 0}\mathcal{K}(1,v) / \sqrt{v} = 0$.\footnote{As a side note, if the covariance kernel satisfies the two displays, then the Gaussian process is two-sided Brownian motion, up to mean shift and scaling by a positive constant.}

We begin by adapting the arguments of \citet[pp. 196-198]{Kim-Pollard_1990_AoS} to show that a maximizer of $\G(v)$ over $v \in \R$ exists and is unique with probability one. Let $\G^{\mu}(v) = \G(v) - \mu(v)$ be the centered process and suppose that, for the same $c>1/2$ as in the hypothesis,
\begin{equation}\label{Appendix eq: BC lemma}
    \P\left[\limsup_{|v|\to\infty}\frac{\G^{\mu}(v)}{|v|^{c}} > \eta \right] = 0 \qquad \text{for any } \eta > 0.
\end{equation}
Then, with probability one, $\G(v) \to -\infty$ as $|v| \to \infty$, implying in turn that a maximizer of $\G(v)$ exists (because sample paths are continuous). Also, since
\begin{equation*}
    \V[\G(v) - \G(v')] = \mathcal{K}(v,v) + \mathcal{K}(v',v')-2\mathcal{K}(v,v-v') = \mathcal{K}(v-v',v-v') > 0
\end{equation*}
for $v \neq v'$, Lemma 2.6 of \cite{Kim-Pollard_1990_AoS} implies that this maximizer is unique with probability one. In turn, \eqref{Appendix eq: BC lemma} follows from the Borel-Cantelli lemma because
\begin{align*}
    \sum_{k = 2}^{\infty} \P \left[ \sup_{k - 1 \leq |v| \leq k} \frac{\G^{\mu}(v)}{|v|^c} > \eta \right] 
    &\leq \sum_{k = 2}^{\infty} \P \left[ \sup_{|v| \leq k} \G^{\mu}(v) > (k / 2)^c \eta \right] \\
    &= \sum_{k = 2}^{\infty} \P \left[ \sup_{|v| \leq 1} \G^{\mu}(v) > k^{c - 1/2} 2^{-c} \eta \right] \\
    &\leq \frac{2^{4c/(2c-1)} \E \left[ \sup_{|v| \leq 1} |\G^{\mu}(v)|^{4c/(2c-1)} \right]} {\eta^{4c/(2c-1)}} \sum_{k = 2}^{\infty}k^{-2} < \infty,
\end{align*}
where the equality uses the rescaling property $\mathcal{K}(v\tau,v'\tau) = \tau\mathcal{K}(v,v')$ and where the last inequality uses \citet[][Corollary 4.7]{Jain-Marcus_1978_BookCh}.

To show continuity of the function $x \mapsto 
 \P[\argmax_{v \in \R}\{\G(v)\} \leq x]$, it suffices to show that $\P[\argmax_{v \in \R}\{\G(v)\} = x] = 0$ for every $x \in \R$.
Fix $x \in \R$ and define $\widetilde{\G}_x(0) = 0$ and
\begin{equation*}
    \widetilde{\G}_x(v) = \frac{\G(x+v) - \G(x) }{\sqrt{\mathcal{K}(v,v)}}, \qquad v \neq 0.
\end{equation*}
Then $\max_{v \in \R}\widetilde{\G}_x(v) \geq 0$ and, for any set $\mathcal{V} \subset \R$,
\begin{equation}\label{Appendix eq: probability bound with subset}
    \P \left[ \argmax_{v \in \R}\{\G(v)\} = x \right] = \P \left[ \max_{v \in \R} \widetilde{\G}_x(v) \leq 0 \right] \leq \P \left[ \max_{v\in\mathcal{V}} \widetilde{\G}_x(v) \leq 0\right].
\end{equation}
In the sequel, we show that the majorant in \eqref{Appendix eq: probability bound with subset} can be made arbitrarily small by choice of $\mathcal{V}$.
In particular, for $\varepsilon \in (0,1)$, we construct $v_{\varepsilon} > 0$ such that
\begin{equation}\label{Appendix eq: mean_bound}
    \E \left[ \widetilde{\G}_x (v_{\varepsilon}^i) \right] \geq - \varepsilon \quad \text{ for every } i \geq 1,
\end{equation}
and
\begin{equation}\label{Appendix eq: cov_bound}
    \left| \Cov \left( \widetilde{\G}_x(v_{\varepsilon}^i),\widetilde{\G}_x(v_{\varepsilon}^j) \right) \right| \leq \varepsilon \quad \text{ for every } j > i \geq 1.
\end{equation}
For any $N \in \N$, defining $\mathcal{V}_{\varepsilon,N} = \{v_{\varepsilon},\dots,v_{\varepsilon}^N\}$, we therefore have
\begin{equation*}
    \P \left[ \max_{v \in \R} \widetilde{\G}_x(v) \leq 0 \right] \leq \liminf_{\varepsilon \downarrow 0} \P \left[ \max_{v \in \mathcal{V}_{\varepsilon,N}} \widetilde{\G}_x(v) \leq 0 \right] \leq \left| \int_{-\infty}^{0} \frac{\exp(-x^2/2)}{\sqrt{2\pi}}dx \right| ^N = 2^{-N},
\end{equation*}
where the second inequality uses the fact that convergence of means and covariances of normal random vectors implies convergence in distribution. Since $N$ is arbitrary, the left-hand side in the preceding display is zero. Letting $\varepsilon \in (0,1)$ be given, the proof can therefore be completed by exhibiting $v_{\varepsilon}$ satisfying \eqref{Appendix eq: mean_bound} and \eqref{Appendix eq: cov_bound}.

Because $\mathcal{K}(\tau,\tau) = \tau \mathcal{K}(1,1)$ and $\lim_{\tau \downarrow 0}[\mu(x + \tau) - \mu(x)]/\sqrt{\tau} = 0$, there exists $\bar{\tau}_\varepsilon^{'} \in (0,1)$ such that
\begin{equation*}
    \E \left[ \widetilde{\G}_x(\tau) \right] = \frac{\mu(x + \tau) - \mu(x)}{\sqrt{\mathcal{K}(\tau,\tau)}}  =  \frac{[\mu(x + \tau ) -\mu(x)] / \sqrt{\tau}}{\sqrt{\mathcal{K}(1,1)}} \geq -\varepsilon
\end{equation*}
for every $\tau \in (0,\bar{\tau}_\varepsilon^{'}]$. Also, because $\mathcal{K}(\tau_i,\tau_j) = \tau_i\mathcal{K}(1,\tau_j/\tau_i)$ and $\lim_{\tau \downarrow 0}\mathcal{K}(1,\tau )/\sqrt{\tau} = 0$, there exists $\bar{\tau}_\varepsilon^{''} \in (0,1)$ such that
\begin{equation*}
    \Cov\left(\widetilde{\G}_x(\tau_i),\widetilde{\G}_x( \tau_j)\right) = \frac{\mathcal{K}(\tau_i,\tau_j)}{\sqrt{\mathcal{K}(\tau_i,\tau_i)\mathcal{K}(\tau_j,\tau_j)}} = \frac{\mathcal{K}(1,\tau_j/\tau_i)/\sqrt{\tau_j/\tau_i} }{\mathcal{K}(1,1)} \in [-\varepsilon,\varepsilon]
\end{equation*}
for all $\tau_i,\tau_j > 0$ with $\tau_j/\tau_i \leq \bar{\tau}_\varepsilon^{''}$. If $v_{\varepsilon} \in (0, \bar{\tau}_\varepsilon^{'} \land \bar{\tau}_\varepsilon^{''})$, then \eqref{Appendix eq: mean_bound} and \eqref{Appendix eq: cov_bound} are satisfied.\qed

\subsection{Technical Lemmas}

In preparation for the proof of Theorem \ref{M-[Theorem] Main result}, this section presents six technical lemmas. The first lemma is a switching lemma, which will be used when characterizing the limiting distributions obtained in Theorem \ref{M-[Theorem] Main result}.

\begin{lemma}\label{[Appendix Lemma] Switching lemma for R}
Let $\Gamma: \R \to \R$ be a lower semi-continuous function that is bounded from below and satisfies $\lim_{|v| \to \infty}\Gamma(v) / |v| = \infty$. Then, for any $\x,t \in \R$,
\begin{equation*}
    \partial_{-} \GCM_{\R}(\Gamma)(\x) > t \qquad \iff \qquad \max \argmax_{v \in \R} \left\{ vt - \Gamma(v) \right\} < \x.
\end{equation*}
\end{lemma}

\begin{proof}
Because $\lim_{|v| \to \infty}\Gamma(v) / |v| = \infty$, there exists a $K > |\x|$ such that if $|v| \geq K$, then $vt - \Gamma(v) < -\Gamma(0)$, implying in particular that
\begin{equation*}
    \argmax_{v \in \R} \left\{ vt - \Gamma(v) \right\} = \argmax_{v \in [-c,c]} \left\{vt - \Gamma(v) \right\}
\end{equation*}
for every $c \geq K$. Also, by Lemma A.1.\ of \cite{Sen-Banerjee-Woodroofe_2010_AoS} there exists a $c > K$ such that $\GCM_{\R}(\Gamma) = \GCM_{[-c,c]}(\Gamma)$ on $[-K,K]$, implying in particular that
\begin{equation*}
    \partial_{-}\GCM_{\R}(\Gamma)(\x) = \partial_{-}\GCM_{[-c,c]}(\Gamma)(\x).
\end{equation*}
For any such $c$, the conclusion of the lemma is equivalent to the statement
\begin{equation*}
    \partial_{-}\GCM_{[-c,c]}(\Gamma)(\x) > t \qquad \iff \qquad \max \argmax_{v \in [-c,c]}\left\{vt - \Gamma(v)\right\} < \x,
\end{equation*}
whose validity follows from Lemma 4.1 of \cite{vanderVaart-vanderLaan_2006_IJB}.
\end{proof}

The proof of Theorem \ref{M-[Theorem] Main result} furthermore employs various approximations to functionals of the form $\LSC_\Phi(f)$. The approximations in question are obtained using Lemmas \ref{[Appendix Lemma] LSC of Gamma_0 (global)}, \ref{[Appendix Lemma] LSC of Phi (global)}, \ref{[Appendix Lemma] LSC of f (local)}, and \ref{[Appendix Lemma] LSC of Phi (local)}. In all cases, the approximations are based on the representation
\begin{equation}\label{Appendix eq: LSC representation}
    \LSC_\Phi(f)(x) = \liminf_{\epsilon \downarrow 0} \inf_{x' \in \mathcal{X}_\Phi^\epsilon(x)}f(x'),
\end{equation}
where $\mathcal{X}_\Phi^\epsilon(x) = \left(\Phi^-(\Phi(x) - \epsilon),\Phi^-(\Phi(x))\right] \cup \left(\Phi^-(\Phi(x)+),\Phi^-(\Phi(x) + \epsilon) \right)$.

The following lemma uses (\ref{Appendix eq: LSC representation}) and the special structure of $\Gamma_0$ to obtain a simple ``global'' bound on the error of the approximation $\LSC_{\Phi}(\Gamma_0) \approx \Gamma_0$.

\begin{lemma}\label{[Appendix Lemma] LSC of Gamma_0 (global)}
Suppose Assumption \ref{M-Assumption A} holds and suppose $\Phi$ is non-decreasing and right-continuous on $I$. Then, for every $x \in I$,
\begin{equation*}
    |\LSC_{\Phi} (\Gamma_0)(x) - \Gamma_0(x)| \leq 2 \left( \sup_{x' \in I} |\theta_0 (x')| \right) \sup_{x' \in I} |\Phi(x') - \Phi_0 (x')|.
\end{equation*}
\end{lemma}

\begin{proof}
By (\ref{Appendix eq: LSC representation}) and continuity of $\Gamma_0$,
\begin{equation*}
    \LSC_{\Phi} (\Gamma_0)(x) = \min \left[ \Gamma_0(\Phi^-(\Phi(x))),\Gamma_0(\Phi^-(\Phi(x)+)) \right],
\end{equation*}
while, by Assumption \ref{M-Assumption A},
\begin{equation*}
    |\Gamma_0(x')-\Gamma_0(x)| \leq \left( \sup_{x'' \in I} |\theta_0 (x'')| \right) |\Phi_0(x') - \Phi_0(x)|.
\end{equation*}
Now, using $\Phi\circ\Phi^-\circ\Phi = \Phi$,
\begin{align*}
    |\Phi_0(\Phi^-(\Phi(x))) - \Phi_0(x)| &= |\Phi_0(\Phi^-(\Phi(x))) - \Phi(\Phi^-(\Phi(x))) + \Phi(x) - \Phi_0(x)| \\
    &\leq 2 \sup_{x' \in I}|\Phi(x') - \Phi_0 (x')|.    
\end{align*}
Also,
\begin{equation*}
    |\Phi_0(\Phi^-(\Phi(x)+)) - \Phi_0(x)| \leq 2 \sup_{x' \in I} |\Phi(x') - \Phi_0 (x')|
\end{equation*}
because, for every $\eta > 0$,
\begin{align*}
    0
    &\leq \Phi_0(\Phi^-(\Phi(x)+)) - \Phi_0(x) \leq \Phi_0(\Phi^-(\Phi(x) + \eta)) - \Phi_0(x) \\
    &\leq \Phi_0\left( \Phi^- \left( \Phi_0(x) + \sup_{x' \in I} |\Phi(x') - \Phi_0 (x')| + \eta \right) \right) - \Phi_0(x) \\
    &\leq \Phi_0 \left( \Phi_0^- \left(\Phi_0(x)+2\sup_{x' \in I} |\Phi(x') - \Phi_0 (x')| + \eta \right) \right) - \Phi_0(x) \\
    &\leq 2 \sup_{x' \in I} |\Phi(x') - \Phi_0 (x')| + \eta,
\end{align*}
where the last inequality uses continuity of $\Phi_0$.
\end{proof}

A simple ``global'' bound on the error of the approximation $\LSC_{\Phi}(f) \approx f$ is available also in the important special case where $f$ is proportional to $\Phi$.

\begin{lemma}\label{[Appendix Lemma] LSC of Phi (global)}
Suppose Assumption \ref{M-Assumption A} holds and suppose $\Phi$ is non-decreasing and right-continuous on $I$. Then, for every $x \in I$ and every $\theta \in \R$,
\begin{equation*}
    |\LSC_{\Phi}(\theta \Phi)(x) - \theta \Phi(x)| \leq |\theta| \sup_{x' \in I} |\Phi(x') - \Phi(x'-)|.
\end{equation*}
\end{lemma}

\begin{proof}
First, if $\theta > 0$, then the result follows from the fact that, by (\ref{Appendix eq: LSC representation}),
\begin{equation*}
    \LSC_{\Phi}(\theta \Phi)(x) = \theta \Phi(\Phi^-(\Phi(x))-) \leq \theta \Phi(\Phi^-(\Phi(x))) = \theta \Phi(x),
\end{equation*}
where
\begin{equation*}
    \Phi(x) = \Phi(\Phi^-(\Phi(x))) \leq \Phi(\Phi^-(\Phi(x))-) + \sup_{x' \in I} |\Phi(x') - \Phi(x'-)|.
\end{equation*}

Next, if $\theta < 0$, then the result follows from the fact that, by (\ref{Appendix eq: LSC representation}),
\begin{equation*}
    \LSC_{\Phi}(\theta \Phi)(x) = \theta \Phi(\Phi^-(\Phi(x)+)) \leq \theta \Phi(\Phi^-(\Phi(x))) = \theta \Phi(x),
\end{equation*}
where
\begin{align*}
    \Phi(\Phi^-(\Phi(x)+)) &\leq \liminf_{\eta \downarrow 0} \Phi(\Phi^-(\Phi(x) + \eta)-) + \sup_{x' \in I}|\Phi(x') - \Phi(x'-)| \\
    &\leq \Phi(x)+ \sup_{x' \in I}|\Phi(x') - \Phi(x'-)|.
\end{align*}
\end{proof}

Next, we give a ``local'' approximation to $\Phi^- \circ \Phi$. That approximation will later be used in combination with (\ref{Appendix eq: LSC representation}) to obtain ``local'' approximations to $\LSC_\Phi(f)$, but the approximation is also useful in its own right and we therefore state it as a separate lemma.
\begin{lemma}\label{[Appendix Lemma] Phi- circ Phi}
Suppose Assumption \ref{M-Assumption A} holds and suppose $\Phi$ is non-decreasing and right-continuous on $I$. Also, suppose
\begin{equation*}
    \frac{2 \sup_{x' \in I_{\x}^{\delta + \epsilon}} |\Phi(x') - \Phi_0(x')|}{\inf_{x' \in I_{\x}^{\delta + \epsilon}}\partial\Phi_0(x')} < \epsilon
\end{equation*}
for some $\delta,\epsilon > 0$ with $I_{\x}^{\delta + \epsilon} \subseteq I$. Then, for every $x \in I_{\x}^{\delta}$,
\begin{equation*}
    |\Phi^-(\Phi(x)) - x| < \epsilon \qquad \text{and} \qquad |\Phi^-(\Phi(x)+) - x| < \epsilon.
\end{equation*}
\end{lemma}

\begin{proof}
First, suppose $|\Phi^-(\Phi(x))-x| \geq \epsilon$. Then $\Phi^-(\Phi(x)) \leq x - \epsilon$, implying in particular that $\Phi(x-\epsilon) = \Phi(x)$ and therefore also
\begin{equation*}
    \left[\Phi(x - \epsilon) - \Phi_0(x - \epsilon)\right] - \left[\Phi(x) - \Phi_0(x)\right] = \Phi_0(x) - \Phi_0(x - \epsilon).
\end{equation*}
Now, if $x \in I_{\x}^{\delta}$, then $[x-\epsilon,x] \subseteq I_{\x}^{\delta + \epsilon}$, so
\begin{equation*}
    \Phi_0(x) - \Phi_0(x - \epsilon) \geq \epsilon \inf_{x' \in I_{\x}^{\delta + \epsilon}} \partial \Phi_0(x') > 2 \sup_{x' \in I_{\x}^{\delta + \epsilon}}|\Phi(x') - \Phi_0(x')|,
\end{equation*}
whereas
\begin{equation*}
    \left| \left[ \Phi(x - \epsilon) - \Phi_0(x - \epsilon) \right] - \left[ \Phi(x) - \Phi_0(x) \right] \right| \leq 2 \sup_{x' \in I_{\x}^{\delta + \epsilon}} |\Phi(x') - \Phi_0(x')|.
\end{equation*}
In other words, $x \not \in I_{\x}^{\delta}$.

Next, suppose $|\Phi^-(\Phi(x)+) - x| \geq \epsilon$. Then, for every $\eta,\eta' > 0$, $\Phi^-(\Phi(x) + \eta) \geq x + \epsilon$ and therefore
\begin{equation*}
    \Phi(x + \epsilon - \eta')-\sup_{x' \in I_{\x}^{\delta + \epsilon}}|\Phi(x') - \Phi_0(x')| - \eta < \Phi(x)-\sup_{x' \in I_{\x}^{\delta + \epsilon}}|\Phi(x') - \Phi_0(x')|,
\end{equation*}
where, for $x \in I_{\x}^{\delta}$,
\begin{equation*}
    \Phi(x) - \sup_{x' \in I_{\x}^{\delta + \epsilon}}|\Phi(x') - \Phi_0(x')| \leq \Phi_0(x),
\end{equation*}
whereas
\begin{align*}
    &\liminf_{\eta,\eta' \downarrow 0} \left[ \Phi(x + \epsilon - \eta') - \sup_{x' \in I_{\x}^{\delta + \epsilon}} |\Phi(x') - \Phi_0(x')| - \eta  \right] \\
    &\geq \liminf_{\eta' \downarrow 0} \left[ \Phi_0(x + \epsilon - \eta') - 2 \sup_{x' \in I_{\x}^{\delta + \epsilon}} |\Phi(x') - \Phi_0(x')| \right] \\
    &\geq \Phi_0(x) + \left( \inf_{x' \in I_{\x}^{\delta + \epsilon}} \partial \Phi_0(x') \right) \liminf_{\eta' \downarrow 0} \left[ \epsilon - \frac{2 \sup_{x' \in I_{\x}^{\delta + \epsilon}} |\Phi(x') - \Phi_0(x')|} {\inf_{x' \in I_{\x}^{\delta + \epsilon}} \partial \Phi_0 x')} - \eta' \right] > \Phi_0(x).
\end{align*}
In other words, $x \not \in I_{\x}^{\delta}$.
\end{proof}

Next, we obtain two ``local'' approximations to $\LSC_\Phi(f)$. The first of these is a generic approximation obtained by simply combining (\ref{Appendix eq: LSC representation}) and Lemma \ref{[Appendix Lemma] Phi- circ Phi}, but for later reference we state the result as a separate lemma.

\begin{lemma}\label{[Appendix Lemma] LSC of f (local)}
Suppose the assumptions of Lemma \ref{[Appendix Lemma] Phi- circ Phi} hold. Then, for every $x \in I_{\x}^{\delta}$ and every $f : I \to \R$,
\begin{equation*}
    \inf_{|x' - x| \leq \epsilon}f(x') \leq \LSC_{\Phi}(f)(x) \leq \sup_{|x' - x| \leq \epsilon}f(x')
\end{equation*}
and
\begin{equation*}
    |\LSC_{\Phi}(f)(x) - f(x)| \leq \sup_{|x' - x| \leq \epsilon} |f(x') - f(x)|.
\end{equation*}
\end{lemma}

The final lemma is concerned with the special case where $f$ is proportional to $\Phi$. In that case, the following ``local'' analog of Lemma \ref{[Appendix Lemma] LSC of Phi (global)} shows that the bound(s) obtained in Lemma \ref{[Appendix Lemma] LSC of f (local)} can be improved under mild conditions on $\Phi$.

\begin{lemma}\label{[Appendix Lemma] LSC of Phi (local)}
Suppose the assumptions of Lemma \ref{[Appendix Lemma] Phi- circ Phi} hold. Then, for every $x \in I_{\x}^{\delta}$ and every $\theta \in \R$,
\begin{equation*}
    |\LSC_{\Phi} (\theta \Phi) (x) - \theta \Phi(x)| \leq |\theta| \sup_{x' \in I_{\x}^{\delta + \epsilon}}|\Phi(x') - \Phi(x'-)|.
\end{equation*}
\end{lemma}

\begin{proof}
First, if $\theta > 0$, then the result follows from the fact that, by (\ref{Appendix eq: LSC representation}),
\begin{equation*}
    \LSC_{\Phi} (\theta \Phi) (x) = \theta \Phi(\Phi^-(\Phi(x))-) \leq \theta \Phi(\Phi^-(\Phi(x))) = \theta \Phi(x),
\end{equation*}
where, if $x \in I_{\x}^{\delta}$, then $\Phi^-(\Phi(x)) \in I_{\x}^{\delta + \epsilon}$ by Lemma \ref{[Appendix Lemma] Phi- circ Phi}, and therefore
\begin{equation*}
    \Phi(x) = \Phi(\Phi^-(\Phi(x))) \leq \Phi(\Phi^-(\Phi(x))-) + \sup_{x' \in I_{\x}^{\delta + \epsilon}}|\Phi(x') - \Phi(x'-)|.
\end{equation*}

Next, if $\theta < 0$, then the result follows from the fact that, by (\ref{Appendix eq: LSC representation}),
\begin{equation*}
    \LSC_{\Phi} (\theta \Phi)(x) = \theta \Phi(\Phi^-(\Phi(x)+)) \leq \theta \Phi(\Phi^-(\Phi(x))) = \theta \Phi(x),
\end{equation*}
where, if $x \in I_{\x}^{\delta}$, then $\Phi^-(\Phi(x)+) \in I_{\x}^{\delta + \epsilon}$ by Lemma \ref{[Appendix Lemma] Phi- circ Phi}, and therefore
\begin{align*}
    \Phi(\Phi^-(\Phi(x)+)) &\leq \liminf_{\eta \downarrow 0} \Phi(\Phi^-(\Phi(x) + \eta)-) + \sup_{x' \in I_{\x}^{\delta + \epsilon}}|\Phi(x') - \Phi(x'-)| \\
    &\leq \Phi(x)+ \sup_{x' \in I_{\x}^{\delta + \epsilon}}|\Phi(x') - \Phi(x'-)|.
\end{align*}
\end{proof}

\subsection{Proof of Theorem \ref{M-[Theorem] Main result}}

\paragraph*{Proof of \eqref{M-[Asymptotics] thetahat}}
Let $t \in \R$ be given. By Lemma \ref{M-[Lemma] Generalized Switch Relation} and change of variables,
\begin{align*}
    &\P \left[ r_n(\widehat{\theta}_n(\x) - \theta_0(\x) ) > t \right] \\
    &= \P \left[ \max\argmax_{x \in \widehat{\Phi}_n^{-}([0,\widehat{u}_n])} \left\{ [\theta_0(\x) + t r_n^{-1}] \widehat{\Phi}_n(x) - \LSC_{\widehat{\Phi}_n} (\widehat{\Gamma}_n) (x) \right\} < \widehat{\Phi}_n^{-} \circ \widehat{\Phi}_n (\x) \right] \\
    &= \P \left[ \max\argmin_{v \in \widehat{V}_{\x,n}^{\q}} \widehat{H}_{\x,n}^{\q}(v;t) < \widehat{Z}_{\x,n}^{\q} \right],
\end{align*}
where
\begin{equation*}
    \widehat{V}_{\x,n}^{\q} = \left\{a_n (x - \x) : x \in \widehat{\Phi}_{n}^{-}([0,\widehat{u}_{n}])\right\},
\end{equation*}
\begin{equation*}
    \widehat{H}_{\x,n}^{\q}(v;t) = \LSC_{\widehat{L}_{\x,n}^{\q}} \left(\widehat{G}_{\x,n}^{\q} + M_{\x,n}^{\q} + r_n \theta_0(\x) \widehat{L}_{\x,n}^{\q} \right)(v)  - [r_n \theta_0(\x) + t] \widehat{L}_{\x,n}^{\q}(v),
\end{equation*}
and
\begin{equation*}
    \widehat{Z}_{\x,n}^{\q} = a_n \left[\widehat{\Phi}_{n}^{-} \circ \widehat{\Phi}_{n}(\x) - \x \right],
\end{equation*}
with
\begin{align*}
    \widehat{G}_{\x,n}^{\q}(v) &= \sqrt{n a_n} \left[\widehat{\Gamma}_n(\x + v a_n^{-1}) - \widehat{\Gamma}_n(\x) - \Gamma_0(\x + v a_n^{-1}) + \Gamma_0(\x) \right] \\
    &\qquad - \theta_0(\x) \sqrt{n a_n} \left[\widehat{\Phi}_n(\x + v a_n^{-1})-\widehat{\Phi}_n(\x) - \Phi_0(\x + v a_n^{-1}) + \Phi_0(\x)\right], \\
    M_{\x,n}^{\q}(v) 
    &= \sqrt{n a_n} \left[\Gamma_{0}(\x + v a_n^{-1}) - \Gamma_{0}(\x)\right] - \theta_{0}(\x) \sqrt{n a_n} \left[\Phi_{0}(\x + v a_n^{-1})-\Phi_{0}(\x)\right], \\
    \widehat{L}_{\x,n}^{\q}(v) &= a_n \left[\widehat{\Phi}_{n}(\x + v a_n^{-1}) - \widehat{\Phi}_{n}(\x)\right].
\end{align*}
By \ref{M-Assumption B: uniform convergence Phi} and Lemma \ref{[Appendix Lemma] Phi- circ Phi}, $\widehat{Z}_{\x,n}^{\q} = o_{\P}(1)$. Suppose also that
\begin{equation}\label{Appendix eq: argmax convergence}
    \max\argmin_{v \in \widehat{V}_{\x,n}^{\q}} \widehat{H}_{\x,n}^{\q}(v;t) \leadsto \argmin_{v \in \R} \mathcal{H}_{\x}^{\q} (v;t),
\end{equation}
where $\mathcal{H}_{\x}^{\q}(v;t) = \mathcal{G}_{\x}(v) + \mathcal{M}_{\x}^{\q}(v)  - t \partial \Phi_0(\x) v$. Then
\begin{align*}
    \P \left[ r_n (\widehat{\theta}_n(\x) - \theta_0(\x)) > t \right]
    &= \P \left[ \max\argmin_{v \in \widehat{V}_{\x,n}^{\q}} \widehat{H}_{\x,n}^{\q}(v;t) < \widehat{Z}_{\x,n}^{\q} \right] \\
    &\to \P \left[ \argmin_{v \in \R} \mathcal{H}_{\x}^{\q}(v;t) < 0 \right] \\
    &= \P \left[ \frac{1}{\partial\Phi_0(\x)}\partial_{-}\GCM_{\R}(\mathcal{G}_{\x} + \mathcal{M}_{\x}^{\q}(0) > t \right],
\end{align*}
where the second line uses Lemma \ref{M-[Lemma] Continuity of argmax} and where the last equality uses Lemma \ref{[Appendix Lemma] Switching lemma for R}. The proof of \eqref{M-[Asymptotics] thetahat} can therefore be completed by showing \eqref{Appendix eq: argmax convergence}.

We shall do so by means of the argmax continuous mapping theorem of \citet{Cox_2022}. To be specific, using that theorem it can be shown that \eqref{Appendix eq: argmax convergence} holds if
\begin{equation}\label{Appendix eq: vhat tight}
    \widehat{v}_n(t) = \max\argmin_{v \in \widehat{V}_{\x,n}^{\q}} \widehat{H}_{\x,n}^{\q}(v;t) = O_{\P}(1)
\end{equation}
and if
\begin{equation}\label{Appendix eq: Hhat convergence}
    \widehat{H}_{\x,n}^{\q}(\cdot;t) \leadsto \mathcal{H}_{\x}^{\q} (\cdot;t).
\end{equation}

We begin by showing \eqref{Appendix eq: Hhat convergence}. First, by \ref{M-Assumption B: weak convergence}, $\widehat{G}_{\x,n}^{\q} \leadsto \mathcal{G}_{\x}$. Also, by \ref{M-Assumption A: theta0} and \ref{M-Assumption A: Phi0}, as $u \to 0$,
\begin{align*}
    \frac{\Gamma_0(\x+u) - \theta_0(\x) \Phi_0(\x + u) -\Gamma_0(\x) + \theta_0(\x) \Phi_0(\x)}{u^{\q+1}} 
    &\to \lim_{u \to 0}\frac{\left[\theta_0(\x + u)-\theta_0(\x)\right] \partial \Phi_0(\x + u)}{(\q+1) u^{\q}} \\
    &= \frac{\partial^{\q} \theta_0(\x) \partial \Phi_0(\x)}{(\q+1)!},
\end{align*}
where the first equality uses L'H\^{o}pital's rule and
\begin{equation*}
    \partial \left[ \Gamma_0(\x + u)-\theta_0(\x) \Phi_0(\x + u) -\Gamma_0(\x) + \theta_0(\x) \Phi_0(\x) \right] = \left[ \theta_0(\x + u)-\theta_0(\x) \right] \partial \Phi_0(\x + u).
\end{equation*}
As a consequence, $M_{\x,n}^{\q} \leadsto \mathcal{M}_{\x}^{\q}$. Moreover, $\widehat{L}_{\x,n}^{\q} - L_{\x,n}^{\q} \leadsto 0$ by \ref{M-Assumption B: uniform convergence Phi} and $L_{\x,n}^{\q} \leadsto \mathcal{L}_{\x}$ by \ref{M-Assumption A: Phi0}, where  
$L_{\x,n}^{\q}(v) = a_n \left[\Phi_0(\x + v a_n^{-1}) - \Phi_0(\x)\right]$. 
In particular, $\widehat{L}_{\x,n}^{\q} \leadsto \mathcal{L}_{\x}$ and therefore
\begin{equation*}
    (\widehat{G}_{\x,n}^{\q},M_{\x,n}^{\q},\widehat{L}_{\x,n}^{\q}) \leadsto (\mathcal{G}_{\x},\mathcal{M}_{\x}^{\q},\mathcal{L}_{\x}).
\end{equation*}

Because $\widehat{G}_{\x,n}^{\q}+M_{\x,n}^{\q}$ is asymptotically equicontinuous,
\begin{equation*}
    \LSC_{\widehat{L}_{\x,n}^{\q}}\left(\widehat{G}_{\x,n}^{\q} + M_{\x,n}^{\q}\right) - \left(\widehat{G}_{\x,n}^{\q} + M_{\x,n}^{\q}\right) \leadsto 0
\end{equation*}
by \ref{M-Assumption B: uniform convergence Phi} and Lemma \ref{[Appendix Lemma] LSC of f (local)}. Also, by \ref{M-Assumption B: jump of Phi hat} and Lemma \ref{[Appendix Lemma] LSC of Phi (local)},
\begin{equation*}
    \LSC_{\widehat{L}_{\x,n}^{\q}}\left(r_n \theta_0(\x)  \widehat{L}_{\x,n}^{\q}\right) - \left(r_n \theta_0(\x) \widehat{L}_{\x,n}^{\q}\right) \leadsto 0.
\end{equation*}
The result \eqref{Appendix eq: Hhat convergence} follows from the three preceding displays and the fact that, on $\widehat{V}_{\x,n}^{\q}$,
\begin{align*}
    0 
    &\geq \LSC_{\widehat{L}_{\x,n}^{\q}}\left(\widehat{G}_{\x,n}^{\q} + M_{\x,n}^{\q} + r_n \theta_0(\x) \widehat{L}_{\x,n}^{\q}\right) - \left(\widehat{G}_{\x,n}^{\q} + M_{\x,n}^{\q} + r_n \theta_0(\x) \widehat{L}_{\x,n}^{\q}\right) \\
    &\geq \LSC_{\widehat{L}_{\x,n}^{\q}}\left(\widehat{G}_{\x,n}^{\q} + M_{\x,n}^{\q}\right) - \left(\widehat{G}_{\x,n}^{\q} + M_{\x,n}^{\q}\right) + \LSC_{\widehat{L}_{\x,n}^{\q}}\left(r_n \theta_0(\x) \widehat{L}_{\x,n}^{\q}\right) - \left(r_n \theta_0(\x) \widehat{L}_{\x,n}^{\q}\right).    
\end{align*}

Next, to show \eqref{Appendix eq: vhat tight}, we first define $\theta_n(\x;t) = \theta_0(\x) + t r_n^{-1}$ and note that
\begin{equation*}
    \widehat{v}_n(t) = a_n \left[\max\argmax_{x \in \widehat{\Phi}_n^{-}([0,\widehat{u}_n])} \left\{\theta_n(\x;t) \widehat{\Phi}_n(x) - \LSC_{\widehat{\Phi}_n}\left(\widehat{\Gamma}_n\right)(x)\right\} - \x\right].
\end{equation*}
Now, if $|\widehat{v}_n(t)| > a_n \delta > 0,$ then
\begin{equation*}
    \sup_{x \not \in I_{\x}^{\delta}} \left\{\theta_n(\x;t) \widehat{\Phi}_n(x) - \LSC_{\widehat{\Phi}_n} \left(\widehat{\Gamma}_n\right)(x) \right\} \geq \theta_n(\x;t) \widehat{\Phi}_n(\x) - \LSC_{\widehat{\Phi}_n}\left(\widehat{\Gamma}_n\right)(\x),
\end{equation*}
where $|\theta_n(\x;t) - \theta_0(\x)| = O(r_n^{-1}) = o(1)$, and, by \ref{M-Assumption B: uniform convergence Phi},
\begin{equation*}
    \sup_{x \in I}\left|\widehat{\Phi}_n(x) - \Phi_0(x)\right| = o_{\P}(1).
\end{equation*}
Also, using \ref{M-Assumption B: uniform consistency Gamma}, \ref{M-Assumption B: uniform convergence Phi}, and Lemma \ref{[Appendix Lemma] LSC of Gamma_0 (global)},
\begin{align*}
    &\sup_{x \in I}\left|\LSC_{\widehat{\Phi}_n}\left(\widehat{\Gamma}_n\right)(x) - \Gamma_0(x)\right| \\
    &\leq \sup_{x \in I} \left|\LSC_{\widehat{\Phi}_n}\left(\widehat{\Gamma}_n\right)(x) - \LSC_{\widehat{\Phi}_n}\left(\Gamma_0\right)(x)\right| + \sup_{x \in I}\left|\LSC_{\widehat{\Phi}_n}\left(\Gamma_0\right)(x) - \Gamma_0(x)\right| \\
    &\leq \sup_{x \in I} \left|\widehat{\Gamma}_n(x) - \Gamma_0(x)\right| + 2 \left(\sup_{x \in I}|\theta_0 (x)|\right) \sup_{x \in I}|\widehat{\Phi}_n(x) - \Phi_0(x)| = o_{\P}(1).
\end{align*}
As a consequence, $\widehat{v}_n(t) = o_{\P}(a_n)$: For any $\delta > 0$,
\begin{align*}
    &\P\left[|\widehat{v}_n(t)| > a_n \delta \right] \\
    &\leq \P\left[\sup_{x \not \in I_{\x}^{\delta}} \left\{\theta_0(\x) \Phi_0(x) - \Gamma_0(x) \right\} \geq \theta_0(\x) \Phi_0(\x) - \Gamma_0(\x) + o_{\P}(1)\right] = o(1),
\end{align*}
where the equality uses the fact, noted by \citet[Supplement, proof of Lemma 3]{Westling-Carone_2020_AoS}, that the function  $v \mapsto \theta_0(\x) \Phi_0(v) - \Gamma_0(v)$ is unimodal and maximized at $v = \x$.

Next, defining $\widehat{V}_{\x,n}^{\q}(j) = \{v\in \widehat{V}_{\x,n}^{\q} : 2^j < |v| \leq 2^{j+1}\}$ and using $\widehat{v}_n(t) = o_{\P}(a_n)$, we have, for any $K$, any positive $\delta'$, and any sequence of events $\{\mathcal{A}_n'\}$ with $\lim_{n \to \infty}\P[\mathcal{A}_n'] = 1$,
\begin{align*}
    \limsup_{n \to \infty}\P \left[ |\widehat{v}_n(t)| > 2^K \right]
    &\leq \limsup_{n \to \infty}\sum_{j \geq K : 2^j \leq a_n \delta'} \P \left[ 2^j < |\widehat{v}_n(t)| \leq 2^{j+1} \cap \mathcal{A}_n' \right] \\
    &\leq \limsup_{n \to \infty} \sum_{j \geq K : 2^j \leq a_n \delta'} \P\left[\inf_{v \in \widehat{V}_{\x,n}^{\q}(j)} \widehat{H}_{\x,n}^{\q}(v;t) \leq \widehat{H}_{\x,n}^{\q}(0;t) \cap \mathcal{A}_n' \right].
\end{align*}
The proof of \eqref{Appendix eq: vhat tight} can therefore be completed by showing that the majorant side in the display can be made arbitrarily small by choice of $K$, $\delta'$, and $\{\mathcal{A}_n'\}$.

To do so, we begin by analyzing each term in the basic bound 
\begin{align*}
    \widehat{H}_{\x,n}^{\q}(v;t) - \widehat{H}_{\x,n}^{\q}(0;t)
    &\geq \LSC_{\widehat{L}_{\x,n}^{\q}}\left(\widehat{G}_{\x,n}^{\q}\right)(v) + \LSC_{\widehat{L}_{\x,n}^{\q}}\left(M_{\x,n}^{\q}\right)(v) - t \widehat{L}_{\x,n}^{\q}(v) \\
    &+ \LSC_{\widehat{L}_{\x,n}^{\q}}\left(r_n \theta_0(\x) \widehat{L}_{\x,n}^{\q}\right)(v) - r_n \theta_0(\x) \widehat{L}_{\x,n}^{\q} (v) - \widehat{H}_{\x,n}^{\q}(0;t).
\end{align*}
Because $\widehat{H}_{\x,n}^{\q}(0;t) \leadsto \mathcal{H}_{\x}^{\q}(0;t) = 0$ and because, by \ref{M-Assumption B: jump of Phi hat} and Lemma \ref{[Appendix Lemma] LSC of Phi (local)}, there is a positive $\delta'$ such that
\begin{equation*}
    \sup_{|v|\leq a_n \delta'} \left| \LSC_{\widehat{L}_{\x,n}^{\q}} \left( r_n\theta_0(\x) \widehat{L}_{\x,n}^{\q} \right) (v) - r_n\theta_0(\x) \widehat{L}_{\x,n}^{\q} (v) \right| = o_{\P}(1),
\end{equation*}
we may assume that, on $\{\mathcal{A}_n'\}$ and for some $C_0$,
\begin{equation*}
    \sup_{|v| \leq a_n \delta'} \left| \LSC_{\widehat{L}_{\x,n}^{\q}} \left( r_n \theta_0(\x) \widehat{L}_{\x,n}^{\q} \right) (v) - r_n \theta_0(\x) \widehat{L}_{\x,n}^{\q}(v) - \widehat{H}_{\x,n}^{\q}(0;t) \right| \leq C_0.
\end{equation*}
Also, because, by \ref{M-Assumption B: uniform convergence Phi} and \ref{M-Assumption A: Phi0}, there is a positive $\delta'$ such that
\begin{equation*}
    \sup_{|v| \leq a_n \delta'} \left| \widehat{L}_{\x,n}^{\q}(v) - L_{\x,n}^{\q} (v) \right| = o_{\P}(1) \qquad \text{and} \qquad \sup_{1 \leq |v| \leq a_n \delta'} \left| \frac{L_{\x,n}^{\q}(v)}{\mathcal{L}_{\x} (v)} \right| < \infty,
\end{equation*}
we may assume that, on $\{\mathcal{A}_n'\}$ and for some $C_L$,
\begin{equation*}
    \sup_{1 \leq |v| \leq a_n \delta'}\left|\frac{\widehat{L}_{\x,n}^{\q}(v)}{v}\right| \leq C_L.
\end{equation*}
Next, by \ref{M-Assumption B: uniform convergence Phi} and Lemma \ref{[Appendix Lemma] LSC of f (local)}, with probability approaching one,
\begin{equation*}
    \LSC_{\widehat{L}_{\x,n}^{\q}} \left( M_{\x,n}^{\q} \right) (v) \geq \inf_{|v|/2 \leq |v'| \leq 2|v|} M_{\x,n}^{\q}(v') \qquad \text{for every } |v| \geq 2,
\end{equation*}
while, by \ref{M-Assumption A: theta0} and \ref{M-Assumption A: Phi0}, there is a positive $\delta'$ such that
\begin{equation*}
    \inf_{1 \leq |v| \leq \eta a_n \delta'}\frac{M_{\x,n}^{\q}(v)}{\mathcal{M}_{\x}^{\q} (v)} > 0.
\end{equation*}
We may therefore assume that, on $\{\mathcal{A}_n'\}$ and for some positive $C_M$,
\begin{equation*}
    \inf_{2 \leq |v| \leq a_n \delta'}\frac{\LSC_{\widehat{L}_{\x,n}^{\q}} \left (M_{\x,n}^{\q} \right) (v)}{v^{\q+1}} \geq C_M.
\end{equation*}
Finally, by \ref{M-Assumption B: uniform convergence Phi} and Lemma \ref{[Appendix Lemma] LSC of f (local)}, with probability approaching one,
\begin{equation*}
    \left| \LSC_{\widehat{L}_{\x,n}^{\q}} \left( \widehat{G}_{\x,n}^{\q} \right) (v) \right| \leq \sup_{|v| / 2 \leq |v'| \leq 2 |v|} \left| \widehat{G}_{\x,n}^{\q}(v') \right| \qquad \text{for every } |v| \geq 2,
\end{equation*}
and we may therefore assume that, on $\{\mathcal{A}_n'\}$,
\begin{equation*}
    \sup_{v \in \widehat{V}_{\x,n}^{\q}(j)} \left| \LSC_{\widehat{L}_{\x,n}^{\q}} \left( \widehat{G}_{\x,n}^{\q} \right)(v) \right| \leq \sup_{v \in V_2(j)} \left| \widehat{G}_{\x,n}^{\q}(v) \right| \qquad \text{for every } j \geq 2 \text{ with } 2^j \leq a_n \delta',
\end{equation*}
where $V_\eta(j) = \{v \in \R: \eta^{-1}2^j \leq |v| \leq \eta 2^{j+1}\}$.

As a consequence, by the Markov inequality,
\begin{align*}
    &\P \left[ \inf_{v \in \widehat{V}_{\x,n}^{\q}(j)} \widehat{H}_{\x,n}^{\q}(v;t) \leq \widehat{H}_{\x,n}^{\q}(0;t) \cap \mathcal{A}_n' \right] \\
    &\leq \P \left[ \sup_{v \in V_2(j)} \left| \widehat{G}_{\x,n}^{\q}(v) \right| \geq \inf_{v \in V_1(j)} \left[ C_M v^{\q+1} - C_L |v| - C_0 \right] \cap \mathcal{A}_n' \right] \\
    &\leq \frac{\E \left[ \sup_{v \in V_2(j)} \left| \widehat{G}_{\x,n}^{\q}(v) \right| \1_{\mathcal{A}_n'} \right]}{\inf_{v \in V_1(j)} \left[C_M v^{\q+1} - C_L |v| - C_0\right]} \qquad \text{for every } j \geq 2 \text{ with } 2^j \leq a_n \delta',
\end{align*}
where, by \ref{M-Assumption B: moment bound}, we may assume that, for some $C_G$,
\begin{equation*}
    \E \left[ \sup_{v \in V_2(j)} \left| \widehat{G}_{\x,n}^{\q}(v) \right| \1_{\mathcal{A}_n'} \right] \leq \sum_{j-1 \leq j' \leq j+1} \E \left[ \sup_{v \in V_1(j')} \left| \widehat{G}_{\x,n}^{\q}(v) \right| \1_{\mathcal{A}_n'} \right] \leq C_G 2^{j \beta},
\end{equation*}
and where, for all sufficiently large $j$,
\begin{equation*}
    \inf_{v \in V_1(j)} \left[ C_M v^{\q+1} - C_L |v| - C_0 \right] \geq \frac{1}{2} C_M 2^{j(\q+1)}.
\end{equation*}
In other words, for large $K$,
\begin{equation*}
    \limsup_{n \to \infty} \sum_{j \geq K : 2^j \leq a_n \delta'} \P \left[ \inf_{v \in \widehat{V}_{\x,n}^{\q}(j)} \widehat{H}_{\x,n}^{\q}(v;t) \leq \widehat{H}_{\x,n}^{\q}(0;t) \cap \mathcal{A}_n' \right] \leq \frac{2 C_G}{C_M}\sum_{j \geq K} 2^{j[\beta - (\q + 1)]},
\end{equation*}
which can be made arbitrarily small by choice of $K$.

\paragraph*{Proof of \eqref{M-[Asymptotics] thetatildestar}}
We proceed as in the proof of \eqref{M-[Asymptotics] thetahat}. Let $t \in \R$ be given. By Lemma \ref{M-[Lemma] Generalized Switch Relation} and change of variables,
\begin{align*}
    &\P_n^* \left[ r_n(\widetilde{\theta}_n^*(\x) - \widehat{\theta}_n(\x) ) > t \right] \\
    &= \P_n^* \left[ \max\argmax_{x \in \widehat{\Phi}_n^{*-}([0,\widehat{u}_n^*])} \left\{[\widehat{\theta}_n(\x) + t r_n^{-1}] \widehat{\Phi}_n^*(x) - \LSC_{\widehat{\Phi}_n^*} (\widetilde{\Gamma}_n^*)(x) \right\} < \widehat{\Phi}_n^{*-} \circ \widehat{\Phi}_n^*(\x) \right] \\
    &= \P_n^* \left[ \max\argmin_{v \in \widehat{V}_{\x,n}^{\q,*}} \widehat{H}_{\x,n}^{\q,*}(v;t) < \widehat{Z}_{\x,n}^{\q,*} \right],
\end{align*}
where
\begin{equation*}
    \widehat{V}_{\x,n}^{\q,*} = \left\{a_n (x - \x) : x \in \widehat{\Phi}_{n}^{*-}([0,\widehat{u}_{n}^*])\right\},
\end{equation*}
\begin{equation*}
    \widehat{H}_{\x,n}^{\q,*}(v;t) = \LSC_{\widehat{L}_{\x,n}^{\q,*}} \left( \widehat{G}_{\x,n}^{\q,*} + \widetilde{M}_{\x,n}^{\q} + r_n \widehat{\theta}_n(\x) \widehat{L}_{\x,n}^{\q,*} \right) (v)  - [r_n \widehat{\theta}_n(\x) + t] \widehat{L}_{\x,n}^{\q,*}(v)
\end{equation*}
and
\begin{equation*}
    \widehat{Z}_{\x,n}^{\q,*} = a_n \left[ \widehat{\Phi}_{n}^{*-} \circ \widehat{\Phi}_{n}^*(\x) - \x \right],
\end{equation*}
with
\begin{align*}
    \widehat{G}_{\x,n}^{\q,*}(v)
    &= \sqrt{n a_n} \left[\widehat{\Gamma}_n^*(\x + v a_n^{-1}) - \widehat{\Gamma}_n^*(\x) - \widehat{\Gamma}_n(\x + v a_n^{-1}) + \widehat{\Gamma}_n(\x) \right] \\
    &\qquad - \widehat{\theta}_n(\x) \sqrt{n a_n} \left[\widehat{\Phi}_n^*(\x + v a_n^{-1}) - \widehat{\Phi}_n^*(\x) - \widehat{\Phi}_n(\x + v a_n^{-1}) + \widehat{\Phi}_n(\x)\right], \\
    \widetilde{M}_{\x,n}^{\q}(v) 
    &= \sqrt{n a_n} \widetilde{M}_{\x,n}(v a_n^{-1}), \qquad \widehat{L}_{\x,n}^{\q,*}(v) = a_n \left[ \widehat{\Phi}_{n}^*(\x + v a_n^{-1}) - \widehat{\Phi}_{n}^*(\x) \right].
\end{align*}
By \ref{M-Assumption B: uniform convergence Phi} and Lemma \ref{[Appendix Lemma] Phi- circ Phi}, $\widehat{Z}_{\x,n}^{\q,*} = o_{\P}(1)$. Suppose also that
\begin{equation}\label{Appendix eq: bootstrap argmax convergence}
    \max\argmin_{v \in \widehat{V}_{\x,n}^{\q,*}} \widehat{H}_{\x,n}^{\q,*}(v;t) \leadsto_{\P} \argmin_{v \in \R} \mathcal{H}_{\x}^{\q} (v;t).
\end{equation}
Then, as in the proof of \eqref{M-[Asymptotics] thetahat},
\begin{equation*}
    \P_n^* \left[ r_n(\widetilde{\theta}_n^*(\x) - \widehat{\theta}_n(\x) ) > t \right] \to_{\P} \P \left[ \frac{1}{\partial\Phi_0(\x)}\partial_{-}\GCM_{\R}(\mathcal{G}_{\x} + \mathcal{M}_{\x}^{\q}(0) > t \right].
\end{equation*}
The proof of \eqref{M-[Asymptotics] thetatildestar} can therefore be completed by showing \eqref{Appendix eq: bootstrap argmax convergence}.

We shall do so by showing that
\begin{equation}\label{Appendix eq: vhatstar tight}
    \widehat{v}_n^*(t) = \max\argmin_{v \in \widehat{V}_{\x,n}^{\q,*}} \widehat{H}_{\x,n}^{\q,*}(v;t) = O_{\P}(1)
\end{equation}
and
\begin{equation}\label{Appendix eq: Hhatstar convergence}
    \widehat{H}_{\x,n}^{\q,*}(\cdot;t) \leadsto_{\P} \mathcal{H}_{\x}^{\q} (\cdot;t).
\end{equation}

We begin by showing \eqref{Appendix eq: Hhatstar convergence}. First, by \ref{M-Assumption B: weak convergence}, $\widehat{G}_{\x,n}^{\q,*} \leadsto_{\P} \mathcal{G}_{\x}$. Also, by Assumption \ref{M-Assumption C}, $\widetilde{M}_{\x,n}^{\q,*} \leadsto_{\P} \mathcal{M}_{\x}^{\q}$. Moreover, $\widehat{L}_{\x,n}^{\q,*} - \widehat{L}_{\x,n}^{\q} \leadsto_{\P} 0$ by \ref{M-Assumption B: uniform convergence Phi}, where, as shown in the proof of (\ref{Appendix eq: argmax convergence}), $\widehat{L}_{\x,n}^{\q} \leadsto \mathcal{L}_{\x}$. In particular, $\widehat{L}_{\x,n}^{\q,*} \leadsto_{\P} \mathcal{L}_{\x}$ and therefore
\begin{equation*}
    (\widehat{G}_{\x,n}^{\q,*} , \widetilde{M}_{\x,n}^{\q} , \widehat{L}_{\x,n}^{\q,*}) \leadsto_{\P} (\mathcal{G}_{\x}^{\q} , \mathcal{M}_{\x}^{\q} , \mathcal{L}_{\x}).
\end{equation*}

Because $\widehat{G}_{\x,n}^{\q,*} + \widetilde{M}_{\x,n}^{\q}$ is asymptotically equicontinuous,
\begin{equation*}
    \LSC_{\widehat{L}_{\x,n}^{\q,*}} \left( \widehat{G}_{\x,n}^{\q,*} + \widetilde{M}_{\x,n}^{\q} \right) - \left( \widehat{G}_{\x,n}^{\q,*} + \widetilde{M}_{\x,n}^{\q} \right) \leadsto_{\P} 0
\end{equation*}
by \ref{M-Assumption B: uniform convergence Phi} and Lemma \ref{[Appendix Lemma] LSC of f (local)}. Also, by \ref{M-Assumption B: jump of Phi hat} and Lemma \ref{[Appendix Lemma] LSC of Phi (local)},
\begin{equation*}
    \LSC_{\widehat{L}_{\x,n}^{\q,*}} \left( r_n \widehat{\theta}_n(\x)  \widehat{L}_{\x,n}^{\q,*} \right) - \left( r_n \widehat{\theta}_n(\x) \widehat{L}_{\x,n}^{\q,*} \right) \leadsto_{\P} 0.
\end{equation*}
The result \eqref{Appendix eq: Hhatstar convergence} follows from the three preceding displays and the fact that, on $\widehat{V}_{\x,n}^{\q,*}$,
\begin{align*}
    0
    &\geq \LSC_{\widehat{L}_{\x,n}^{\q,*}} \left( \widehat{G}_{\x,n}^{\q,*} + \widetilde{M}_{\x,n}^{\q} + r_n \widehat{\theta}_n(\x) \widehat{L}_{\x,n}^{\q,*} \right) - \left( \widehat{G}_{\x,n}^{\q,*} + \widetilde{M}_{\x,n}^{\q} + r_n \widehat{\theta}_n(\x) ,\widehat{L}_{\x,n}^{\q,*} \right) \\
    &\geq \LSC_{\widehat{L}_{\x,n}^{\q,*}} \left(\widehat{G}_{\x,n}^{\q,\,*} + \widetilde{M}_{\x,n}^{\q}\right) - \left( \widehat{G}_{\x,n}^{\q,*} + \widetilde{M}_{\x,n}^{\q} \right) \\
    &\qquad + \LSC_{\widehat{L}_{\x,n}^{\q,*}} \left(r_n \widehat{\theta}_n(\x) \widehat{L}_{\x,n}^{\q,*} \right) - \left( r_n \widehat{\theta}_n(\x) \widehat{L}_{\x,n}^{\q,*} \right).    
\end{align*}

Next, to show \eqref{Appendix eq: vhatstar tight}, we first define $\widehat{\theta}_n(\x;t) = \widehat{\theta}_n(\x) + t r_n^{-1}$ and note that
\begin{equation*}
    \widehat{v}_n^*(t) = a_n \left[ \max\argmax_{x \in \widehat{\Phi}_n^{*-}([0,\widehat{u}_n^*])} \left\{ \widehat{\theta}_n(\x;t) \widehat{\Phi}_n^*(x) - \LSC_{\widehat{\Phi}_n^*}\left( \widetilde{\Gamma}_n^* \right)(x) \right\} - \x \right].
\end{equation*}
Now, if $|\widehat{v}_n^*(t)| > a_n \delta > 0,$ then
\begin{equation*}
    \sup_{x \not \in I_{\x}^{\delta}} \left\{ \widehat{\theta}_n(\x;t) \widehat{\Phi}_n^*(x) - \LSC_{\widehat{\Phi}_n^*} \left( \widetilde{\Gamma}_n^* \right)(x) \right\} \geq \widehat{\theta}_n(\x;t) \widehat{\Phi}_n^*(\x) - \LSC_{\widehat{\Phi}_n^*} \left( \widetilde{\Gamma}_n^* \right) (\x),
\end{equation*}
where $|\widehat{\theta}_n(\x;t) - \theta_0(\x)| = O_{\P}(r_n^{-1}) = o_{\P}(1)$, and, by \ref{M-Assumption B: uniform convergence Phi},
\begin{equation*}
    \sup_{x \in I} \left| \widehat{\Phi}_n^*(x) - \Phi_0(x) \right| = o_{\P}(1).
\end{equation*}
Therefore, defining $\widehat{\x}_n^* = \widehat{\Phi}_n^{*-}(\widehat{\Phi}_n^{*}(\x)) = \x + o_{\P}(1)$ and using (\ref{Appendix eq: LSC representation}),
\begin{align*}
    &\widehat{\theta}_n(\x;t) \widehat{\Phi}_n^*(\x) - \LSC_{\widehat{\Phi}_n^*} \left( \widetilde{\Gamma}_n^* \right) (\x) \\
    &\geq \widehat{\theta}_n(\x;t) \widehat{\Phi}_n^*(\x) - \widetilde{\Gamma}_n^*(\widehat{\x}_n^*) \\
    &= \widehat{\theta}_n(\x;t) \widehat{\Phi}_n^*(\x) - \widehat{\theta}_n(\x) \widehat{\Phi}_n(\widehat{\x}_n^*) - \widehat{\Gamma}_n^*(\widehat{\x}_n^*) + \widehat{\Gamma}_n(\widehat{\x}_n^*) - \widetilde{M}_{\x,n}^\q(\widehat{Z}_{\x,n}^{\q,*}) / \sqrt{n a_n} = o_{\P}(1),  
\end{align*}
where the last equality uses \ref{M-Assumption B: uniform consistency Gamma}, $\widehat{Z}_{\x,n}^{\q,*} = o_{\P}(1)$, and Assumption \ref{M-Assumption C}. Also, using \ref{M-Assumption B: uniform consistency Gamma}, \ref{M-Assumption B: uniform convergence Phi}, and Assumption \ref{M-Assumption C}, we have, uniformly in $x \not \in I_\x^\delta$ and for some $c > 0$,
\begin{align*}
    \widetilde{\Gamma}_n^*(x) - \widehat{\theta}_n(\x) \widehat{\Phi}_n^*(x)
    &= \widehat{\Gamma}_n^*(x) - \widehat{\Gamma}_n(x) - \widehat{\theta}_n(\x) \left(\widehat{\Phi}_n^*(x) - \widehat{\Phi}_n(x)\right) + \widetilde{M}_n(x - \x) \\
    &\geq c \delta^{\q + 1} + o_{\P}(1),
\end{align*}
and therefore, by Lemma \ref{[Appendix Lemma] LSC of Phi (global)},
\begin{equation*}
    \sup_{x \not \in I_{\x}^{\delta}} \left\{ \widehat{\theta}_n(\x;t) \widehat{\Phi}_n^*(x) - \LSC_{\widehat{\Phi}_n^*} \left( \widetilde{\Gamma}_n^* \right) (x) \right\} \leq - c \delta^{\q + 1} + o_{\P}(1).
\end{equation*}

As a consequence, $\widehat{v}_n^*(t) = o_{\P}(a_n)$: For any $\delta > 0$,
\begin{equation*}
    \P \left[ |\widehat{v}_n^*(t)| > a_n \delta \right] \leq \P \left[- c \delta^{\q + 1} \geq o_{\P}(1) \right] = o(1).
\end{equation*}

Next, defining $\widehat{V}_{\x,n}^{\q,*}(j) = \{v\in \widehat{V}_{\x,n}^{\q,*} : 2^j < |v| \leq 2^{j+1}\}$ and using $\widehat{v}_n^*(t) = o_{\P}(a_n)$, we have, for any $K$, any positive $\delta'$, and any sequence of events $\{\mathcal{A}_n'\}$ with $\lim_{n \to \infty}\P[\mathcal{A}_n'] = 1$,
\begin{align*}
    \limsup_{n \to \infty} \P \left[ |\widehat{v}_n^*(t)| > 2^K \right] 
    &\leq \limsup_{n \to \infty}\sum_{j \geq K : 2^j \leq a_n \delta'} \P \left[ 2^j < |\widehat{v}_n^*(t)| \leq 2^{j+1} \cap \mathcal{A}_n' \right] \\
    &\leq \limsup_{n \to \infty} \sum_{j \geq K : 2^j \leq a_n \delta'} \P \left[ \inf_{v \in \widehat{V}_{\x,n}^{\q,*}(j)} \widehat{H}_{\x,n}^{\q,*}(v;t) \leq \widehat{H}_{\x,n}^{\q,*}(0;t) \cap \mathcal{A}_n' \right].
\end{align*}
The proof of \eqref{Appendix eq: vhatstar tight} can therefore be completed by showing that the majorant side in the display can be made arbitrarily small by choice of $K$, $\delta'$, and $\{\mathcal{A}_n'\}$.

To do so, we begin by analyzing each term in the basic bound 
\begin{align*}
    \widehat{H}_{\x,n}^{\q,*}(v;t) - \widehat{H}_{\x,n}^{\q,*}(0;t) 
    &\geq \LSC_{\widehat{L}_{\x,n}^{\q,*}} \left( \widehat{G}_{\x,n}^{\q,*} \right) (v) + \LSC_{\widehat{L}_{\x,n}^{\q,*}} \left( \widetilde{M}_{\x,n}^{\q} \right) (v) - t \widehat{L}_{\x,n}^{\q,*}(v) \\
    &+ \LSC_{\widehat{L}_{\x,n}^{\q,*}} \left( r_n \widehat{\theta}_n(\x) \widehat{L}_{\x,n}^{\q,*} \right) (v) - r_n \widehat{\theta}_n(\x) \widehat{L}_{\x,n}^{\q,*} (v) - \widehat{H}_{\x,n}^{\q,*}(0;t).
\end{align*}
Because $\widehat{H}_{\x,n}^{\q,*}(0;t) \leadsto_{\P} \mathcal{H}_{\x}^{\q}(0;t) = 0$ and because, by \ref{M-Assumption B: jump of Phi hat} and Lemma \ref{[Appendix Lemma] LSC of Phi (local)}, there is a positive $\delta'$ such that
\begin{equation*}
    \sup_{|v|\leq a_n \delta'}\left|\LSC_{\widehat{L}_{\x,n}^{\q,*}}\left(r_n \widehat{\theta}_n(\x) \widehat{L}_{\x,n}^{\q,*}\right)(v) - r_n \widehat{\theta}_n(\x) \widehat{L}_{\x,n}^{\q,*} (v)\right| = o_{\P}(1),
\end{equation*}
we may assume that, on $\{\mathcal{A}_n'\}$ and for some $C_0$,
\begin{equation*}
    \sup_{|v| \leq a_n \delta'} \left| \LSC_{\widehat{L}_{\x,n}^{\q,*}} \left( r_n \widehat{\theta}_n(\x) \widehat{L}_{\x,n}^{\q,*} \right) (v) - r_n \widehat{\theta}_n(\x) \widehat{L}_{\x,n}^{\q,*}(v) - \widehat{H}_{\x,n}^{\q,*}(0;t) \right| \leq C_0.
\end{equation*}
Also, because, by \ref{M-Assumption B: uniform convergence Phi} and \ref{M-Assumption A: Phi0}, there is a positive $\delta'$ such that
\begin{equation*}
    \sup_{|v| \leq a_n \delta'} \left| \widehat{L}_{\x,n}^{\q,*}(v) - L_{\x,n}^{\q} (v) \right| = o_{\P}(1) \qquad \text{and} \qquad \sup_{1 \leq |v| \leq a_n \delta'} \left| \frac{L_{\x,n}^{\q}(v)}{\mathcal{L}_{\x} (v)} \right| < \infty,
\end{equation*}
we may assume that, on $\{\mathcal{A}_n'\}$ and for some $C_L$,
\begin{equation*}
    \sup_{1 \leq |v| \leq a_n \delta'} \left| \frac{\widehat{L}_{\x,n}^{\q,*}(v)}{v} \right| \leq C_L.
\end{equation*}
Next, by \ref{M-Assumption B: uniform convergence Phi} and Lemma \ref{[Appendix Lemma] LSC of f (local)}, with probability approaching one,
\begin{equation*}
    \LSC_{\widehat{L}_{\x,n}^{\q,*}} \left( \widetilde{M}_{\x,n}^{\q} \right)(v) \geq \inf_{|v|/2 \leq |v'| \leq 2|v|} \widetilde{M}_{\x,n}^{\q}(v') \qquad \text{for every } |v| \geq 2,
\end{equation*}
while, by Assumption \ref{M-Assumption C}, there is a positive $c$ such that, with probability approaching one,
\begin{equation*}
    \inf_{|v| \geq 1}\frac{\widetilde{M}_{\x,n}^{\q}(v)}{v^{\q + 1}} > c.
\end{equation*}
We may therefore assume that, on $\{\mathcal{A}_n'\}$ and for some positive $C_M$,
\begin{equation*}
    \inf_{2 \leq |v| \leq a_n \delta'}\frac{\LSC_{\widehat{L}_{\x,n}^{\q,*}} \left( \widetilde{M}_{\x,n}^{\q} \right) (v)}{v^{\q+1}} \geq C_M.
\end{equation*}
Finally, by \ref{M-Assumption B: uniform convergence Phi} and Lemma \ref{[Appendix Lemma] LSC of f (local)}, with probability approaching one,
\begin{equation*}
    \left| \LSC_{\widehat{L}_{\x,n}^{\q,*}} \left( \widehat{G}_{\x,n}^{\q,*} \right) (v) \right| \leq \sup_{|v| / 2 \leq |v'| \leq 2 |v|} \left| \widehat{G}_{\x,n}^{\q,*}(v') \right| \qquad \text{for every } |v| \geq 2,
\end{equation*}
and we may therefore assume that, on $\{\mathcal{A}_n'\}$,
\begin{equation*}
    \sup_{v \in \widehat{V}_{\x,n}^{\q,*}(j)} \left| \LSC_{\widehat{L}_{\x,n}^{\q,*}} \left( \widehat{G}_{\x,n}^{\q,*} \right)(v) \right| \leq \sup_{v \in V_2(j)} \left| \widehat{G}_{\x,n}^{\q,*}(v) \right| \qquad \text{for every } j \geq 2 \text{ with } 2^j \leq a_n \delta'.
\end{equation*}

As a consequence, by the Markov inequality,
\begin{align*}
    &\P \left[ \inf_{v \in \widehat{V}_{\x,n}^{\q,*}(j)} \widehat{H}_{\x,n}^{\q,*}(v;t) \leq \widehat{H}_{\x,n}^{\q,*}(0;t) \cap \mathcal{A}_n' \right] \\
    &\leq \P \left[ \sup_{v \in V_2(j)} \left| \widehat{G}_{\x,n}^{\q,*}(v) \right| \geq \inf_{v \in V_1(j)} \left[ C_M v^{\q+1} - C_L |v| - C_0 \right] \cap \mathcal{A}_n' \right] \\
    &\leq \frac{\E \left[ \sup_{v \in V_2(j)} \left| \widehat{G}_{\x,n}^{\q,*}(v) \right| \1_{\mathcal{A}_n'} \right]}{\inf_{v \in V_1(j)} \left[C_M v^{\q+1} - C_L |v| - C_0 \right]} \qquad \text{for every } j \geq 2 \text{ with } 2^j \leq a_n \delta',
\end{align*}
where, by \ref{M-Assumption B: uniform convergence Phi}, we may assume that, for some $C_G$,
\begin{equation*}
    \E \left[ \sup_{v \in V_2(j)} \left| \widehat{G}_{\x,n}^{\q,*}(v) \right| \1_{\mathcal{A}_n'} \right] \leq \sum_{j-1 \leq j' \leq j+1} \E \left[ \sup_{v \in V_1(j')} \left| \widehat{G}_{\x,n}^{\q,*}(v) \right| \1_{\mathcal{A}_n'} \right] \leq C_G 2^{j \beta},
\end{equation*}
and where, for all sufficiently large $j$,
\begin{equation*}
    \inf_{v \in V_1(j)} \left[ C_M v^{\q+1} - C_L |v| - C_0 \right] \geq \frac{1}{2} C_M 2^{j(\q+1)}.
\end{equation*}
In other words, for large $K$,
\begin{equation*}
    \limsup_{n \to \infty} \sum_{j \geq K : 2^j \leq a_n \delta'} \P \left[ \inf_{v \in \widehat{V}_{\x,n}^{\q,*}(j)} \widehat{H}_{\x,n}^{\q,*}(v;t) \leq \widehat{H}_{\x,n}^{\q,*}(0;t) \cap \mathcal{A}_n' \right] \leq \frac{2 C_G}{C_M}\sum_{j \geq K} 2^{j[\beta - (\q + 1)]},
\end{equation*}
which can be made arbitrarily small by choice of $K$.

\paragraph*{Proof of \eqref{M-Bootstrap consistency}} The bootstrap consistency result \eqref{M-Bootstrap consistency} follows from \eqref{M-[Asymptotics] thetahat}, \eqref{M-[Asymptotics] thetatildestar}, Polya's theorem, and the fact that, by Lemma \ref{M-[Lemma] Continuity of argmax},  the limiting distribution in \eqref{M-[Asymptotics] thetahat} and \eqref{M-[Asymptotics] thetatildestar} has a continuous cdf. \qed


\subsection{Proof of Lemma \ref{M-[Lemma] Mean function estimation with known q}}

For the monomial approximation estimator, we have
\begin{align*}
    \widetilde{\mathcal{D}}_{\q,n}^{\mathtt{MA}}(\x) 
    &= \epsilon_n^{-(\q + 1)}[\Gamma_0(\x + \epsilon_n) - \Gamma_0(\x) - \theta_0(\x) \{\Phi_0(\x + \epsilon_n) - \Phi_0(\x)\} ] \\
    &\quad + \epsilon_n^{-(\q + 1/2)}n^{-1/2} \widehat{G}_{\x,n}(1;\epsilon_n) \\
    &\quad - \epsilon_n^{-\q}[\widehat{\theta}_n(\x) -\theta_0(\x)] \widehat{R}_{\x,n}(1;\epsilon_n) \\
    &\quad -\epsilon_n^{-(\q + 1)}[\widehat{\theta}_n(\x) - \theta_0(\x)] [\Phi_0(\x + \epsilon_n) - \Phi_0(\x)] \\
    &= \frac{\partial^{\q} \theta_0(\x) \partial \Phi_0(\x)}{(\q + 1)!} + o(1) + O_{\P}[(n \epsilon_n^{1 + 2 \q})^{-1/2} + (n \epsilon_n^{1 + 2 \q})^{-\q / (1 + 2 \q)}] \\
    &= \frac{\partial^{\q} \theta_0(\x) \partial \Phi_0(\x)}{(\q + 1)!} + o_{\P}(1),
\end{align*}
where the second equality uses $\epsilon_n \to 0$ and the last equality uses $n \epsilon_n^{1 + 2 \q} \to \infty$.

Similarly, for the forward difference estimator, we have
\begin{align*}
    \widetilde{\mathcal{D}}_{\q,n}^{\mathtt{FD}}(\x) 
    &= \epsilon_n^{-(\q+1)}\sum_{k = 1}^{\q + 1} (-1)^{k + \q + 1} \binom{\q + 1}{k} [\Upsilon_0(\x + k \epsilon_n) - \Upsilon_0(\x) ] \\
    &\quad + \epsilon_n^{-(\q + 1/2)} n^{- 1/2}\sum_{k = 1}^{\q + 1}(-1)^{k + \q + 1} \binom{\q + 1}{k} \widehat{G}_{\x,n}(k;\epsilon_n) \\
    &\quad - \epsilon_n^{-\q}[\widehat{\theta}_n(\x) - \theta_0(\x)] \sum_{k = 1}^{\q + 1} (-1)^{k + \q + 1}\binom{\q + 1}{k} \widehat{R}_{\x,n}(k;\epsilon_n) \\
    &\quad -\epsilon_n^{-(\q + 1)}[\widehat{\theta}_n(\x) - \theta_0(\x)] \sum_{k = 1}^{\q + 1} (-1)^{k + \q + 1} \binom{\q + 1}{k} [\Phi_0(\x + k \epsilon_n) - \Phi_0(\x)] \\
    &= \frac{\partial^{\q} \theta_0(\x) \partial \Phi_0(\x)}{(\q + 1)!} + o(1) + O_{\P}[(n \epsilon_n^{1 + 2 \q})^{-1/2} + (n \epsilon_n^{1 + 2 \q})^{-\q / (1 + 2 \q)}] \\
    &= \frac{\partial^{\q} \theta_0(\x) \partial \Phi_0(\x)}{(\q + 1)!} + o_{\P}(1). \qed
\end{align*}

\subsection{Proof of Lemma \ref{M-[Lemma] Mean function estimation unknown q}}

Proceeding as in the proof of Lemma \ref{M-[Lemma] Mean function estimation with known q}, we have
\begin{align*}
    \widetilde{\mathcal{D}}_{j,n}^{\mathtt{BR}}(\x) &= \epsilon_n^{-(j + 1)} \sum_{k = 1}^{\underline{\s} + 1} \lambda_j^{\mathtt{BR}}(k) [\Upsilon_0(\x + c_k \epsilon_n) - \Upsilon_0(\x)] \\
    &\quad + \epsilon_n^{-(j + 1/2)} n^{-1/2} \sum_{k = 1}^{\underline{\s} + 1}\lambda_j^{\mathtt{BR}}(k) \widehat{G}_{\x,n}(c_k;\epsilon_n) \\
    &\quad - \epsilon_n^{-j}[\widehat{\theta}_n(\x) -\theta_0(\x)] \sum_{k = 1}^{\underline{\s} + 1} \lambda_j^{\mathtt{BR}}(k)\widehat{R}_{\x,n}(c_k;\epsilon_n) \\
    &\quad -\epsilon_n^{-(j + 1)}[\widehat{\theta}_n(\x)-\theta_0(\x)] \sum_{k = 1}^{\underline{\s} + 1} \lambda_j^{\mathtt{BR}}(k) [\Phi_0(\x + c_k \epsilon_n) - \Phi_0(\x)] \\
    &= \frac{\partial^{j + 1} \Upsilon_0(\x)}{(j + 1)!} + O(\epsilon_n^{\min \{\s,\underline{\s} + 1\} - j}) + O_{\P}(\epsilon_n^{-(j + 1/2)} n^{-1/2} + \epsilon_n^{-j} a_n^{-\q}) \\
    &= \mathcal{D}_j(\x) + O(\epsilon_n^{\min \{\s,\underline{\s} + 1\} - j}) + a_n^{j - \q} O_{\P}[(a_n \epsilon_n)^{-(j + 1/2)} + (a_n \epsilon_n)^{-j}],
\end{align*}
where the second equality uses $\epsilon_n \to 0$ and the defining property of $\{\lambda_j^{\mathtt{BR}}(k) : k = 1,\dots,\underline{\s}\}$.

The second part of the lemma follows from the fact that if
\begin{equation*}
    n \epsilon_n^{(1 + 2 \bar \q) \min(\underline{\s},\s - 1)/(\bar \q - 1)} \to 0 \quad \text{and} \quad n \epsilon_n^{1 + 2 \bar \q} \to \infty,
\end{equation*}
then
\begin{equation*}
    a_n^{\q - j} \epsilon_n^{\min\{\s,\underline{\s} + 1\} - j} \to 0 \quad \text{and} \quad a_n \epsilon_n \to \infty.
\end{equation*}
\qed 

\subsection{Higher-order expansion of the bias-reduced estimator }\label{Appendix Section: higher-order expansion}

In addition to the assumptions of Lemma \ref{M-[Lemma] Mean function estimation unknown q}, suppose that $\widehat{R}_{\x,n}(1;\eta_n) = O_{\P}(a_n^{-1/2})$ for $a_n^{-1}\eta_n^{-1}=O(1)$ and that, for some $\delta > 0$, $\theta_0$ is $(\underline{\s} + 1)$-times continuously differentiable and $\Phi_0$ is $(\underline{\s} + 2)$-times continuously differentiable on $I_{\x}^{\delta}$. Then, the first term in the stochastic expansion of $\widetilde{\mathcal{D}}_{j,n}^{\mathtt{BR}}(\x)$ satisfies
\begin{align*}
    \epsilon_n^{-(j + 1)} &\sum_{k = 1}^{\underline{\s} + 1}\lambda_j^{\mathtt{BR}}(k) [\Upsilon_0(\x + c_k \epsilon_n) - \Upsilon_0(\x) ] - \frac{\partial^{j + 1} \Upsilon_0(\x)}{(j + 1)!} \\
    &= \epsilon_n^{\underline{\s} + 1 - j} \frac{\partial^{\underline{\s} + 2} \Upsilon_0(\x)}{(\underline{\s} + 2)!} \sum_{k = 1}^{\underline{\s} + 1} \lambda_j^{\mathtt{BR}}(k) c_k^{\underline{\s} + 2} + o(\epsilon_n^{\underline{\s} + 1 - j}) =\epsilon_n^{\underline{\s} + 1 - j} \mathsf{B}_{j}^{\mathtt{BR}}(\x) + o(\epsilon_n^{\underline{\s} + 1 - j}).
\end{align*}
Also, the approximate variance of
\begin{equation*}
    \epsilon_n^{-(j + 1/2)} n^{-1/2} \sum_{k = 1}^{\underline{\s} + 1}\lambda_j^{\mathtt{BR}}(k) \widehat{G}_{\x,n}(c_k;\epsilon_n)
\end{equation*}
is
\begin{equation*}
    \frac{1}{n \epsilon_n^{1 + 2 j}} \sum_{k = 1}^{\underline{\s} + 1}\sum_{l = 1}^{\underline{\s} + 1} \lambda_{j}^{\mathtt{BR}}(k) \lambda_{j}^{\mathtt{BR}}(l) \mathcal{C}_\x(c_k,c_l) = \frac{1}{n \epsilon_n^{1 + 2 j}} \mathsf{V}_{j}^{\mathtt{BR}}(\x).
\end{equation*}
Finally, the third term in the stochastic expansion of $\widetilde{\mathcal{D}}_{j,n}^{\mathtt{BR}}(\x)$ is asymptotically negligible under the condition that $\widehat{R}_{\x,n}(1;\eta_n) = O_{\P}(a_n^{-1/2})$ for $a_n^{-1}\eta_n^{-1}=O(1)$, while the fourth term exhibits only a higher-order dependence on $\epsilon_n$ (relative to the dependence exhibited by the first two terms).


\subsection{Proof of Lemma \ref{M-[Lemma] Bootstrapping}}

We verify that Assumptions \ref{M-Assumption A} and \ref{M-Assumption E} imply Assumptions \ref{M-Assumption B: weak convergence}-\ref{M-Assumption B: uniform convergence Phi}. Define
\[\bar{\Gamma}_n^*(x) = \frac{1}{n} \sum_{i = 1}^nW_{i,n} \gamma_0(x;\mathbf{Z}_i) \qquad \text{and} \qquad \bar{\Phi}_n^*(x) = \frac{1}{n} \sum_{i = 1}^nW_{i,n} \phi_0(x;\mathbf{Z}_i).\]
\subsubsection*{Verifying Assumption  \ref{M-Assumption B: weak convergence}}

\paragraph*{Non-bootstrap weak convergence}
We first prove $\widehat{G}_{\x,n}^{\q}\rightsquigarrow \mathcal{G}_{\x}$. 
By Assumption \ref{M-Assumption E: Gamma hat}-\ref{M-Assumption E: Phi hat}, 
\begin{equation*}
    \sqrt{n a_n} \sup_{| v|\leq K}\left|\widehat{\Gamma}_n(\x+va_n^{-1}) -\widehat{\Gamma}_n(\x) - \bar{\Gamma}_n(\x+va_n^{-1}) + \bar{\Gamma}_n(\x)\right| =  o_{\P}(1), 
\end{equation*}
\begin{equation*}
    \sqrt{n a_n} \sup_{| v|\leq K}\left|\widehat{\Phi}_n(\x+va_n^{-1}) -\widehat{\Phi}_n(\x) - \bar{\Phi}_n(\x+va_n^{-1}) + \bar{\Phi}_n(\x)\right| = o_{\P}(1), 
\end{equation*}
for each $K>0$, and thus, defining $\bar{\psi}_{\x,n}(v;\mathbf{Z}_i) = \sqrt{a_n}\psi_{\x}(va_n^{-1};\mathbf{Z}_i)$,
\begin{equation*}
    \sup_{| v|\leq K}\left|\widehat{G}_{\x,n}^{\q}(v) - \frac{1}{\sqrt{n}} \sum_{i=1}^n \left\{\bar{\psi}_{\x,n}(v;\mathbf{Z}_i) - \E[\bar{\psi}_{\x,n}(v;\mathbf{Z})]  \right\}\right|  = o_{\P}(1).
\end{equation*}
Next, we prove weak convergence of the empirical process indexed by $\{\bar{\psi}_{\x,n}(v;\cdot):|v|\leq K\}$ by verifying finite-dimensional weak convergence and stochastic equicontinuity. 

Letting $\eta_n = Ka_n^{-1}$,
\begin{equation*}
    n^{-1}\sup_{|v|\leq K} \E[|\bar{\psi}_{\x,n}(v;\mathbf{Z})|^4] \leq (1+|\theta_0(\x)|)^4 n^{-1} a_n^2 \E[ \bar{D}_{\gamma}^{\eta_n} (\mathbf{Z})^4 +   \bar{D}_{\phi}^{\eta_n} (\mathbf{Z})^4] = o(1).
\end{equation*}
Also, convergence of the covariance kernel is imposed in Assumption \ref{M-Assumption E: covariance kernel}. Thus, the Lyapunov central limit theorem implies finite-dimensional weak convergence. 

For stochastic equicontinuity, proceeding as in \citet[Lemma 4.6]{Kim-Pollard_1990_AoS} and using $a_n\E[\bar{D}_{\gamma}^{\eta_n} (\mathbf{Z})^2 +   \bar{D}_{\phi}^{\eta_n} (\mathbf{Z})^2]=O(1)$, it suffices to show that
\begin{equation*}
    \sup_{\substack{|v-v'|\leq \epsilon_n\\ |v|\lor |v'|\leq K}} \frac{1}{n}\sum_{i=1}^n |\bar{\psi}_{\x,n}(v;\mathbf{Z}_i)-\bar{\psi}_{\x,n}(v';\mathbf{Z}_i)|^2=o_{\P}(1)
\end{equation*}
for any $\epsilon_n=o(1)$. For any $M>0$, using $|\bar{\psi}_{\x,n}(v;\mathbf{z})|\leq \sqrt{a_n}(1+\theta_0(\x))(\bar{D}_{\gamma}^{\eta_n} (\mathbf{z}) + \bar{D}_{\phi}^{\eta_n} (\mathbf{z}) )$
\begin{align*}
    &\sup_{\substack{|v-v'|\leq \epsilon_n\\ |v|\lor |v'|\leq K}} \frac{1}{n}\sum_{i=1}^n |\bar{\psi}_{\x,n}(v;\mathbf{Z}_i)-\bar{\psi}_{\x,n}(v';\mathbf{Z}_i)|^2 \\
    & \leq 4(1+|\theta_0(\x)|)^2\frac{a_n}{n}\sum_{i=1}^n  [\bar{D}_{\gamma}^{\eta_n} (\mathbf{Z}_i) + \bar{D}_{\phi}^{\eta_n} (\mathbf{Z}_i)]^2 \1(\bar{D}_{\gamma}^{\eta_n} (\mathbf{Z}_i) +   \bar{D}_{\phi}^{\eta_n} (\mathbf{Z}_i)> M) \\
    & + \bar{M} \sqrt{a_n} \sup_{\substack{|v-v'|\leq \epsilon_n\\ |v|\lor |v'|\leq K}} \E \left[|\bar{\psi}_{\x,n}(v;\mathbf{Z}_i)-\bar{\psi}_{\x,n}(v';\mathbf{Z}_i)| \right] \\
    & + \bar{M}\sqrt{a_n} \sup_{\substack{|v-v'|\leq \epsilon_n\\ |v|\lor |v'|\leq K}}\frac{1}{n}\sum_{i=1}^n  \left\{|\bar{\psi}_{\x,n}(v;\mathbf{Z}_i) - \bar{\psi}_{\x,n}(v';\mathbf{Z}_i)| -\E \left[|\bar{\psi}_{\x,n}(v;\mathbf{Z}_i) - \bar{\psi}_{\x,n}(v';\mathbf{Z}_i)| \right]\right\}
\end{align*}
where $\bar{M}=2M(1+|\theta_0(\x)|)$, the expectation of the first term after the inequality can be made arbitrarily small by making $M$ large using $a_n \E[ \bar{D}_{\gamma}^{\eta_n} (\mathbf{Z})^4 + \bar{D}_{\phi}^{\eta_n} (\mathbf{Z})^4] = O(1)$, the second term is $o(1)$ by Assumption \ref{M-Assumption E: covariance kernel}, and the third term is $O_{\P}( \sqrt{a_n/n} )$ by Theorem 4.2 of \cite{Pollard_1989_SS}. 

\paragraph*{Bootstrap weak convergence}
We next prove $\widehat{G}_{\x,n}^{\q}\rightsquigarrow_{\P}\mathcal{G}_{\x}$. As shown below, we have
\begin{align}\nonumber
    &\sqrt{n a_n}\sup_{|v|\leq K}|\widehat{\Gamma}_{n}^{* }(\x+va_n^{-1})-\widehat{\Gamma}_{n}^{* }(\x)-\bar{\Gamma}_{n}^{* }(\x+va_n^{-1})+\bar{\Gamma}_{n}^{* }(\x)|=o_{\P}(1),\\
    \label{Appendix eq: high-level bootstrap approximation}
    &\sqrt{n a_n}\sup_{|v|\leq K}|\widehat{\Phi}_{n}^{* }(\x+va_n^{-1})-\widehat{\Phi}_{n}^{* }(\x)-\bar{\Phi}_{n}^{* }(\x+va_n^{-1})+\bar{\Phi}_{n}^{* }(\x)|=o_{\P}(1),
\end{align}
for each $K>0$. Therefore, using $\widehat{\theta}_n(\x)\to_{\P}\theta_0(\x)$ and the fact that, uniformly over $| v| \leq K$,
\begin{equation*}
    \sqrt{\frac{a_n}{n}} \sum_{i=1}^nW_{i,n} \{\phi_0(\x+va_n^{-1};\mathbf{Z}_i) - \phi_0(\x;\mathbf{Z}_i) - \frac{1}{n}\sum_{j=1}^n[\phi_0(\x+va_n^{-1};\mathbf{Z}_j) - \phi_0(\x;\mathbf{Z}_j)]\} = O_{\P}(1),
\end{equation*}
we obtain
\begin{equation*}
    \sup_{| v|\leq K}\left|\widehat{G}_{\x,n}^{\q,*}(v)-\frac{1}{\sqrt{n}} \sum_{i=1}^n W_{i,n} \widehat{\psi}_{\x,n}(v;\mathbf{Z}_i)\right| = o_{\P}(1),
\end{equation*}
where $\widehat{\psi}_{\x,n}(v;\mathbf{Z}) = \bar{\psi}_{\x,n}(v;\mathbf{Z}) - n^{-1}\sum_{j=1}^n \bar{\psi}_{\x,n}(v;\mathbf{Z}_j)$.

To prove finite-dimensional weak convergence, we apply Lemma 3.6.15 of \cite{vanderVaart-Wellner_1996_Book}. Assumption \ref{M-Assumption E: bootstrap weights} implies that $n^{-1}\sum_{i=1}^n(W_{i,n}-1)^2\to_{\P}1$ and $n^{-1}\max_{1\leq i \leq n}W_{i,n}^2=o_{\P}(1)$. Since
\begin{align*}
    \frac{1}{n}\sum_{i=1}^n \widehat{\psi}_{\x,n}(v;\mathbf{Z}_i)\widehat{\psi}_{\x,n}(u;\mathbf{Z}_i)
    &= \frac{1}{n}\sum_{i=1}^n \bar{\psi}_{\x,n}(v;\mathbf{Z}_i)\bar{\psi}_{\x,n}(u;\mathbf{Z}_i) \\
    &\qquad - \frac{1}{n}\sum_{i=1}^n \bar{\psi}_{\x,n}(v;\mathbf{Z}_i)\frac{1}{n}\sum_{i=1}^n \bar{\psi}_{\x,n}(u;\mathbf{Z}_i)
\end{align*}
and $\sup_{|v|\leq\eta } |\psi_{\x}(v;\mathbf{Z})|\leq \bar{D}_{\gamma}^{\eta}(\mathbf{Z}) + |\theta_0(\x)|\bar{D}_{\phi}^{\eta}(\mathbf{Z})$, for any $v,u\in\R$,
\begin{equation*}
    \frac{1}{n}\sum_{i=1}^n \widehat{\psi}_{\x,n}(v;\mathbf{Z}_i)\widehat{\psi}_{\x,n}(u;\mathbf{Z}_i) - \E[\bar{\psi}_{\x,n}(v;\mathbf{Z})\bar{\psi}_{\x,n}(u;\mathbf{Z})] = o_{\P}(1).
\end{equation*}
Also, $n^{-1}\sum_{i=1}^n \widehat{\psi}_{\x,n}^4(v;\mathbf{Z}_i)=O_{\P}(1)$, verifying the hypothesis of the lemma.

For stochastic equicontinuity, let $\epsilon_n=o(1)$ and $\eta_n=K a_n^{-1}$. Lemma 3.6.7 of \cite{vanderVaart-Wellner_1996_Book} implies that for any $n_0 \in \{1,\dots, n\}$, there is a fixed constant $C>0$ such that
\begin{align*}
    &\E\left[ \left. \sup_{\substack{|v-v'|\leq \epsilon_n\\ |v|\lor |v'|\leq K}} \left| \frac{1}{\sqrt{n}} \sum_{i=1}^n W_{i,n}\left[  \widehat{\psi}_{\x,n}(v;\mathbf{Z}_i)-\widehat{\psi}_{\x,n}(v';\mathbf{Z}_i)\right] \right| \right| \{\mathbf{Z}_i\}_{i=1}^n\right] \\
    &\leq C \frac{\sqrt{a_n}}{n}\sum_{i=1}^n \left[\bar{D}_{\gamma}^{\eta_n}(\mathbf{Z}_i) + \bar{D}_{\phi}^{\eta_n}(\mathbf{Z}_i)\right] (n_0-1) \frac{1}{\sqrt{n}}\E\max_{1\leq i \leq n} |W_{i,n}| \\
    & + C \max_{n_0\leq k \leq n} \E\left[ \left. \sup_{\substack{|v-v'|\leq \epsilon_n\\ |v|\lor |v'|\leq K}} \left| \frac{1}{\sqrt{k}} \sum_{i=n_0}^k \left[ \widehat{\psi}_{\x,n}(v;\mathbf{Z}_{R_i}) -\widehat{\psi}_{\x,n}(v';\mathbf{Z}_{R_i}) \right] \right| \right| \{\mathbf{Z}_i\}_{i=1}^n \right]
\end{align*}
where $(R_1,\dots,R_n)$ is uniformly distributed on the set of all permutations of $\{1,\dots,n\}$, independent of $\{\mathbf{Z}_i\}_{i=1}^n$. Choose $n_0$ such that $n^{1/2-1/\mathfrak{r}}/n_0\to \infty$ and $n_0/a_n\to \infty$ (which is possible by $\mathfrak{r}>(4\q+2)/(2\q-1)$), and the first term after the inequality in the above display is $o_{\P}(1)$. For the second term, following the argument of \citet[Theorem 3.6.13]{vanderVaart-Wellner_1996_Book}, it suffices to bound
\begin{equation*}
    \max_{n_0\leq k \leq n} \E_n^* \left[ \sup_{\substack{|v-v'|\leq \epsilon_n\\ |v|\lor |v'|\leq K}} \left| \frac{1}{\sqrt{k}} \sum_{i=1}^k \left[ \widehat{\psi}_{\x,n}(v;\mathbf{Z}_i^*)-\widehat{\psi}_{\x,n}(v';\mathbf{Z}_i^*)\right] \right| \right], 
\end{equation*}
where $\{\mathbf{Z}_i^*\}_{i=1}^k$ denotes a random sample from the empirical cdf and $\E_n^*$ is the expectation under this empirical bootstrap law. Following the argument of \citet[Lemma 4.6]{Kim-Pollard_1990_AoS}, it suffices to show that
\begin{equation*}
    \max_{n_0\leq k \leq n} \E_n^* \left[ \sup_{\substack{|v-v'|\leq \epsilon_n\\ |v|\lor |v'|\leq K}} \frac{1}{k} \sum_{i=1}^k \left| \widehat{\psi}_{\x,n}(v;\mathbf{Z}_i^*)-\widehat{\psi}_{\x,n}(v';\mathbf{Z}_i^*)\right|^2 \right] = o_{\P}(1).
\end{equation*} 
For $k\in \{n_0,\dots,n\}$ and $M>0$,
\begin{align*}
    &\E_n^* \left[ \sup_{\substack{|v-v'|\leq \epsilon_n\\ |v|\lor |v'|\leq K}} \frac{1}{k} \sum_{i=1}^k \left| \widehat{\psi}_{\x,n}(v;\mathbf{Z}_i^*)-\widehat{\psi}_{\x,n}(v';\mathbf{Z}_i^*) \right| ^2 \right] \\
    &\leq 2\E_n^* \left[ \sup_{\substack{|v-v'|\leq \epsilon_n\\ |v|\lor |v'|\leq K}} \frac{1}{k} \sum_{i=1}^k \left| \bar{\psi}_{\x,n}(v;\mathbf{Z}_i^*)-\bar{\psi}_{\x,n}(v';\mathbf{Z}_i^*) \right| ^2 \right] \\
    &\quad + 2 \sup_{\substack{|v-v'|\leq \epsilon_n\\ |v|\lor |v'|\leq K}} \left\vert\frac{1}{n}\sum_{i=1}^n \bar{\psi}_{\x,n}(v;\mathbf{Z}_i)-\bar{\psi}_{\x,n}(v';\mathbf{Z}_i)\right\vert^2
\end{align*}
where the second term after the inequality is shown to be $o_{\P}(1)$ in the non-bootstrap case above. For the first term after the inequality,
\begin{align*}
    &\E_n^* \left[ \sup_{\substack{|v-v'|\leq \epsilon_n\\ |v|\lor |v'|\leq K}} \frac{1}{k} \sum_{i=1}^k \left| \bar{\psi}_{\x,n}(v;\mathbf{Z}_i^*)-\bar{\psi}_{\x,n}(v';\mathbf{Z}_i^*) \right| ^2 \right]\\
    &\leq 4(1+|\theta_0(\x)|)^2 \frac{a_n}{n}\sum_{i=1}^n \left[ \bar{D}_{\gamma}^{\eta_n}(\mathbf{Z}_i) + \bar{D}_{\phi}^{\eta_n}(\mathbf{Z}_i)\right]^2 \1(\bar{D}_{\gamma}^{\eta_n}(\mathbf{Z}_i) + \bar{D}_{\phi}^{\eta_n}(\mathbf{Z}_i) > M) \\
    &+ \bar{M} \sqrt{a_n} \sup_{\substack{|v-v'|\leq \epsilon_n\\ |v|\lor |v'|\leq K}} \frac{1}{n}\sum_{i=1}^n \left| \bar{\psi}_{\x,n}(v;\mathbf{Z}_i)-\bar{\psi}_{\x,n}(v';\mathbf{Z}_i) \right| \\
    & + \bar{M} \sqrt{a_n} \E_n^* \left[ \sup_{\substack{|v-v'|\leq \epsilon_n\\ |v|\lor |v'|\leq K}} \frac{1}{k} \sum_{i=1}^k \left| \bar{\psi}_{\x,n}(v;\mathbf{Z}_i^*)-\bar{\psi}_{\x,n}(v';\mathbf{Z}_i^*) \right| - \E_n^* \left[ \left| \bar{\psi}_{\x,n}(v;\mathbf{Z}^*) - \bar{\psi}_{\x,n}(v';\mathbf{Z}^*) \right| \right] \right]
\end{align*}
where $\bar{M}=2M(1+|\theta_0(\x)|)$.
The first term after the inequality does not depend on $k$ and its expectation can be made arbitrarily small by taking $M$ sufficiently large. The second term is independent of $k$ and we can handle this term by adding and subtracting the expectation inside the summation. For the third term, applying Theorem 4.2 of \cite{Pollard_1989_SS} again, it is bounded by a constant multiple of
\begin{equation*}
    \sqrt{\frac{a_n}{k}} \sqrt{\frac{a_n}{n} \sum_{i=1}^n \left[ \bar{D}_{\gamma}^{\eta_n}(\mathbf{Z}_i)^2 + \bar{D}_{\phi}^{\eta_n}(\mathbf{Z}_i)^2 \right] } = O_{\P}\left( \sqrt{a_n/k} \right),
\end{equation*}
which is $o_{\P}(1)$ by the choice of $n_0$.

\paragraph*{Verifying \eqref{Appendix eq: high-level bootstrap approximation}}
We focus on the first display.
By adding and subtracting the bootstrap means,
\begin{align*}
    &\widehat{\Gamma}_n^*(\x + va_n^{-1}) -\widehat{\Gamma}_n^*(\x) - \bar{\Gamma}_n^*(\x+va_n^{-1})+ \bar{\Gamma}_n^*(\x)\\
    &= \frac{1}{n}\sum_{i=1}^n W_{i,n} \Breve{\gamma}_n(v;\mathbf{Z}) + \left[ \check{\Gamma}_n(\x+va_n^{-1}) - \check{\Gamma}_n(\x) - \bar{\Gamma}_n(\x+va_n^{-1}) + \bar{\Gamma}_n(\x) \right]
\end{align*}
where 
\begin{align*}
    \Breve{\gamma}_n(v;\mathbf{Z}) &= \widehat{\gamma}_n(\x+va_n^{-1};\mathbf{Z}_i)-\widehat{\gamma}_n(\x;\mathbf{Z}_i) - \gamma_0(\x+va_n^{-1};\mathbf{Z}_i)+\gamma_0(\x;\mathbf{Z}_i)\\
    &\qquad - \check{\Gamma}_n(\x+va_n^{-1}) + \check{\Gamma}_n(\x) - \bar{\Gamma}_n(\x+va_n^{-1}) + \bar{\Gamma}_n(\x).
\end{align*}
By Assumption \ref{M-Assumption E: Gamma hat},
\begin{equation*}
    \sqrt{n a_n} \sup_{|v|\leq K} | \check{\Gamma}_n(\x+va_n^{-1}) - \check{\Gamma}_n(\x) - \bar{\Gamma}_n(\x+va_n^{-1}) + \bar{\Gamma}_n(\x) | = o_{\P}(1).
\end{equation*}
Identical to above, Lemma 3.6.7 and the argument in Theorem 3.6.13 of \cite{vanderVaart-Wellner_1996_Book} imply that for some fixed $C>0$,
\begin{align*}
    &\sqrt{n a_n} \E\left[ \left. \sup_{ |v|\leq K}  \left| \frac{1}{n}\sum_{i=1}^n W_{i,n} \Breve{\gamma}_n(v;\mathbf{Z}_i) \right| \right| \{\mathbf{Z}_i\}_{i=1}^n \right] \\
    \label{Appendix eq: high-level bootstrap approximation remainder 1}
    &\leq C  \frac{\sqrt{a_n}}{n} \sum_{i=1}^n \sup_{ |v|\leq K} |\Breve{\gamma}_n(v;\mathbf{Z}_i) | \frac{n_0 n^{\mathfrak{r}}}{\sqrt{n}} +  C \sqrt{a_n}\max_{n_0\leq k\leq n} \E_n^* \left[ \sup_{ |v|\leq K} \left| \frac{1}{\sqrt{k}} \sum_{i=1}^k \Breve{\gamma}_n(v;\mathbf{Z}_i^*) \right| \right].
\end{align*}
By Assumption \ref{M-Assumption E: Gamma hat},
\begin{align*}
    & \frac{\sqrt{a_n}}{n} \sum_{i=1}^n \sup_{ |v|\leq K} |\Breve{\gamma}_n(v;\mathbf{Z}_i) |\\
    &\leq \frac{\sqrt{a_n}}{n} \sum_{i=1}^n \sup_{ |v|\leq K} |\widehat{\gamma}_n(\x+va_n^{-1};\mathbf{Z}_i)-\widehat{\gamma}_n(\x;\mathbf{Z}_i) -\gamma_0(\x+va_n^{-1};\mathbf{Z}_i) + \gamma_0(\x;\mathbf{Z}_i) | \\
    &\quad + \sqrt{a_n} \sup_{ |v|\leq K} \left|\check{\Gamma}_n(\x+va_n^{-1}) - \check{\Gamma}_n(\x) - \bar{\Gamma}_n(\x+va_n^{-1}) + \bar{\Gamma}_n(\x) \right| = o_{\P}(1).
\end{align*}
Also, Corollary 4.3 of \cite{Pollard_1989_SS} implies that for some fixed $C>0$,
\begin{align*}
    &\E_n^*\left[ \sup_{ |v|\leq K} \left|\frac{1}{\sqrt{k}}\sum_{i=1}^k \Breve{\gamma}_n(v;\mathbf{Z}_i^*)  \right| \right]^2 \\
    &\leq C \frac{1}{n}\sum_{i=1}^n\sup_{ |v|\leq K} |\widehat{\gamma}_n(\x+va_n^{-1};\mathbf{Z}_i)-\widehat{\gamma}_n(\x;\mathbf{Z}_i) -\gamma_0(\x+va_n^{-1};\mathbf{Z}_i) + \gamma_0(\x;\mathbf{Z}_i) |^2,
\end{align*}
which is $o_{\P}(a_n^{-1})$ by  Assumption  \ref{M-Assumption E: Gamma hat}.

\subsubsection*{Verifying Assumption \ref{M-Assumption B: moment bound}}
Defining
\begin{align*}
    \bar{G}_{\x,n}^{\q}(v)
    &= \sqrt{n a_n} \left[ \bar{\Gamma}_n(\x +v a_n^{-1} ) - \bar{\Gamma}_n(\x) - \Gamma_0(\x+v a_n^{-1}) + \Gamma_0(\x) \right]  \\
    &\quad -\theta_0(\x) \sqrt{n a_n} \left[ \bar{\Phi}_n(\x +v a_n^{-1} ) - \bar{\Phi}_n(\x) - \Phi_0(\x+v a_n^{-1}) + \Phi_0(\x) \right],
\end{align*}
we have
\begin{align*}
    \widehat{G}_{\x,n}^{\q}(v) - \bar{G}_{\x,n}^{\q}(v)
    &= \sqrt{n a_n} \left[ \widehat{\Gamma}_n(\x+va_n^{-1}) -\widehat{\Gamma}_n(\x) - \bar{\Gamma}_n(\x +va_n^{-1} ) + \bar{\Gamma}_n(\x) \right] \\
    &\quad - \theta_0(\x) \sqrt{n a_n} \left[ \widehat{\Phi}_n(\x+va_n^{-1}) -\widehat{\Phi}_n(\x) - \bar{\Phi}_n(\x +va_n^{-1} ) + \bar{\Phi}_n(\x) \right].
\end{align*}
Assumptions \ref{M-Assumption E: Gamma hat} and \ref{M-Assumption E: Phi hat} imply that for $V \in [1,a_n \delta]$,
\begin{equation*}
    \sqrt{n a_n}\sup_{|v|\in [V,2V]} | \widehat{\Gamma}_n(\x+v a_n^{-1}) -\widehat{\Gamma}_n(\x) - \bar{\Gamma}_n(\x +v a_n^{-1}) + \bar{\Gamma}_n(\x) | \leq (2V)^{\beta} a_n^{-\beta}B_n+A_n
\end{equation*}
and
\begin{equation*}
    \sqrt{n a_n}\sup_{|v|\in [V,2V]} | \widehat{\Phi}_n(\x+v a_n^{-1}) -\widehat{\Phi}_n(\x) - \bar{\Phi}_n(\x +v a_n^{-1}) + \bar{\Phi}_n(\x) | \leq (2V)^{\beta} a_n^{-\beta}B_n+A_n,
\end{equation*}
where $\beta = \beta_\gamma \lor \beta_{\phi} , A_n = A_{\gamma,n} \lor A_{\phi,n} = o_{\P}(1)$, and $B_n = B_{\gamma,n} \lor B_{\phi,n} = o_{\P}(a_n^{\beta})$. As a consequence, there exists $\eta_n' = o(1)$ such that
\begin{equation*}
    \lim_{n\to\infty} \P\left[\bigcap_{V \in [1,a_n\delta]} \left\{ \sup_{| v|\in [V,2V]}|\widehat{G}_{\x,n}^{\q}(v) - \bar{G}_{\x,n}^{\q}(v) | \leq [V^{\beta}+1]\eta_n' \right\}\right] = 1;
\end{equation*}
we take the event in the display to be $\mathcal{A}_n$. Also,
\begin{equation*}
    \limsup_{\eta \downarrow 0}\E[\bar{D}_{\gamma}^{\eta}(\mathbf{Z})^2/\eta] + \E[\bar{D}_{\phi}^{\eta}(\mathbf{Z})^2/\eta] < \infty
\end{equation*}
and, using Corollary 4.3 of \cite{Pollard_1989_SS},
\begin{equation*}
    \E \left[ \sup_{| v|\leq V}|\bar{G}_{\x,n}^{\q}(v)| \right] \leq C \sqrt{V} \quad \text{for } V \in (0,a_n\delta].
\end{equation*}
Therefore,
\begin{equation*}
    \E \left[ \sup_{| v| \in [V,2V]}| \widehat{G}_{\x,n}^{\q}(v)|\1_{\mathcal{A}_n} \right] \leq C \sqrt{V} + [V^{\beta}+1]\eta_n',
\end{equation*}
implying in particular that
\begin{equation*}
    \sup_{V\in[1,a_n\delta]}\E\left[V^{-\beta} \sup_{|v|\in [V,2V] } | \widehat{G}_{\x,n}^{\q}(v)| \1_{\mathcal{A}_n}\right] = O(1).
\end{equation*}

For the bootstrap counterpart, defining
\begin{align*}
    \bar{G}_{\x,n}^{\q,*}(v) &= \sqrt{n a_n}\left[ \bar{\Gamma}_n^*(\x +v a_n^{-1} ) - \bar{\Gamma}_n^*(\x) - \bar{\Gamma}_n(\x +v a_n^{-1} ) + \bar{\Gamma}_n(\x) \right]  \\
    &\quad  -\theta_0(\x) \sqrt{n a_n} \left[ \bar{\Phi}_n^*(\x +v a_n^{-1} ) - \bar{\Phi}_n^*(\x) - \bar{\Phi}_n(\x +v a_n^{-1} ) + \bar{\Phi}_n(\x) \right],
\end{align*}
we have
\begin{align}
\label{Appendix eq: verifying bootstrap rate of convergence}
    &\widehat{G}_{\x,n}^{\q,*}(v) - \bar{G}_{\x,n}^{\q,*}(v) \\
    \nonumber
    &= \sqrt{n a_n} \left[ \widehat{\Gamma}_n^*(\x+va_n^{-1}) -\widehat{\Gamma}_n^*(\x) - \bar{\Gamma}_n^*(\x +va_n^{-1}) + \bar{\Gamma}_n^*(\x) \right] \\
    \nonumber
    & - \widehat{\theta}_n(\x)\sqrt{n a_n} \left[ \widehat{\Phi}_n^*(\x+va_n^{-1}) -\widehat{\Phi}_n^*(\x) - \bar{\Phi}_n^*(\x +va_n^{-1}) + \bar{\Phi}_n^*(\x) \right] \\
    \nonumber
    & + \sqrt{n a_n} \left[ \bar{\Gamma}_n(\x +va_n^{-1}) - \bar{\Gamma}_n(\x) - \widehat{\Gamma}_n(\x+va_n^{-1}) + \widehat{\Gamma}_n(\x) \right] \\
    \nonumber
    & - \widehat{\theta}_n(\x)\sqrt{n a_n} \left[ \bar{\Phi}_n(\x +va_n^{-1}) - \bar{\Phi}_n(\x) - \widehat{\Phi}_n(\x+va_n^{-1}) + \widehat{\Phi}_n(\x) \right] \\
    \nonumber
    & +\sqrt{a_n}[\widehat{\theta}_n(\x) - \theta_0(\x)] \sqrt{n} \left[ \bar{\Phi}_n^*(\x + va_n^{-1}) - \bar{\Phi}_n^*(\x) - \bar{\Phi}_n(\x +va_n^{-1}) + \bar{\Phi}_n(\x) \right],
\end{align}
where bounds for the third and fourth terms were derived in the non-bootstrap case and where the last term is $o_{\P}(1)$ uniformly over $| v| \leq a_n \delta$ as $\sqrt{a_n}[\widehat{\theta}_n(\x) - \theta_0(\x)]=o_{\P}(1)$ and $\sqrt{n} \sup_{| v|\leq \delta}|\bar{\Phi}_n^*(\x + v) - \bar{\Phi}_n(\x +v ) |=O_{\P}(1)$. For the first term after the equality in \eqref{Appendix eq: verifying bootstrap rate of convergence},
\begin{align*}
    & \sqrt{na_n} \left[ \widehat{\Gamma}_n^*(\x+va_n^{-1}) -\widehat{\Gamma}_n^*(\x) - \bar{\Gamma}_n^*(\x +v a_n^{-1}) + \bar{\Gamma}_n^*(\x) \right] \\
    &= \frac{1}{\sqrt{n}}\sum_{i=1}^n W_{i,n}\Breve{\gamma}_n(va_n^{-1};\mathbf{Z}_i) + \sqrt{na_n} \left[ \check{\Gamma}_n(\x+va_n^{-1}) - \check{\Gamma}_n(\x) - \bar{\Gamma}_n(\x+va_n^{-1}) + \bar{\Gamma}_n(\x) \right]
\end{align*}
where, by Assumption \ref{M-Assumption E: Gamma hat},
\begin{equation*}
    \sqrt{n a_n}\sup_{|v|\leq V } | \check{\Gamma}_n(\x+va_n^{-1}) - \check{\Gamma}_n(\x) - \bar{\Gamma}_n(\x+va_n^{-1}) + \bar{\Gamma}_n(\x)| =[1+ V^{\beta}]o_{\P}(1)
\end{equation*}
uniformly over $V\in [1,a_n\delta]$ and where, applying Lemma 3.6.7 of \cite{vanderVaart-Wellner_1996_Book}, we have, for any $n_0\in \{1,\dots, n\}$,
\begin{align*}
    &\E\left[ \left. \sup_{| v|\leq \delta} \left| \frac{1}{\sqrt{n}}\sum_{i=1}^nW_{i,n}\Breve{\gamma}_n(v;\mathbf{Z}_i) \right| \right| \{\mathbf{Z}_i\}_{i=1}^n \right] \\
    \label{Appendix eq: verifying bootstrap moment bound}
    &\leq C \frac{1}{n}\sum_{i=1}^n \sup_{| v| \leq \delta} |\Breve{\gamma}_n(v;\mathbf{Z}_i)| \frac{n_0}{\sqrt{n}} \E \left[ \max_{1\leq i\leq n} | W_{i,n}| \right] + C \max_{n_0\leq k\leq n}\E_n^* \left[\sup_{| v|\leq \delta} \left| \frac{1}{\sqrt{k}} \sum_{i=1}^k \Breve{\gamma}_n(v;\mathbf{Z}_i^*) \right| \right],
\end{align*}
where $C>0$ is a fixed constant. Letting $n_0$ be a diverging sequence (dependent on $n$) such that $n_0 n^{\mathfrak{r}}\sqrt{a_n/n} =o(1)$, the first term on the majorant side is bounded by
\begin{equation*}
    C\frac{1}{n}\sum_{i=1}^n \sup_{| x-\x|\leq \delta} | \widehat{\gamma}_n(x;\mathbf{Z}_i) - \gamma_0(x;\mathbf{Z}_i) | + C\sup_{| v|\leq \delta } | \check{\Gamma}_n(\x+v) - \bar{\Gamma}_n(\x+v)  | =o_{\P}(1),
\end{equation*}
and, by Corollary 4.3 of \cite{Pollard_1989_SS}, the second term is bounded (up to a constant) by
\begin{equation*}
    \sqrt{a_n} \left| \frac{1}{n}\sum_{i=1}^n \sup_{x\in I_{\x}^{\delta} } | \widehat{\gamma}_n(x;\mathbf{Z}_i)- \gamma_0(x;\mathbf{Z}_i) |^2 \right| ^{1/2} 
\end{equation*}
which is $o_{\P}(1)$ by Assumption \ref{M-Assumption E: Gamma hat}.
As a consequence, there exists a sequence of random variables $A_n'=o_{\P}(1)$ such that for $V\in [1,a_n\delta]$,
\begin{equation*}
    \sqrt{n a_n}\sup_{|v|\leq [V,2V]} \left| \widehat{\Gamma}_n^*(\x+va_n^{-1}) -\widehat{\Gamma}_n^*(\x) - \bar{\Gamma}_n^*(\x +va_n^{-1}) + \bar{\Gamma}_n^*(\x) \right| \leq  [V^{\beta}+1]A_n'.
\end{equation*}
By identical arguments, an analogous bound holds for the second term after the inequality in \eqref{Appendix eq: verifying bootstrap rate of convergence}. Therefore, there exists $\eta_n'=o(1)$ and events $\mathcal{A}_n$ such that $\lim_{n\to\infty}\P[\mathcal{A}_n]=1$ and
\begin{equation*}
    \E\left[\sup_{|v|\in [V,2V]} \left| \widehat{G}_{\x,n}^{\q,*}(v) - \bar{G}_{\x,n}^{\q,*}(v) \right| \1_{\mathcal{A}_n} \right]\leq [V^{\beta}+1]\eta_n'
\end{equation*}
for any $V\in [1,a_n\delta]$. Finally, proceeding as in the proof of stochastic equicontinuity, we obtain the bound $\E[\sup_{| v|\in [V,2V] }|\bar{G}_{\x,n}^{\q,*}(v)|]\leq C\sqrt{V}$.

\subsubsection*{Verifying Assumptions  \ref{M-Assumption B: uniform consistency Gamma}-\ref{M-Assumption B: uniform convergence Phi}}
We have
\begin{equation*}
    \sup_{x\in I} \left| \widehat{\Gamma}_n(x)-\Gamma_0(x) \right| \leq \sup_{x\in I} \left| \widehat{\Gamma}_n(x)-\bar{\Gamma}_n(x) \right| + \sup_{x\in I} \left| \bar{\Gamma}_n(x)-\Gamma_0(x) \right|,
\end{equation*}
where the first term after the inequality is assumed to be $o_{\P}(1)$ and the second term is $o_{\P}(1)$ by standard arguments. Similarly, $\sup_{x\in I}|\widehat{\Phi}_n(x)-\Phi_0(x)|=o_{\P}(1)$. Also,
\begin{align*}
    \sup_{x\in I_{\x}^{\delta}} \left| \widehat{\Phi}_n(x) - \Phi_0(x) \right| 
    &\leq \sup_{x\in I_{\x}^{\delta}} \left| \widehat{\Phi}_n(x) - \widehat{\Phi}_n(\x) - \bar{\Phi}_n(x) +\bar{\Phi}_n(\x) \right| \\
    &\quad +\sup_{x\in I_{\x}^{\delta}} \left| \bar{\Phi}_n(x)-\Phi_0(x) \right| + \left| \widehat{\Phi}_n(\x) - \bar{\Phi}_n(\x) \right| 
\end{align*}
where the last two terms are $o_{\P}(a_n^{-1})$ and where Assumption \ref{M-Assumption E: Phi hat} implies
\begin{equation*}
    \sup_{x\in I_{\x}^{\delta}} \left| \widehat{\Phi}_n(x) - \widehat{\Phi}_n(\x) - \bar{\Phi}_n(x) +\bar{\Phi}_n(\x) \right| \leq (na_n)^{-1/2} [ A_{\phi,n}+ \delta^{\beta_{\phi}} B_{\phi,n}]  = o_{\P}(a_n^{-1}),
\end{equation*}
the last equality using $B_{\phi,n}=o_{\P}(a_n^{\beta_{\phi}})$ and $\beta_{\phi}\leq \q$.

Now we look at the bootstrap objects. For $\widehat{\Gamma}_n^*$,
\begin{equation*}
    \sup_{x\in I} \left| \widehat{\Gamma}_n^*(x) - \widehat{\Gamma}_n(x) \right| \leq \sup_{x\in I} \left| \widehat{\Gamma}_n^*(x) - \bar{\Gamma}_n^*(x) \right| + \sup_{x\in I} \left| \bar{\Gamma}_n^*(x) -\bar{\Gamma}_n(x) \right| + \sup_{x\in I} \left| \bar{\Gamma}_n(x) - \widehat{\Gamma}_n(x) \right|
\end{equation*}
where the last term is assumed to be $o_{\P}(1)$ and where
\begin{align*}
    \sup_{x\in I} \left| \widehat{\Gamma}_n^*(x) - \bar{\Gamma}_n^*(x) \right|
    &\leq \frac{1}{n}\sum_{i=1}^n |W_{i,n}| \sup_{x\in I} \left| \widehat{\gamma}_n(x;\mathbf{Z}_i) - \gamma_0(x;\mathbf{Z}_i) \right| \\
    &\leq \sqrt{\frac{1}{n}\sum_{i=1}^n |W_{i,n}|^2 \frac{1}{n}\sum_{i=1}^n  \sup_{x\in I} \left| \widehat{\gamma}_n(x;\mathbf{Z}_i) - \gamma_0(x;\mathbf{Z}_i) \right|^2 } = o_{\P}(1).
\end{align*}
For $\bar{\Gamma}_n^* -\bar{\Gamma}_n$, using Lemma 3.6.7 of \cite{vanderVaart-Wellner_1996_Book} and the same argument as for verifying Assumption \ref{M-Assumption B: weak convergence}, it suffices to show that
\begin{equation*}
    n^{-1/2} \max_{\lfloor\sqrt{n}\rfloor\leq k\leq n} \E_n^* \left[ \sup_{x\in I} \left| \frac{1}{\sqrt{k}} \sum_{i=1}^k \bar{\gamma}_n(x;\mathbf{Z}_i^*) \right| \right] = o_{\P}(1),
\end{equation*}
as can be done using Corollary 4.3 of \cite{Pollard_1989_SS}.

For $\widehat{\Phi}_n^*$, $\sup_{x\in I}|\widehat{\Phi}_n^*(x)-\widehat{\Phi}_n(x)|=o_{\P}(1)$ follows from the same argument as for $\widehat{\Gamma}_n^*$. For $a_n \sup_{x\in I_{\x}^{\delta}}|\widehat{\Phi}_n^*(x) -\widehat{\Phi}_n(x)|=o_{\P}(1)$,
\begin{align*}
    \sup_{x\in I_{\x}^{\delta}} \left| \widehat{\Phi}_n^*(x) -\widehat{\Phi}_n(x) \right| 
    &\leq \sup_{x\in I_{\x}^{\delta}} \left| \widehat{\Phi}_n^*(x) -\bar{\Phi}_n^*(x) \right| +\sup_{x\in I_{\x}^{\delta}} \left| \bar{\Phi}_n^*(x) - \bar{\Phi}_n(x) \right| \\
    &\quad + \sup_{x\in I_{\x}^{\delta}} \left| \bar{\Phi}_n(x) - \Phi_0(x) \right| + \sup_{x\in I_{\x}^{\delta}} \left| \widehat{\Phi}_n(x) - \Phi_0(x) \right|
\end{align*}
where the last term is $o_{\P}(a_n^{-1})$ as shown above and the second and third terms after the inequality are $O_{\P}(n^{-1/2})$ by standard arguments. For the remaining term,
\begin{align*}
    \sup_{x\in I_{\x}^{\delta}} \left| \widehat{\Phi}_n^*(x) -\bar{\Phi}_n^*(x) \right| 
    &\leq \sup_{x\in I_{\x}^{\delta}} \left| \frac{1}{n}\sum_{i=1}^n W_{i,n}\Breve{\phi}_n(x;\mathbf{Z}_i) \right| + \left| \check{\Phi}_n(\x) - \bar{\Phi}_n(\x) \right| \\
    &\quad + \sup_{x\in I_{\x}^{\delta}} \left| \check{\Phi}_n(x) - \check{\Phi}_n(\x) - \bar{\Phi}_n(x) + \bar{\Phi}_n(\x) \right|
\end{align*}
where $\Breve{\phi}_n(x;\mathbf{Z}) = \widehat{\phi}_n(x;\mathbf{Z})-\phi_0(x;\mathbf{Z}) - [\check{\Phi}_n(x)-\bar{\Phi}_n(x)]$. The last two terms are $o_{\P}(a_n^{-1})$ by Assumption \ref{M-Assumption E: Phi hat}. Using Lemma 3.6.7 of \cite{vanderVaart-Wellner_1996_Book} and the argument similar to above, the remaining term is $o_{\P}(a_n^{-1})$.


\subsection{Proof of Theorem \ref{M-[Theorem]: Bootstrap inconsistency}}
The proof is by contradiction and follows \cite{Kosorok_2008_BookCh}. We omit some details in cases where the arguments are almost identical to those for Theorem \ref{M-[Theorem] Main result} and Lemma \ref{M-[Lemma] Bootstrapping}.

Suppose that the bootstrap approximation is consistent; that is, suppose
\begin{equation*}
    r_n \left( \widehat{\theta}_n^*(\x) - \widehat{\theta}_n(\x) \right) \leadsto_{\P} \frac{1}{\partial \Phi_0(\x)} \partial_{-}\GCM_{\R}\{\mathcal{G}_{\x}+\mathcal{M}_{\x}^{\q}\}(0) = Y.
\end{equation*}
Then, by Theorem 2.2 of \cite{Kosorok_2008_BookCh}, we have 
\begin{equation}\label{eqBootFail:contradiction}
    r_n \left( \widehat{\theta}_n^*(\x) - \theta_0(\x) \right) \leadsto Y_1 + Y_2 =_d \sqrt{2} Y
\end{equation}
where $=_d$ denotes the distributional equality, $Y_1$ and $Y_2$ are independent copies of $Y$, and where the convergence in distribution is unconditional.

Using the switching lemma, $\P \left[r_n \left( \widehat{\theta}_n^*(\x) - \theta_0(\x) \right)  > t \right]$ equals
\begin{equation*}
    \P \left[ \max \argmax_{x\in \widehat{\Phi}_n^{*-} ([0,\widehat{u}_n^*]) } \left\{ [\theta_0(\x) + r_n^{-1}t]\widehat{\Phi}_n^*(x) - \LSC_{\widehat{\Phi}_n^{*}} \left( \widehat{\Gamma}_n^* \right) (x) \right\} < \widehat{\Phi}_n^{*-} \circ \widehat{\Phi}_n^*(\x) \right].
\end{equation*}
By the arguments used in the proof of Theorem \ref{M-[Theorem] Main result}, to characterize the limiting distribution of $r_n(\widehat{\theta}_n^*(\x) - \theta_0(\x))$, it suffices to look at
\begin{equation*}
    -\check{G}_{\x,n}^{\q,*}(v) - \check{M}_{\x,n}^{\q}(v) + t\widehat{L}_{\x,n}^{\q,*}(v)
\end{equation*}
where 
\begin{align*}
    \check{G}_{\x,n}^{\q,*}(v) &= \sqrt{na_n} \left[ \widehat{\Gamma}_n^*(\x+va_n^{-1})  - \widehat{\Gamma}_n^*(\x) - \check{\Gamma}_n (\x+va_n^{-1}) + \check{\Gamma}_n (\x) \right] \\
    &\quad - \theta_0(\x)\sqrt{na_n} \left[ \widehat{\Phi}_n^*(\x+va_n^{-1}) - \widehat{\Phi}_n^*(\x)-\check{\Phi}_n(\x+va_n^{-1}) + \check{\Phi}_n(\x) \right]
\end{align*}
and
\begin{equation*}
    \check{M}_{\x,n}^{\q}(v)=\sqrt{n a_n} \left[ \check{\Gamma}_n(\x+va_n^{-1}) -\check{\Gamma}_n(\x) - \theta_0(\x)\{ \check{\Phi}_n(\x+va_n^{-1}) -\check{\Phi}_n(\x)\} \right].
\end{equation*}
It can be shown that $\check{G}_{\x,n}^{\q,*} \leadsto_{\P} \mathcal{G}_{\x}, \widehat{L}_{\x,n}^{\q,*} \leadsto \mathcal{L}_{\x}$, and that $\check{M}_{\x,n}^{\q} \leadsto \mathcal{G}_{\x} + \mathcal{M}_{\x}^{\q}$. Thus, 
\begin{equation*}
    \P \left[ r_n \left( \widehat{\theta}_n^*(\x) - \theta_0(\x) \right) >t \right] \to \P \left[ \argmin_{v\in \R} \left\{ \mathcal{G}_{\x,1}(v) + \mathcal{G}_{\x,2}(v) + \mathcal{M}_{\x}^{\q}(v)  - t \partial\Phi_0(\x) v \right\} < 0 \right],
\end{equation*}
where $\mathcal{G}_{\x,1}$ and $\mathcal{G}_{\x,2}$ are independent copies of $\mathcal{G}_{\x}$.
Noting that $\mathcal{G}_{\x}(a v)=_d \sqrt{|a|} \mathcal{G}_{\x}(v)$ and using the change of variable $v=u2^{\frac{1}{2\q+1}}$, the limit distribution equals
\begin{align*}
    &\P \left[ 2^{\frac{1}{2\q+1}} \argmin_{u\in \R} \left\{   \mathcal{G}_{\x}(u) + \mathcal{M}_{\x}^{\q}(u) -  2^{-\frac{\q}{2\q+1}} t \partial\Phi_0(\x)u \right\} < 0 \right] \\
    &= \P \left[ 2^{\frac{\q}{2\q+1}} (\partial\Phi_0(\x))^{-1}  \partial_{-} \GCM_{\R} \{\mathcal{G}_{\x}+\mathcal{M}_{\x}^{\q}\}(0) >  t \right].
\end{align*} 
As a consequence,
\begin{equation*}
    r_n \left( \widehat{\theta}_n^*(\x) - \theta_0(\x) \right) \leadsto 2^{\frac{\q}{2\q+1}} (\partial \Phi_0(\x))^{-1} \partial_{-}\GCM_{\R}\{\mathcal{G}_{\x}+\mathcal{M}_{\x}^{\q}\}(0),
\end{equation*}
contradicting \eqref{eqBootFail:contradiction} because $2^{\frac{\q}{2\q+1}} \neq \sqrt{2}$.

In other words, the bootstrap estimator $\widehat{\theta}_n^*(\x)$ fails to approximate the limit distribution. \hfill$\qedsymbol$


\subsection{Remarks on verifying conditions in applications}\label{Appendix: section remark}

Below we verify the hypothesis of Theorem \ref{M-[Theorem] Main result} for various examples. For this purpose, it suffices to verify Assumptions \ref{M-Assumption A}, \ref{M-Assumption B: uhat}-\ref{M-Assumption B: jump of Phi hat}, and \ref{M-Assumption E} since Assumption \ref{M-Assumption E} implies \ref{M-Assumption B: weak convergence}-\ref{M-Assumption B: uniform convergence Phi} by Lemma \ref{M-[Lemma] Bootstrapping}.

When $\gamma_0$ is known, it is natural to take $\widehat{\Gamma}_n=\check{\Gamma}_n=\bar{\Gamma}_n$, in which case \ref{M-Assumption E: Gamma hat} reduces to the requirement that, for some $\rho_{\gamma}\in (0,2)$, 
\begin{equation}\label{Appendix Eq: simpler D gamma hat}
    \limsup_{\varepsilon\downarrow0}\frac{\log N_{U}(\varepsilon,\mathfrak{F}_{\gamma})}{\varepsilon^{-\rho_{\gamma}}}<\infty,\quad \E[\bar{F}_{\gamma}(\mathbf{Z})^2] <\infty,\quad \limsup_{\eta\downarrow0}\frac{\E[\bar{D}_{\gamma}^{\eta}(\mathbf{Z})^2+\bar{D}_{\gamma}^{\eta}(\mathbf{Z})^4]}{\eta}<\infty.
\end{equation} 
An identical remark applies to $\phi_0$ and \ref{M-Assumption E: Phi hat}.

\subsection{Proof of Corollary \ref{M-[Corollary] monotone density without censoring}}
Assumption \ref{M-Assumption A} and \ref{M-Assumption E: iid}-\ref{M-Assumption E: bootstrap weights} follow from the hypothesis. 
\newline 
\ref{M-Assumption E: Gamma hat}:
In this example, $\gamma_0(x;\mathbf{Z})=\1(X\leq x)$ is known, so it suffices to verify \eqref{Appendix Eq: simpler D gamma hat}. The uniform covering number of $\{\1(\cdot\leq x) :x\in\R\}$ grows linearly, and an envelope function can be taken to be $1$. For an envelope function of $\{\1(\cdot\leq x)-\1(\cdot\leq \x) : |x-\x|\leq \eta \}$, we can take $\1(-\eta +\x\leq \cdot\leq \x+ \eta )$ and the moment bound is satisfied as $\E[\1(-\eta +\x\leq X\leq \x+ \eta )]\leq C \eta$.\smallskip
\newline 
\ref{M-Assumption E: Phi hat} trivially holds as $\widehat{\Phi}_n(x)=\widehat{\Phi}_n^*(x)=x$.
\newline
\ref{M-Assumption E: covariance kernel}: 
Here $\psi_{\x}(v;\mathbf{Z})=\1(X\leq \x+v)-\1(X\leq\x) -\Phi_0(\x)v$. Then,
\begin{equation*}
    \frac{\E[|\psi_{\x}(v;\mathbf{Z}) - \psi_{\x}(v';\mathbf{Z})|]}{|v-v'|}\leq \frac{\E[\1(\x + v \land v'<X\leq \x + v \lor v')]}{|v-v'|} +f_0(\x) \leq C.
\end{equation*}
Also, $\psi_{\x_n}(s\eta_n;\mathbf{Z})=\1(\x_n \land (\x_n+s\eta_n) < X\leq \x_n \lor (\x_n+s\eta_n) ) - f_0(\x) s\eta_n$ and
\begin{align*}
    \psi_{\x_n}(s\eta_n;\mathbf{Z})\psi_{\x_n}(t\eta_n;\mathbf{Z})&= \1( \x_n < X\leq \x_n + \eta_n (s \land t) )\1(s>0,t>0)\\
    &\quad + \1( \x_n + (s \lor t) < X\leq \x_n )\1(s<0,t<0) \\
    &\quad - \1(\x_n \land (\x_n+s\eta_n) < X\leq \x_n \lor (\x_n+s\eta_n) ) f_0(\x_n) t\eta_n\\
    &\quad - \1(\x_n \land (\x_n+t\eta_n) < X\leq \x_n \lor (\x_n+t\eta_n) ) f_0(\x_n) s\eta_n\\
    &\quad +f_0(\x_n)^2 st \eta_n^2.
\end{align*}
Then, for any $s,t\in\R$ and $\x_n\to\x$, using continuity of $f_0$ at $\x$,
\begin{align*}
    &\eta_n^{-1}\E[\psi_{\x_n}(s\eta_n;\mathbf{Z})\psi_{\x_n}(t\eta_n;\mathbf{Z})]\\ &= \eta_n^{-1} \int_{\x_n}^{\x_n + \eta_n (s \land t)} f_0(u)du \1(s>0,t>0)  \\
    &\quad +\eta_n^{-1} \int_{\x_n + \eta_n (s \lor t)}^{\x_n} f_0(u)du \1(s<0,t<0) + o(1)\\
    &= f_0(\x) [(s \land t)\1(s>0,t>0) - (s \lor t)\1(s<0,t<0)] + o(1) \\
    &= f_0(\x) (|s| \land |t|)\1(\sign(s)=\sign(t)) + o(1).
\end{align*}

\noindent
\ref{M-Assumption B: uhat} holds since $\widehat{u}_n=\widehat{u}_n^*$ converges in probability to $u_0$, the supremum of the support of $X$ by i.i.d.\ assumption.\smallskip
\newline
\ref{M-Assumption B: Phi uhat} and \ref{M-Assumption B: jump of Phi hat} hold trivially since $\widehat{\Phi}_n$ and $\widehat{\Phi}_n^*$ are the identity map.\smallskip
\newline
Assumption \ref{M-Assumption D} follows from \ref{M-Assumption E: Gamma hat} and empirical process theory arguments.

\subsection{Proof of Corollary \ref{M-[Corollary] monotone density independent right-censoring}}
\label{Appendix subsection: proof of Theorem monotone density independent right-censoring}
Assumption \ref{M-Assumption A} and \ref{M-Assumption E: iid}-\ref{M-Assumption E: bootstrap weights} follow from the hypothesis. 
\newline
\ref{M-Assumption E: Gamma hat}:
We have $\widehat{\Gamma}_n=1-\widehat{S}_n$ with $\widehat{S}_n$ the Kaplan-Meier estimator.
By Theorem 1 of \cite{Lo-Singh_1986_PTRF},
\begin{equation*}
    \sup_{x\in I}\bigg|\widehat{\Gamma}_n(x) -\frac{1}{n}\sum_{i=1}^n \gamma_0(x;\mathbf{Z}) \bigg| = O_{\P}\Big( \Big|\frac{\log n}{n}\Big|^{3/4}\Big).
\end{equation*}
Since $\sqrt{n a_n} \leq n^{2/3}$ for $\q\geq 1$, $\sup_{x \in I}|\widehat{\Gamma}_n(x)-\Gamma_0(x)|=o_{\P}(1)$ and $\sqrt{na}\sup_{|v|\leq \delta}|\widehat{\Gamma}_n(\x+v)-\widehat{\Gamma}_n(\x) - \bar{\Gamma}_n(\x+v)+\bar{\Gamma}_n(\x)|=o_{\P}(1)$ hold. 

We have
\begin{align*}
    &\widehat{\gamma}_n(x;\mathbf{Z}) - \gamma_0(x;\mathbf{Z}) \\
    &= \widehat{F}_n(x) -F_0(x) + \big[\widehat{S}_n(x)-S_0(x)\big]\bigg[ \frac{\1(\check{X}\leq x)\Delta}{\widehat{S}_n(\check{X})\widehat{G}_n(\check{X})} -\int_0^{\check{X}\land x} \frac{\widehat{\Lambda}_n(du)}{\widehat{S}_n(u)\widehat{G}_n(u)} \bigg] \\
    &\quad + S_0(x)\1(\check{X}\leq x)\Delta\frac{S_0(\check{X})G_0(\check{X})-\widehat{S}_n(\check{X})\widehat{G}_n(\check{X})}{S_0(\check{X})G_0(\check{X})\widehat{S}_n(\check{X})\widehat{G}_n(\check{X})}\\
    &\quad - S_0(x) \int_0^{\check{X}\land x} \frac{S_0(u)G_0(u)-\widehat{S}_n(u)\widehat{G}_n(u)}{S_0(u)G_0(u)\widehat{S}_n(u)\widehat{G}_n(u)}\widehat{\Lambda}_n(du)  - S_0(x)\int_0^{\check{X}\land x} \frac{[\widehat{\Lambda}_n-\Lambda_0](du)}{S_0(u)G_0(u)}.
\end{align*}
By $S_0(u_0)G_0(u_0)>0$, we have $\sqrt{n}\sup_{x\in I}|\widehat{S}_n(x)-S_0(x)|=O_{\P}(1)$, $\sqrt{n}\sup_{x\in I}|\widehat{G}_n(x)-G_0(x)|=O_{\P}(1)$, and $\sqrt{n}\sup_{x\in I}|\widehat{\Lambda}_n(x)-\Lambda_0(x)|=O_{\P}(1)$, which in turn implies
\begin{equation*}
    \frac{1}{n}\sum_{i=1}^n \sup_{x\in I}\big|\widehat{\gamma}_n(x;\mathbf{Z}_i)-\gamma_0(x;\mathbf{Z}_i)\big|^2=o_{\P}(1),\ a_n \frac{1}{n}\sum_{i=1}^n \sup_{x\in I_{\x}^{\delta}}\big|\widehat{\gamma}_n(x;\mathbf{Z}_i)-\gamma_0(x;\mathbf{Z}_i)\big|^2=o_{\P}(1).
\end{equation*}

Let $\delta = \min\{\x, (u_0-\x)\}/4$, $R_{1n}(v)=| \widehat{F}_n(\x+v)-\widehat{F}_n(\x) - F_0(\x+v)+F_0(\x)|$, $R_{2n}=| \widehat{F}_n(\x)-F_0(\x)|$, $R_{3n}= \sup_{x\in I}| [\widehat{S}_n(x)\widehat{G}_n(x)]^{-1} - [S_0(x)G_0(x)]^{-1} |$, and $R_{4n}(x_1,x_2) = |\int_{x_1}^{x_2} \frac{\widehat{\Lambda}_n(du)}{\widehat{S}_n(u)\widehat{G}_n(u)} -\int_{x_1}^{x_2} \frac{\Lambda_0(du)}{S_0(u)G_0(u)}| $. For $|v|\leq \delta$,
\begin{align*}
    &\big|\check{\Gamma}_n(\x+v)-\check{\Gamma}_n(\x) - \bar{\Gamma}_n(\x+v) + \bar{\Gamma}_n(\x)\big| \\
    &\leq R_{1n}(v) \left[1 + \frac{1 + \widehat{\Lambda}_n(x) }{\widehat{S}_n(x)\widehat{G}_n(x)}\right] + R_{2n}\big[\widehat{S}_n(\x-\delta)\widehat{G}_n(\x-\delta) \big]^{-1} \frac{1}{n}\sum_{i=1}^n \1(|X_i-\x|\leq v )\\
    &+ R_{2n} \int_{\x-v}^{\x+v} \frac{\widehat{\Lambda}_n(du)}{\widehat{S}_n(u)\widehat{G}_n(u)} + \sup_{x\in I_{\x}^{\delta}} f_0(x)| v | \bigg[R_{3n} +  \frac{1}{n}\sum_{i=1}^n R_{4n}(0,\check{X}_i\land (\x+v) )   \bigg] \\
    & +  R_{3n} \frac{1}{n}\sum_{i=1}^n \1(|X_i-\x|\leq v ) + \frac{1}{n}\sum_{i=1}^n R_{4n}(\check{X}_i\land\x,\check{X}_i\land (\x+v) ).
\end{align*}
As noted above, $\sup_{|v|\leq \delta}R_{1n}(v)=o_{\P}( (na_n)^{-1/2})$, $R_{2n}=O_{\P}(n^{-1/2})$, and $R_{3n}=O_{\P}(n^{-1/2})$ by $S_0(u_0)G_0(u_0)>0$.
Also, uniformly over $V\in (0, 2\delta]$, $\sup_{| v|\leq V}| \frac{1}{n}\sum_{i=1}^n  (\1(X_i\leq \x +v) - \1(X_i\leq \x) )| \leq C V + O_{\P}(n^{-1/2})$. If 
\begin{equation}\label{Appendix: Example monotone density right-censoring R3 1}
    \sup_{|v|\leq \delta}\frac{1}{n}\sum_{i=1}^n R_{4n}(0,\check{X}_i\land (\x+v) )=O_{\P}(n^{-1/2})
\end{equation}
and 
\begin{equation}\label{Appendix: Example monotone density right-censoring R3 2}
    \frac{1}{n}\sum_{i=1}^n R_{4n}(\check{X}_i\land\x,\check{X}_i\land (\x+v) ) \leq |v| O_{\P}(n^{-1/2}) + o_{\P}\Big(( na_n)^{-1/2} \Big)
\end{equation}
uniformly over $|v|\leq 2\delta$, then there exist random variables $A_n=o_{\P}(1)$ and $B_n=O_{\P}(\sqrt{a_n})$ independent of $v$ such that for $V\in (0, 2\delta]$,
\begin{equation*}
    \sqrt{na_n} \sup_{|v|\leq V} \big|\check{\Gamma}_n(\x+v)-\check{\Gamma}_n(\x) - \bar{\Gamma}_n(\x+v) + \bar{\Gamma}_n(\x)\big| \leq A_n + V B_n;
\end{equation*}
that is,\ $\beta_{\gamma}=1$ in the notation of \ref{M-Assumption E: Gamma hat}. To show \eqref{Appendix: Example monotone density right-censoring R3 1} and \eqref{Appendix: Example monotone density right-censoring R3 2}, for $x_1,x_2\in I$,
\begin{equation*}
    R_{4n}(x_1,x_2) \leq R_{3n}\big\vert \widehat{\Lambda}_n(x_2)- \widehat{\Lambda}_n(x_1)\big\vert  + \left\vert \int_{x_1}^{x_2}\frac{[\widehat{\Lambda}_n-\Lambda_0](du)}{S_0(u)G_0(u)}\right\vert.
\end{equation*}
Let $J_0(u)= 1/[S_0(u)G_0(u)]$ and integration by parts implies
\begin{align*}
    \int_{x_1}^{x_2} \frac{[\widehat{\Lambda}_n-\Lambda_0](du)}{S_0(u)G_0(u)} &=  \frac{\widehat{\Lambda}_n(x_2)-\Lambda_0(x_2)}{S_0(x_2)G_0(x_2)} - \frac{\widehat{\Lambda}_n(x_1)-\Lambda_0(x_1)}{S_0(x_1)G_0(x_1)} - \int_{x_1}^{x_2} [\widehat{\Lambda}_n(u)-\Lambda_0(u)]J_0(du) \\
    &= [\widehat{\Lambda}_n(x_2)-\Lambda_0(x_2)]\left[\frac{1}{S_0(x_2)G_0(x_2)}-\frac{1}{S_0(x_1)G_0(x_1)}\right] \\
    & + \frac{\widehat{\Lambda}_n(x_2)-\Lambda_0(x_2) -\widehat{\Lambda}_n(x_1)+\Lambda_0(x_1)}{S_0(x_1)G_0(x_1)}  - \int_{x_1}^{x_2} [\widehat{\Lambda}_n(u)-\Lambda_0(u)]J_0(du).
\end{align*}
The first term after the second equality is bounded by $O_{\P}(n^{-1/2}) \vert x_2-x_1\vert$. For the second term, Theorem 1 of \cite{Burketal} implies that on a suitable probability space there exists a sequence of standard Brownian motion $W_n$ such that $\sqrt{a_n}\sup_{\vert x_1-x_2\vert\leq v} \vert \sqrt{n}[\widehat{\Lambda}_n(x_1)-\Lambda_0(x_2) -\widehat{\Lambda}_n(x_1)+\Lambda_0(x_2)] - W_n(d(x_1))+W_n(d(x_2))\vert = o_{\P}(1)$, where $d(x)=\int_0^x \frac{F_0(du)}{S_0(u)^2G_0(u)}$. By Theorem 3.2 of \cite{Pollard_1989_SS}, there is some fixed constant $C>0$ such that $\E[\sup_{\vert x_1-x_2\vert\leq v}\vert W_n(d(x_1))-W_n(d(x_2))\vert]\leq C v$. Finally, $\vert \int_{x_1}^{x_2} [\widehat{\Lambda}_n(u)-\Lambda_0(u)]J_0(du)\vert\leq \sup_{x'\in [x_1,x_2]}\vert\widehat{\Lambda}_n(x')-\Lambda_0(x')\vert [J_0(x_2)-J_0(x_1)]\leq O_{\P}(n^{-1/2}) \vert x_2-x_1\vert$. Thus, \eqref{Appendix: Example monotone density right-censoring R3 1} and \eqref{Appendix: Example monotone density right-censoring R3 2} hold.

For the function class $\mathfrak{F}_{\gamma}$, we can take $\bar{F}_{\gamma}(\mathbf{Z}) = 1 + [S_0(u_0)G_0(u_0)]^{-1}[1 + \Lambda_0(u_0)]$ as a constant envelope. 
For the function class $\{S_0(x):x\in I\}$, given $m\in\mathbbm{N}$, there exists $\{x_1,\dots, x_{m+1}\}\subset I$ such that $\sup_{x\in I}\min_{l=1,\dots,m+1}|S_0(x_l)-S_0(x)|\leq 1/m$, which implies the uniform covering number is bounded by a linear function. The covering numbers of $\{\1(\cdot\leq s):s\in I\}$ and $\{\int_0^{\cdot\land s}[S_0(u)G_0(u)]^{-1}\Lambda_0(du):s\in I\}$ are also bounded by a linear function. By Lemma 5.1 of \cite{vanderVaart-vanderLaan_2006_IJB}, there exists $\rho\in (0,2)$ such that $\limsup_{\eta\downarrow 0}\log N_U(\eta,\mathfrak{F}_{\gamma})\eta^{\rho} < \infty$ holds.

Now consider the uniform covering number of $\hat{\mathfrak{F}}_{\gamma}$. Given a realization of $(\widehat{S}_n,\widehat{G}_n)$, the mapping $x\mapsto\int_0^{x\land s}[\widehat{S}_n(u)\widehat{G}_n(u)]^{-1}\widehat{\Lambda}_n(du)$ is a composition of $x\mapsto x\land s$ and $x\mapsto \int_0^x [\widehat{S}_n(u)\widehat{G}_n(u)]^{-1}\widehat{\Lambda}_n(du)$. The latter mapping is monotone, and the first mapping is a VC-subgraph class, and Lemma 2.6.18 of \cite{vanderVaart-Wellner_1996_Book} implies $\{\int_0^{\cdot\land s}[\widehat{S}_n(u)\widehat{G}_n(u)]^{-1}\widehat{\Lambda}_n(du):s\in I\}$ is a VC-subgraph class. Note that since $S_0,G_0$ are bounded away from zero, $\widehat{S}_n,\widehat{G}_n$ are bounded away from zero with probability approaching one. Thus, for some $\rho\in (0,2)$, $\limsup_{\eta\downarrow 0}\log N_U(\eta,\hat{\mathfrak{F}}_{\gamma})\eta^{\rho} = O_{\P}(1)$ holds.

For $s\leq t\in I$,
\begin{equation*}
    |\gamma_0(s;\mathbf{Z})-\gamma_0(t;\mathbf{Z})|\leq C|F_0(s)-F_0(t)| + C|\1(\check{X}\leq s)-\1(\check{X}\leq t)|\Delta + \int_{\check{X}\land s}^{\check{X}\land t}\frac{\Lambda_0(du)}{S_0(du)G_0(du)}
\end{equation*}
and we can take $D_{\gamma}^{\eta}(\mathbf{Z})$ to be a constant multiple of $\sup_{|s|\leq \eta}|F_0(\x +s)-F_0(\x)| + \Delta \1(|\check{X}-\x|\leq \eta) + \int_{ x-\eta}^{x+\eta}\Lambda_0(du)/S_0(u)G_0(u)$. For $\eta >0$ small enough, there is some fixed $C>0$ with
\begin{equation*}
    \E[D_{\gamma}^{\eta}(\mathbf{Z})^2 + D_{\gamma}^{\eta}(\mathbf{Z})^4] \leq C f_0(\x+\eta) \eta.
\end{equation*}
\ref{M-Assumption E: Phi hat} trivially holds as $\widehat{\Phi}_n(x)=\widehat{\Phi}_n^*(x)=x$.
\newline
\ref{M-Assumption E: covariance kernel}:
We have
\begin{equation*}
    \psi_{\x}(v;\mathbf{Z}) = S_0(\x)\frac{(\1(\check{X}\leq \x +v)-\1(\check{X}\leq \x))\Delta}{S_0(\check{X})G_0(\check{X})}  + O(|v|)
\end{equation*}
where $O(|v|)$ is uniformly over small enough $|v|$. Since
\begin{equation*}
    \E[|\1(\check{X}\leq \x+v) - \1(\check{X}\leq \x+v')|\Delta]=\int_{\x + v\land v'}^{\x+v\lor v'} G_0(u) f_0(u)du \leq C|v-v'|,
\end{equation*}
the first display in \ref{M-Assumption E: covariance kernel} is satisfied. For the covariance kernel,
\begin{align*}
    &\E[\psi_{\x_n}(s\eta_n;\mathbf{Z})\psi_{\x_n}(t\eta_n;\mathbf{Z})] \\
    &= S_0(\x_n)^{2}\bigg( \int_{\x_n}^{\x_n + \eta_n(s \land t)}\frac{f_0(u)}{S_0(u)^2G_0(u)}du  \1(s>0,t>0)\\
    &\qquad\qquad+\int_{\x_n + \eta_n(s \lor t)}^{\x_n}\frac{f_0(u)}{S_0(u)^2G_0(u)}du\1(s<0,t<0) \bigg)+O(\eta_n^2)\\
    &= \frac{S_0(\x_n)^{2}f_0(\x)}{S_0(\x)^2G_0(\x)} (s \land t)\eta_n \1(s>0,t>0) \\
    &\quad - \frac{S_0(\x_n)^{2}f_0(\x)}{S_0(\x)^2G_0(\x)} (s \lor t)\eta_n \1(s<0,t<0) + o(\eta_n)
\end{align*}
where the last equality uses continuity of $(S_0,G_0,f_0)$ at $\x$ i.e.,\ $\int_{\x_n}^{\x_n +\eta_n} [\frac{f_0(u)}{S_0(u)^2G_0(u)}-\frac{f_0(\x)}{S_0(\x)^2G_0(\x)}]du =o(1)\eta_n$.

\noindent
\ref{M-Assumption B: uhat}, \ref{M-Assumption B: Phi uhat}, and \ref{M-Assumption B: jump of Phi hat} hold since in this example, $\widehat{u}_n=\widehat{u}_n^*=u_0$ and $\widehat{\Phi}_n,\widehat{\Phi}_n^*$ are the identity map.

\paragraph*{Assumption \ref{M-Assumption D}}
As noted when verifying \ref{M-Assumption E: Gamma hat}, $\widehat{G}_n(1;\eta_n)=o_{\P}(1)$ for any $\eta_n=o(1)$ with $a_n^{-1}\eta_n^{-1}=O(1)$. $\widehat{\Phi}_n=\Phi_0$ is the identity map and the desired result holds.

\subsection{Proof of Corollary \ref{M-[Corollary] classical monotone regression}}\label{Appendix subsection: proof theorem classical monotone regression}
Assumption \ref{M-Assumption A} and \ref{M-Assumption E: iid}-\ref{M-Assumption E: bootstrap weights} follow from the hypothesis.
\newline
\ref{M-Assumption E: Gamma hat}:
In this example, $\gamma_0(x;\mathbf{Z})=Y\1(X\leq x)$ is known, so it suffices to verify \eqref{Appendix Eq: simpler D gamma hat}. The uniform covering number bound is straightforward as $\{\1(\cdot\leq x):x\in\R\}$ is a VC-subgraph class. An envelope function is $|Y|$, whose second moment is finite. For $x\in I_{\x}^{\eta}$, $|\gamma_0(x;\mathbf{Z})-\gamma_0(\x;\mathbf{Z})|\leq |Y| \1(\x-\eta \leq X\leq \x +\eta)$, which we can take as $\bar{D}_{\gamma}^{\eta}(\mathbf{Z})$. Then, for $j=2,4$,
\begin{equation*}
    \E[\bar{D}_{\gamma}^{\eta}(\mathbf{Z})^j] \leq 2^{j-1} \int_{\x-\eta}^{\x+\eta} \big(|\mu_0(x)|^j +\E[\varepsilon^j|X=x]\big) f_0(x)dx \leq C \eta
\end{equation*}
and the desired bound holds.
\newline
\ref{M-Assumption E: Phi hat}:
$\phi_0(x;\mathbf{Z})=\1(X\leq x)$ is known, so it suffices to verify the analogue of \eqref{Appendix Eq: simpler D gamma hat}. The argument is the same as for checking \ref{M-Assumption E: Gamma hat} in monotone density estimation with no censoring.

\noindent
\ref{M-Assumption E: covariance kernel}:
We have
\begin{equation*}
    \psi_{\x}(v;\mathbf{Z}) = \varepsilon(\1(X\leq \x +v)-\1(X\leq \x))  + (\mu_0(X)-\mu_0(\x)) (\1(X\leq \x +v)-\1(X\leq \x)).
\end{equation*}
Then,
\begin{equation*}
    \E[|\psi_{\x}(v;\mathbf{Z})-\psi_{\x}(v';\mathbf{Z})|]\leq \int_{\x+v \land v'}^{\x + v\lor v'} [\sigma_0(x) +|\mu_0(x)-\mu_0(\x)|] f_0(x)dx \leq C |v-v'|
\end{equation*}
and the first display holds. 
For the covariance kernel, note $|(\mu_0(X)-\mu_0(\x_n))(\1(X\leq \x_n+v)-\1(X\leq \x_n)|\leq |v| \sup_{|x-\x|\leq 2\eta} |\partial\mu_0(x)|$ for $|x_n-\x|\lor |v|\leq \eta$ for $\eta>0$ small enough. Then,
\begin{align*}
    &\E[\psi_{\x_n}(s\eta_n;\mathbf{Z})\psi_{\x_n}(t\eta_n;\mathbf{Z})]\\
    &= \E[\varepsilon^2(\1(X\leq \x_n +s\eta_n)-\1(X\leq \x_n))(\1(X\leq \x_n +t\eta_n)-\1(X\leq \x_n)) ] + O(\eta_n^2)\\
    &= \int_{\x_n}^{\x_n+\eta_n (s\land t)} \sigma_0^2(x)f_0(x)dx\1(s>0,t>0) \\
    &\quad + \int^{\x_n}_{\x_n+\eta_n (s\lor t)} \sigma_0^2(x)f_0(x)dx\1(s<0,t<0) + O(\eta_n^2)
\end{align*}
and
\begin{align*}
   &\eta_n^{-1}\E[\psi_{\x_n}(s\eta_n;\mathbf{Z})\psi_{\x_n}(t\eta_n;\mathbf{Z})] \\
   &\to \sigma_0^2(\x)f_0(\x) \big[ (s\land t) \1(s>0,t>0) - (s\lor t) \1(s<0,t<0) \big]\\
   &= \sigma_0^2(\x)f_0(\x) (|s|\land |t|) \1(\sign(s)=\sign(t)),
\end{align*}
as desired.

\noindent
\ref{M-Assumption B: uhat} trivially holds since $\widehat{u}_n=\widehat{u}_n^*=1$ in this example.

\noindent \ref{M-Assumption B: Phi uhat}:
$\widehat{\Phi}_n,\widehat{\Phi}_n^*$ are (empirical) cdfs, so they are non-negative, non-decreasing, and right-continuous. $\{0,\widehat{u}_n\}\subset \widehat{\Phi}_n(I)$ and $\{0,\widehat{u}_n^*\}\subset \widehat{\Phi}_n^*(I)$ hold as $\widehat{u}_n=\widehat{u}_n^*=1$, $\widehat{\Phi}_n(\min_iX_i-)=0=\widehat{\Phi}_n^*(\min_iX_i-)$, and $\widehat{\Phi}_n(\max_iX_i)=1=\widehat{\Phi}_n^*(\max_iX_i)$. The sets $\widehat{\Phi}_n(I),\widehat{\Phi}_n^*(I)$ are finite and thus closed.

\noindent \ref{M-Assumption B: jump of Phi hat}:
With probability one, all $X_i$'s are distinct. If $x$ is one of $X_i$'s, $\widehat{\Phi}_n^*( x ) - \widehat{\Phi}_n^*( x- ) = n^{-1}W_{j,n}$ for some $j$, and $\widehat{\Phi}_n^*( x ) - \widehat{\Phi}_n^*( x- ) =0$ otherwise.

Given $\E|W_{1,n}|^{\mathfrak{r}}<\infty$, $n^{-1}\max_{1\leq i\leq n}|W_{i,n}|=o_{\P}(n^{-5/6})$, which implies $$\sqrt{na_n}\sup_{x\in I}|\widehat{\Phi}_n^*( x ) - \widehat{\Phi}_n^*( x- ) |=o_{\P}(1).$$ The argument for $\widehat{\Phi}_n$ is similar.
\smallskip
\newline
Assumption \ref{M-Assumption D} follows from \ref{M-Assumption E: Gamma hat}-\ref{M-Assumption E: Phi hat} and empirical process theory arguments.

\subsection{Proof of Corollary \ref{M-[Corollary] monotone regression with covariates}}\label{Appendix subsection: proof theorem monotone regression with covariates}
Assumption \ref{M-Assumption A} and \ref{M-Assumption E: iid}-\ref{M-Assumption E: bootstrap weights} follow from the hypothesis.

\noindent \ref{M-Assumption E: Gamma hat}:
In this example, $\check{\Gamma}_n=\widehat{\Gamma}_n$.
\begin{align*}
    |\widehat{\gamma}_n(x;\mathbf{Z}) - \gamma_0(x;\mathbf{Z})|&\leq \1(X\leq  x) \bigg[|\varepsilon|\Big| \widehat{g}_n(X,\mathbf{A})^{-1}-g_0(X,\mathbf{A})^{-1}\Big| + \frac{|\widehat{\mu}_n(X,\mathbf{A})-\mu_0(X,\mathbf{A})|}{\widehat{g}_n(X,\mathbf{A})}\\
    & + \frac{1}{n}\sum_{j=1}^n |\widehat{\mu}_n(X,\mathbf{A}_j) - \mu_0(X,\mathbf{A}_j)|+ \bigg|\frac{1}{n}\sum_{j=1}^n \mu_0(X,\mathbf{A}_j) - \theta_0(X)\bigg| \bigg].
\end{align*}
The last sum is bounded by $\sup_{x\in I}|\frac{1}{n}\sum_{j=1}^n\mu_0(x,\mathbf{A}_j)-\theta_0(x)|$, and this object is $O_{\P}(n^{-1/2})$: to see this claim, first note that Assumption MRC (iv) and Theorem 2.7.11 of \cite{vanderVaart-Wellner_1996_Book} imply $\limsup_{\epsilon\downarrow0}\log N_U(\epsilon,\{\mu(x,\cdot):x\in I\}) \epsilon^V <\infty$ for some $V\in (0,2)$ and Theorem 4.2 of \cite{Pollard_1989_SS} implies $\sup_{x\in I}|\frac{1}{n}\sum_{j=1}^n\mu_0(x,\mathbf{A}_j)-\theta_0(x)|=O_{\P}(n^{-1/2})$.
Together with Assumption MRC (iii), $a_n \frac{1}{n}\sum_{i=1}^n \sup_{x\in I}|\widehat{\gamma}_n(x;\mathbf{Z})-\gamma_0(x;\mathbf{Z})|^2=o_{\P}(1)$ holds.

By Assumption MRC (iii), uniformly over $V\in (0,2\delta]$
\begin{equation*}
    \sqrt{na_n}\sup_{|v|\leq V}\big\vert\widehat{\Gamma}_n(\x+v)-\widehat{\Gamma}_n(\x) - \bar{\Gamma}_n(\x+v)+\bar{\Gamma}_n(\x)\big\vert \leq o_{\P}(1) + O_{\P}(\sqrt{a_n}) V
\end{equation*}
and the desired inequality holds.

The uniform covering numbers of $\mathfrak{F}_{\gamma},\hat{\mathfrak{F}}_{\gamma,n}$ are the same order as for $\{\1(\cdot\leq x):x\in I\}$. For $x \in I_{\x}^{\eta}$, $|\gamma_0(x;\mathbf{Z})-\gamma_0(\x;\mathbf{Z})|\leq \1(\x-\eta\leq X\leq \x+\eta)(|\varepsilon| c^{-1}+\theta_0(\x+\eta) )$. Then, $\limsup_{\eta\downarrow0}\E[\bar{D}_{\gamma}^{\eta}(\mathbf{Z})^j]\eta^{-1} <\infty$ holds for $j=2,4$.

\noindent \ref{M-Assumption E: Phi hat}:
$\phi_0(x;\mathbf{Z})=\1(X\leq x)$ is known and the same as in the classical case, so the same argument applies.

\noindent \ref{M-Assumption E: covariance kernel}:
We have
\begin{equation*}
    \psi_{\x}(v;\mathbf{Z}) = (\1(X\leq \x+v)-\1(X\leq \x))\bigg[\frac{\varepsilon}{g_0(X,\mathbf{A})} + \theta_0(X)-\theta_0(\x) \bigg].
\end{equation*}
Then, for $v,v'\in [-\eta,\eta]$ with sufficiently small $\eta>0$,
\begin{align*}
    |\psi_{\x}(v;\mathbf{Z}) - \psi_{\x}(v';\mathbf{Z})|&\leq |\1(X\leq \x+ v)-\1(X\leq \x+ v')|\big( c^{-1} |\varepsilon| + |X-\x| \sup_{x\in I_{\x}^{\eta}}|\partial \theta_0(x)|\big) 
\end{align*}
and $\sup_{v\neq v'\in [-\eta_n,\eta_n]}\E[|\psi_{\x}(v;\mathbf{Z}) - \psi_{\x}(v';\mathbf{Z})|]/|v-v'|=O(1)$ holds.

For $s\eta_n$ small enough, $\psi_{\x}(s\eta_n;\mathbf{Z})= (\1(X\leq \x+ s\eta_n)-\1(X\leq \x))\varepsilon g_0(X,\mathbf{A})^{-1} + O(\eta_n)$ and
\begin{align*}
    &\E[\psi_{\x}(s\eta_n;\mathbf{Z})\psi_{\x}(t\eta_n;\mathbf{Z})] \\
    &= \E\bigg[\frac{\sigma_0^2(X,\mathbf{A})}{g_0(X,\mathbf{A})^2}(\1(X\leq \x+ s\eta_n)-\1(X\leq \x))(\1(X\leq \x+ t\eta_n)-\1(X\leq \x)) \bigg] + O(\eta_n^2)
\end{align*}
and
\begin{align*}
    & \E\bigg[\frac{\sigma_0^2(X,\mathbf{A})}{g_0(X,\mathbf{A})^2}(\1(X\leq \x+ s\eta_n)-\1(X\leq \x))(\1(X\leq \x+ t\eta_n)-\1(X\leq \x)) \bigg]\\
    &= \E \bigg[\int_{\x}^{\x+\eta_n (s\land t)} \frac{\sigma_0^2(x,\mathbf{A})}{g_0(x,\mathbf{A})^2} f_{X|\mathbf{A}}(x|\mathbf{A})dx \bigg]\1(s>0,t>0)\\
    &\quad + \E \bigg[\int_{\x+\eta_n (s\lor t)}^{\x} \frac{\sigma_0^2(x,\mathbf{A})}{g_0(x,\mathbf{A})^2} f_{X|\mathbf{A}}(x|\mathbf{A})dx \bigg]\1(s<0,t<0)\\
    &= \E\bigg[\frac{\sigma_0^2(\x,\mathbf{A})}{g_0(\x,\mathbf{A})^2}f_{X|\mathbf{A}}(\x|\mathbf{A})\bigg]\Big(\eta_n (s\land t) \1(s>0,t>0) - \eta_n (s\lor t) \1(s<0,t<0) \Big) + o(\eta_n).
\end{align*}
Since $\frac{f_{X|\mathbf{A}}(x|\mathbf{A})}{g_0(x,\mathbf{A})^2} = \frac{f_0(x)}{g_0(x,\mathbf{A})}$, we have 
\begin{equation*}
    \eta_n^{-1}\E[\psi_{\x}(s\eta_n;\mathbf{Z})\psi_{\x}(t\eta_n;\mathbf{Z})] \to f_0(\x) \E\bigg[\frac{\sigma_0^2(\x,\mathbf{A})}{g_0(\x,\mathbf{A})}\bigg]\Big( (s\land t) \1(s,t>0) -  (s\lor t) \1(s,t<0) \Big),
\end{equation*}
as desired.

\noindent \ref{M-Assumption B: uhat} \ref{M-Assumption B: Phi uhat} \ref{M-Assumption B: jump of Phi hat}:
Verifying these conditions is the same as in the classical monotone regression case.

\section{Primitive sufficient conditions for Assumption MRC (iii)}

Here we provide primitive sufficient conditions for Assumption MRC (iii) by focusing on specific estimators $\widehat{\mu}_n$ and $\widehat{g}_n$. As discussed by \cite{Westling-Gilbert-Carone_2020_JRRSB}, cross-fitting avoids restrictions on uniform entropy, allowing for a large class of flexible preliminary estimators. Here we use sample splitting to simplify exposition, but the proposed procedure can be straightforwardly modified for cross-fitting.

Suppose there is a separate random sample $\tilde{\mathbf{Z}}_1,\dots,\tilde{\mathbf{Z}}_n$ drawn from the distribution of $\mathbf{Z}$, which is independent of $\mathbf{Z}_1,\dots,\mathbf{Z}_n$. Preliminary estimators $\widehat{\mu}_n$ and $\widehat{g}_n$ are constructed from $\tilde{\mathbf{Z}}_1,\dots,\tilde{\mathbf{Z}}_n$.
For concreteness, we consider a partitioning-based least squares estimator $\widehat{\mu}_n$ \citep{Cattaneo-Farrell-Feng_2020_AoS} and local polynomial kernel-based estimators $\widehat{f}_{X|\mathbf{A},n}(x|\mathbf{a})$ and $\widehat{f}_n(x)$ of $f_{X|\mathbf{A}}(x|\mathbf{a})$ and $f_0(x)$ \citep*{Cattaneo-Jansson-Ma_2020_JASA,Cattaneo-Chandak-Jansson-Ma_2024}, from which we construct $\widehat{g}_n(x,\mathbf{a})=\widehat{f}_{X|\mathbf{A},n}(x|\mathbf{a})/\widehat{f}_n(x)$.  

Let $d=\mathrm{dim}(\mathbf{A})$. For simplicity, suppose the support of $(X,\mathbf{A}')'$ equals $[0,1]^{1+d}$. Let $\mathbf{p}(x,\mathbf{a})$ be a $k_n$-dimensional vector of bounded basis functions of order $m$ on $\mathcal{S}$ which are locally supported  e.g.,\ splines \cite[see][for details and examples of basis functions]{Cattaneo-Farrell-Feng_2020_AoS}. We consider the estimator 
\begin{equation*}
    \widehat{\mu}_n(x,\mathbf{a}) = \mathbf{p}(x,\mathbf{a})'\bigg(\sum_{i=1}^n\mathbf{p}(\tilde{X}_i,\tilde{\mathbf{A}}_i) \mathbf{p}(\tilde{X}_i,\tilde{\mathbf{A}}_i)'\bigg)^{-1} \sum_{i=1}^n \mathbf{p}(\tilde{X}_i,\tilde{\mathbf{A}}_i)\tilde{Y}_i.
\end{equation*}
For the estimator of $f_{X|\mathbf{A}}(x\vert\mathbf{a})$, letting $\widehat{F}_{X|\mathbf{A},n}(\cdot\vert \mathbf{a})$ be an estimator of $\P[X\leq \cdot \vert \mathbf{A}=\mathbf{a}]$ specified below, $\widehat{f}_{X|\mathbf{A},n}(x|\mathbf{a})$ is obtained by local polynomial regression:
\begin{equation*}
    \widehat{f}_{X|\mathbf{A},n}(x|\mathbf{a}) = \mathbf{e}_2' \argmin_{\mathbf{u}\in\R^{\mathfrak{p}_1+1}}\sum_{i=1}^n\Big(\widehat{F}_{X|\mathbf{A},n}(\tilde{X}_i|\mathbf{a}) - \mathbf{q}_1(\tilde{X}_i-x)'\mathbf{u}\Big)^2K_h(\tilde{X}_i-x) 
\end{equation*}
where $\mathfrak{p}_1\geq 1$ is the order of the polynomial basis $\mathbf{q}_1(x)=(1,x/1!, x^2/2!,\dots,x^{\mathfrak{p}_1}/\mathfrak{p}_1!)'$, $\mathbf{e}_l$ is the conformable unit vector whose $l$th element is unity, and $K_h(x)=K(x/h)/h$ for some kernel function $K$ and some positive bandwidth $h$. The estimator $\widehat{F}_{X|\mathbf{A},n}(x|\mathbf{a})$ is constructed via local polynomial regression of order $\mathfrak{p}_2=\mathfrak{p}_1-1$:
\begin{equation*}
    \widehat{F}_{X|\mathbf{A},n}(x|\mathbf{a}) =\mathbf{e}_1' \argmin_{\mathbf{v}\in\R^{k_{\mathfrak{p}_2}}}\sum_{i=1}^n \big(\1(\tilde{X}_i\leq x) - \mathbf{q}_2(\tilde{\mathbf{A}}_i-\mathbf{a})'\mathbf{v}\big)^2L_h(\tilde{\mathbf{A}}_i-\mathbf{a})
\end{equation*}
where, using standard multi-index notation, $\mathbf{q}_2(\mathbf{a})$ denotes the $k_{\mathfrak{p}_2}$-dimensional vector collecting the polynomials $\mathbf{a}^{\mathbf{m}}/\mathbf{m}!$ for $0\leq \vert \mathbf{m}\vert\leq \mathfrak{p}_2$ with $\mathbf{a}^{\mathbf{m}} = a_1^{m_1}a_2^{m_2}\dots a_d^{m_d}$, $\vert\mathbf{m}\vert=\sum_{j=1}^dm_j$, and $k_{\mathfrak{p}_2}=\frac{(d+\mathfrak{p}_2)!}{d!\mathfrak{p}_2!}+1$, and $L_h(\mathbf{a})=L(\mathbf{a}/h)/h^d$ for $L(\mathbf{a})=\prod_{j=1}^dK(a_j)$ i.e.,\ product kernel. 
The estimator $\widehat{f}_n(x)$ is constructed in a similar manner. First, the empirical cdf $\widehat{F}_n$ of $\{\tilde{X}_i\}$ is constructed and then $\widehat{f}_n(x)$ is formed via local polynomial regression:
\begin{equation*}
    \widehat{f}_n(x) =\mathbf{e}_2'\argmin_{\mathbf{w}\in\R^{\mathfrak{p}_1}+1} \sum_{i=1}^n\Big(\widehat{F}_n(\tilde{X}_i)-\mathbf{q}_1(\tilde{X}_i-x)'\mathbf{w}\Big)^2K_b(\tilde{X}_i-x)
\end{equation*}
where $b>0$ is some bandwidth.

Now we state sufficient conditions for Assumption MRC (ii) based on the partitioning-based series estimator $\widehat{\mu}_n$ and the kernel-based estimator $\widehat{g}_n$.
\begin{description}
    \item[Primitive Conditions MRC]
    \item[(i)] $\tilde{\mathbf{Z}}_1,\dots,\tilde{\mathbf{Z}}_n$ are independent of $\mathbf{Z}_1,\dots,\mathbf{Z}_n$.
    
    \item[(ii)] The support of $(X,\mathbf{A}')'$ is $[0,1]^{1+d}$, and the distribution of $(X,\mathbf{A}')'$ is absolutely continuous. The Lebesgue density of $(X,\mathbf{A}')'$ and the conditional variance of $Y$ given $(X,\mathbf{A}')'$ are bounded away from zero and continuous on $[0,1]^{1+d}$. $\mu_0$ is $(m+1)$-times continuously differentiable on $[0,1]^{1+d}$.

    \item[(iii)] The vector of basis functions $p$ satisfies Assumptions 2, 3, and 4 of \cite{Cattaneo-Farrell-Feng_2020_AoS}.

    \item[(iv)] $f_{X|\mathbf{A}}(x|\mathbf{a})$ and $f_0(x)$ are $\mathfrak{p}_1$-times continuously differentiable in $x$, and $f_{X|\mathbf{A}}(x|\mathbf{a})$ is $\mathfrak{p}_1$-times continuously differentiable in $\mathbf{a}$.

    \item[(v)] $K$ is a symmetric, Lipschitz continuous probability density function supported on $[-1,1]$.
\end{description}
As verified by \cite{Cattaneo-Farrell-Feng_2020_AoS}, (iii) holds for widely used local basis functions such as splines and wavelets.

\begin{lemma}
    Suppose $\tilde{\mathbf{Z}}_1,\dots,\tilde{\mathbf{Z}}_n$ is a random sample drawn from the distribution of $\mathbf{Z}$ and Primitive Conditions MRC hold. In addition, with $\tau_n= n^{-1}\log n$, $k_n=O(\tau_n^{-\frac{d+1}{2m+d+1}})$, $h=O(\tau_n^{\frac{d+1}{2\mathfrak{p}_1+d+1}})$, $b=O(\tau_n^{\frac{1}{2\mathfrak{p}_1+1}})$, $\frac{m}{2m+d+1}+\frac{\mathfrak{p}_1}{2\mathfrak{p}_1+d+1} \geq \frac{1}{2}$, and $\min\{\frac{2m}{2m+d+1},\frac{2\mathfrak{p}_1}{2\mathfrak{p}_1+d+1} \} > \frac{1}{2\q+1}$. Then, $\widehat{\mu}_n$ and $\widehat{g}_n$ described above and $\widehat{\Gamma}_n$ based on the $\widehat{\mu}_n,\widehat{g}_n$ satisfy Assumption MRC (iii) with $\delta = \min\{\x,1-\x\}/4$. In particular,
    \begin{equation}\label{Appendix eq: result - monotone regression with covariates primitive conditions}
        \sqrt{n} \sup_{|v|\leq V}\big\vert \widehat{\Gamma}_n(\x+v) - \widehat{\Gamma}_n(\x) - \bar{\Gamma}_n(\x+v) +\bar{\Gamma}_n(\x)\big\vert \leq V O_{\P}(1) + o_{\P}(a_n^{-1})
    \end{equation}
    uniformly over $V\in (0, 2\delta]$.
\end{lemma}
\begin{proof}
    By Theorem 4.3 of \cite{Cattaneo-Farrell-Feng_2020_AoS},
    \begin{equation*}
        \sup_{(x,\mathbf{a})\in [0,1]^{1+d}} \vert \widehat{\mu}_n(x,\mathbf{a})-\mu_0(x,\mathbf{a})\vert = O_{\P}\bigg( \sqrt{\frac{k_n\log n}{n}} + k_n^{-\frac{m}{d+1}}\bigg)
    \end{equation*}
    and by Theorem 1 of \citet*{Cattaneo-Chandak-Jansson-Ma_2024},
    \begin{equation*}
        \sup_{(x,\mathbf{a})\in [0,1]^{1+d}}\vert\widehat{f}_{X|\mathbf{A},n}(x|\mathbf{a})- f_{X|\mathbf{A}}(x|\mathbf{a})\vert = O_{\P}\bigg( \sqrt{\frac{\log n}{n h^{1+d}}} + h^{\mathfrak{p}_1}\bigg).
    \end{equation*}
    Also, one can show
    \begin{equation*}
        \sup_{x\in [0,1]}\vert\widehat{f}_n(x) -f_0(x)\vert = O_{\P}\bigg( \sqrt{\frac{\log n}{n b}} + b^{\mathfrak{p}_1}\bigg).
    \end{equation*}
    Then, with the specified rate of $k_n,h,b$,
    \begin{equation*}
        \sup_{(x,\mathbf{a})\in [0,1]^{1+d}} \vert \widehat{\mu}_n(x,\mathbf{a})-\mu_0(x,\mathbf{a})\vert = O_{\P}\Big( \tau_n^{\frac{m}{2m+d+1}} \Big),
    \end{equation*}
    \begin{equation*}
        \sup_{(x,\mathbf{a})\in [0,1]^{1+d}} \vert \widehat{g}_n(x,\mathbf{a})-g_0(x,\mathbf{a})\vert = O_{\P}\Big( \tau_n^{\frac{\mathfrak{p}_1}{2\mathfrak{p}_1+d+1}} \Big).
    \end{equation*}
    Since $\min\{\frac{2m}{2m+d+1},\frac{2\mathfrak{p}_1}{2\mathfrak{p}_1+d+1} \} > \frac{1}{2\q+1}$, it follows $a_n\frac{1}{n}\sum_{i=1}^n|\widehat{\mu}_n(X_i,\mathbf{A}_i)-\mu_0(X_i,\mathbf{A}_i)|^2=o_{\P}(1)$, $a_n\frac{1}{n^2}\sum_{i=1}^n\sum_{j=1}^n|\widehat{\mu}_n(X_i,\mathbf{A}_j)-\mu_0(X_i,\mathbf{A}_j)|^2=o_{\P}(1)$, and $a_n\frac{1}{n}\sum_{i=1}^n\varepsilon_i^2|\widehat{g}_n(X_i,\mathbf{A}_i)-g_0(X_i,\mathbf{A}_i)|^2=o_{\P}(1)$. Also, by $\frac{m}{2m+d+1}+\frac{\mathfrak{p}_1}{2\mathfrak{p}_1+d+1} \geq \frac{1}{2}$,
    \begin{align}\label{Appendix: Primitive Condition Covariate Adjusted Regression sqrt N rate}
        \sup_{(x,\mathbf{a})\in [0,1]^{1+d}} \vert \widehat{\mu}_n(x,\mathbf{a})-\mu_0(x,\mathbf{a})\vert \sup_{(x,\mathbf{a})\in [0,1]^{1+d}} \vert \widehat{g}_n(x,\mathbf{a})-g_0(x,\mathbf{a})\vert = O_{\P}\Big( n^{-1/2}\Big).
    \end{align}
    
    Decompose $\widehat{\gamma}_n$ into $\widehat{\gamma}_{1,n}$ and $\widehat{\gamma}_{2,n}$ where
    \begin{equation*}
        \widehat{\gamma}_{1,n}(x;\mathbf{Z}) = \1(X\leq x) \frac{Y-\widehat{\mu}_n(X,\mathbf{A})}{\widehat{g}_n(X,\mathbf{A})},\quad \widehat{\gamma}_{2,n}(x;\mathbf{Z}) = \1(X\leq x) \frac{1}{n}\sum_{j=1}^n\widehat{\mu}_n(X,\mathbf{A}_j)
    \end{equation*}
    and let $\widehat{\Gamma}_{k,n}(x)=\frac{1}{n}\sum_{i=1}^n \widehat{\gamma}_{1,n}(x;\mathbf{Z}_i)$ for $k=1,2$. Define $\gamma_{k,0},\bar{\Gamma}_{k,n}$ $k=1,2$ in the same manner. 
    
    Letting $\tilde{\mathfrak{Z}}_n$ be the $\sigma$-field generated by $\tilde{\mathbf{Z}}_1,\dots,\tilde{\mathbf{Z}}_n$,
    \begin{align*}
        &\V\Big[ \widehat{\Gamma}_{1,n}(\x+v) - \widehat{\Gamma}_{1,n}(\x) - \bar{\Gamma}_{1,n}(\x+v) + \bar{\Gamma}_{1,n}(\x) \Big\vert \tilde{\mathfrak{Z}}_n\Big] \\
        &\leq n^{-1} \E\left[\left. \big(\1(X\leq\x+v) -\1(X\leq \x) \big)^2 \bigg(\frac{Y-\widehat{\mu}_n(X,\mathbf{A})}{\widehat{g}_n(X,\mathbf{A})} - \frac{Y-\mu_0(X,\mathbf{A})}{g_0(X,\mathbf{A})}\bigg)^2 \right\vert\tilde{\mathfrak{Z}}_n\right] \\
        &\leq n^{-1} C \E\big[\big\vert\1(X\leq\x+v) -\1(X\leq \x) \big\vert \varepsilon^2\big] \sup_{\vert x-\x\vert\leq \vert v\vert}\sup_{\mathbf{a}\in [0,1]^d} \vert \widehat{g}_n(x,\mathbf{a})-g_0(x,\mathbf{a})\vert^2 \\
        &\quad + n^{-1} C \E\big[\big\vert\1(X\leq\x+v) -\1(X\leq \x) \big\vert \big]\sup_{\vert x-\x\vert\leq \vert v\vert}\sup_{\mathbf{a}\in [0,1]^d} \vert \widehat{\mu}_n(x,\mathbf{a})-\mu_0(x,\mathbf{a})\vert^2\\
        &= o_{\P}\big( (n a_n)^{-1} \big)
    \end{align*}
    where the inequalities hold with probability one. 
    
    Note $\widehat{\Gamma}_{2,n}(x) = \frac{1}{n^2}\sum_{1\leq i\neq j\leq n} \1(X_i\leq x)\widehat{\mu}_n(X_i,\mathbf{A}_j) + O_{\P}(n^{-1})$, and
    \begin{align*}
        &\V\left[ \left.\frac{1}{n(n-1)}\sum_{1\leq i\neq j\leq n} \big(\1(X_i\leq \x+v)-\1(X_i\leq \x)\big)\big(\widehat{\mu}_n(X_i,\mathbf{A}_j) -\mu_0(X_i,\mathbf{A}_j)\big) \right\vert \tilde{\mathfrak{Z}}_n\right] \\
        &\leq n^{-1}C \vert v\vert \sup_{\vert x-\x\vert\leq \vert v\vert}\sup_{\mathbf{a}\in [0,1]^d} \vert \widehat{\mu}_n(x,\mathbf{a})-\mu_0(x,\mathbf{a})\vert^2 \\
        &\quad + C n^{-2}\vert v\vert \sup_{\vert x-\x\vert\leq \vert v\vert}\sup_{\mathbf{a}\in [0,1]^d} \vert \widehat{\mu}_n(x,\mathbf{a})-\mu_0(x,\mathbf{a})\vert^2 =  o_{\P}\big((n a_n)^{-1}\big)
    \end{align*}
    where we use Hoeffding decomposition and the inequality holds with probability one.

    To complete the proof, it suffices to show that there is a sequence of random variables $A_n'=O_{\P}(1)$ such that for $V \in (0,a_n\delta]$,
    \begin{equation*}
        \sqrt{n } \sup_{\vert v\vert \leq V}\left\vert  \E\big[\widehat{\Gamma}_{n}(\x+v) - \widehat{\Gamma}_{n}(\x) - \bar{\Gamma}_{n}(\x+v) + \bar{\Gamma}_{n}(\x)\big\vert \tilde{\mathfrak{Z}}_n\big] \right\vert \leq  A_n' V.
    \end{equation*}
    With $f_{\mathbf{A}}(\mathbf{a})$ denoting the Lebesgue density of $\mathbf{A}$,
    \begin{align*}
        &\E\big[\widehat{\Gamma}_{n}(\x+v) - \widehat{\Gamma}_{n}(\x) \big\vert \tilde{\mathfrak{Z}}_n\big] \\
        &= \int \int \big(\1(u\leq \x+v)-\1(u\leq \x)\big)  \frac{g_0(u,\mathbf{a})}{\widehat{g}_n(u,\mathbf{a})}\big[\mu_0(u,\mathbf{a})-\widehat{\mu}_n(u,\mathbf{a})\big]f_0(u)du f_{\mathbf{A}}(\mathbf{a})d\mathbf{a} \\
        &\qquad+ \int \big(\1(u\leq \x+v)-\1(u\leq \x)\big)\int\widehat{\mu}_n(u,\mathbf{a})f_{\mathbf{A}}(\mathbf{a})d\mathbf{a} du\\
        & \E\big[\bar{\Gamma}_{n}(\x+v) - \bar{\Gamma}_{n}(\x)\big\vert \tilde{\mathfrak{Z}}_n\big]\\
        &=\int\big(\1(u\leq \x+v)-\1(u\leq \x)\big)\int \mu_0(u,\mathbf{a})f_{\mathbf{A}}(\mathbf{a})d\mathbf{a}f_0(u)du.
    \end{align*}
    Then, for $v\geq 0$,
    \begin{align*}
        &\E\big[\widehat{\Gamma}_{n}(\x+v) - \widehat{\Gamma}_{n}(\x) - \bar{\Gamma}_{n}(\x+v) + \bar{\Gamma}_{n}(\x)\big\vert \tilde{\mathfrak{Z}}_n\big]\\
        &= \int_{\x\leq u\leq \x +v} \big[\widehat{\mu}_n(u,\mathbf{a})-\mu_0(u,\mathbf{a})\big] \Big[1-\frac{g_0(u,\mathbf{a})}{\widehat{g}_n(u,\mathbf{a})}\Big]f_0(u) f_{\mathbf{A}}(\mathbf{a})d(u,\mathbf{a}).
    \end{align*}
    A similar expression holds for $v<0$. Then, for some fixed $C>0$,
    \begin{align*}
        &\left\vert \E\big[\widehat{\Gamma}_{n}(\x+v) - \widehat{\Gamma}_{n}(\x) - \bar{\Gamma}_{n}(\x+v) + \bar{\Gamma}_{n}(\x)\big\vert \tilde{\mathfrak{Z}}_n\big]\right\vert \\
        &\leq C \vert v\vert \sup_{\vert x-\x\vert\leq \vert v\vert}\sup_{\mathbf{a}\in [0,1]^d} \vert \widehat{\mu}_n(x,\mathbf{a})-\mu_0(x,\mathbf{a})\vert \sup_{\vert x-\x\vert\leq \vert v\vert}\sup_{\mathbf{a}\in [0,1]^d} \vert \widehat{g}_n(x,\mathbf{a})-g_0(x,\mathbf{a})\vert
    \end{align*}
    and the desired result follows from \eqref{Appendix: Primitive Condition Covariate Adjusted Regression sqrt N rate}.
\end{proof}
\paragraph*{Remark}
Using \eqref{Appendix eq: result - monotone regression with covariates primitive conditions}, one can verify the first part of Assumption \ref{M-Assumption D} i.e.,\ $\widehat{G}_{\x,n}(1;\eta_n)=O_{\P}(1)$. The second part is easy to verify using standard empirical process theory arguments.

\section{Additional example: monotone density function with conditionally independent right-censoring}

\label{Appendix Section: Example monotone density conditionally independent}
We consider the problem of estimating the density of a non-negative, continuously distributed random variable with censoring. We use the same notation as in Section \ref{M-[Subsection] Examples: Monotone Density Estimation} of the main paper. Relative to Section \ref{M-[Subsection] Examples: Monotone Density Estimation}, we consider the additional complication of censoring being informative about the ``survival time'' $X$ i.e.,\ $X\not\indep C$. With covariates $\mathbf{A}$, we consider the setting of censoring at random: $X \indep C|\mathbf{A}$. See \cite{vanderlaan-Robins_2003_Book,Zeng_2004_AoS} and references therein for existing analysis of this problem. We have
\begin{equation*}
        \gamma_0(x;\mathbf{Z}) = F_0(x|\mathbf{A}) + S_0(x|\mathbf{A}) \bigg[ \frac{\Delta\1(\check{X}\leq x)}{S_0(\check{X}|\mathbf{A})G_0(\check{X}|\mathbf{A})} - \int_0^{\check{X}\land x}\frac{\Lambda_0(du|\mathbf{A})}{S_0(u|\mathbf{A})G_0(u|\mathbf{A})}\bigg]
\end{equation*}
where $F_0(x|A)=1-S_0(x|A)$, $S_0(x|\mathbf{A})=\P[X>x|\mathbf{A}]$, $G_0(c|\mathbf{A})=\P[C>c|\mathbf{A}]$, and $\Lambda_0(x|\mathbf{A}) = \int_0^{x} \frac{f_0(u|\mathbf{A})}{S_0(u|\mathbf{A})}du$ with $f_0$ being the Lebesgue density of $X$.
Denote by $\widehat{S}_n(\cdot|\cdot)$, $\widehat{G}_n(\cdot|\cdot)$, $\widehat{\Lambda}_n(\cdot|\cdot)$ preliminary estimates of $S_0,G_0,\Lambda_0$, respectively.

\paragraph*{Assumption \thesection}\label{Appendix Assumption: Monotone Density with Conditional Independent Censoring}
Let $\mathfrak{S}_n$, $\mathfrak{G}_n,\mathfrak{L}_n$ be sequences of function classes that contain $S_0(\cdot|\cdot)$, $G_0(\cdot|\cdot)$, $\Lambda_0(\cdot|\cdot)$, respectively. 
\begin{enumerate}[label=\normalfont(\roman*),noitemsep,itemindent=*]
    \item $\x$ is in the interior of $I=[0,u_0]$, $X\indep C |A$, and $\theta_0=f_0$ satisfies Assumption \ref{M-Assumption A: theta0}.
    
    \item \label{Appendix Example: monotone density conditionally independent function classes} There exist $c,c_1,c_2>0,\rho_{\gamma}\in (0,2)$such that for $n\geq 1$, for any $S\in\mathfrak{S}_n$, $G\in \mathfrak{G}_n$, and $\Lambda \in \mathfrak{L}_n$, the following hold: $\log N_U(\varepsilon,\{S(x|\cdot):x\in I\}) \leq c \varepsilon^{-\rho_{\gamma}}$ for $\varepsilon\in (0,1)$, where $N_U$ is as defined in Section 4.2 of the main paper, and $c_1\leq S(x|\mathbf{A})\leq c_2$, $c_1\leq G(x|\mathbf{A})\leq c_2$ for $x\in I$, and $\Lambda(u_0|\mathbf{A})\leq c_2$ with probability one. 

    \item \label{Appendix Example: monotone density conditionally independent high-level}
    There exist $\delta>0, \beta_{\gamma} \in [1/2,2)$ such that for $V \in (0,2\delta]$, $\sqrt{n a_n}\sup_{|v|\leq V}|\widehat{\Gamma}_n(\x+v)-\widehat{\Gamma}_n(\x)-\bar{\Gamma}_n(\x+v)+\bar{\Gamma}_n(\x)|\leq o_{\P}(1) +V^{\beta_{\gamma}}o_{\P}(a_n^{\beta_{\gamma}})$ where $o_{\P}$ terms do not depend on $V$.
    
    \item With probability approaching one, $\widehat{S}_n\in\mathfrak{S}_n$, $\widehat{G}_{n}\in\mathfrak{G}_n$, $\widehat{\Lambda}_n\in\mathfrak{L}_n$. For $(\widehat{h}_n,h_0) \in \{(\widehat{S}_n,S_0),(\widehat{G}_n,G_0),(\widehat{\Lambda}_n,\Lambda_0) \}$,
    \begin{equation*}
        a_n\frac{1}{n}\sum_{i=1}^n\sup_{x\in I}|\widehat{h}_n(x|\mathbf{A}_i)-h_0(x|\mathbf{A}_i)|^2=o_{\P}(1).
    \end{equation*}
 
    \item The conditional distribution of $X$ given $\mathbf{A}$ has bounded Lebesgue density $f_{X|\mathbf{A}}$, $\E[\frac{f_{X|\mathbf{A}}(\x|\mathbf{A})}{G_0(\x|\mathbf{A})}]>0$, and there are real-valued functions $B,\omega$ such that $\E[B(\mathbf{A})]<\infty$, $\lim_{\eta\downarrow0}\omega(\eta)=0$, and for $|x-\x|$ sufficiently small, $|\frac{f_{X|\mathbf{A}}(x|\mathbf{A})}{S_0(x|\mathbf{A})G_0(x|\mathbf{A})}-\frac{f_{X|\mathbf{A}}(\x|\mathbf{A})}{S_0(\x|\mathbf{A})G_0(\x|\mathbf{A})}|\leq \omega(|x-\x|)B(\mathbf{A})$. 
\end{enumerate}
The condition \ref{Appendix Example: monotone density conditionally independent high-level} is high-level, and there are a few different approaches to verify them. See \cite{Westling-Carone_2020_AoS} for details.
\begin{corollary}\label{Appendix Corollary: Monotone Density with Conditional Independent Censoring}
    Under Assumption \ref{Appendix Assumption: Monotone Density with Conditional Independent Censoring}, Assumptions \ref{M-Assumption A} and \ref{M-Assumption B} hold with
    \begin{equation*}
        \widehat{\Gamma}_n(x) = \frac{1}{n}\sum_{i=1}^n \widehat{\gamma}_n(x;\mathbf{Z}_i),\quad  \widehat{\Gamma}_n^{*}(x) = \frac{1}{n}\sum_{i=1}^nW_{i,n} \widehat{\gamma}_n(x;\mathbf{Z}_i)
    \end{equation*}
    \[ \widehat{\gamma}_n(x;\mathbf{Z})
  = \widehat{F}_n(x|\mathbf{A}) + \widehat{S}_n(x|\mathbf{A})
    \bigg[\frac{\Delta\1(\check{X}\leq x)}{\widehat{S}_n(\check{X}|\mathbf{A})\widehat{G}_n(\check{X}|\mathbf{A})} 
    - \int_0^{\check{X}\land x} \frac{\widehat{\Lambda}_n(du|\mathbf{A})}{\widehat{S}_n(u|\mathbf{A})\widehat{G}_n(u|\mathbf{A})}\bigg],
\]
where $\widehat{F}_n=1-\widehat{S}_n$, 
    \begin{equation*}
        \widehat{\Phi}_n(x)=\widehat{\Phi}_n^*(x)=x,\quad \widehat{u}_n=\widehat{u}_n^*=u_0,
    \end{equation*}
    \begin{equation*}
        \mathcal{C}_{\x}(s,t) = \E\Big[\frac{f_{X|\mathbf{A}}(\x|\mathbf{A})}{G_0(\x|\mathbf{A})}\Big] (|s| \land |t|) \1(\sign(s)=\sign(t)), \quad \mathcal{D}_{\q}(\x) = \frac{\partial^\q f_0(\x)}{(\q+1)!}. 
\end{equation*}
\end{corollary}

\subsection{Proof of Corollary \ref{Appendix Corollary: Monotone Density with Conditional Independent Censoring}}
\label{Appendix Section: verification monotone density conditionally independent cernsoring}
In this example, $\widehat{\Phi}_n(x)=\widehat{\Phi}_n^*(x)=x=\Phi_0(x)$. Assumptions \ref{M-Assumption A} and \ref{M-Assumption E: iid}-\ref{M-Assumption E: bootstrap weights} follow from the hypothesis.
\newline
\ref{M-Assumption E: Gamma hat}:
In this example, $\check{\Gamma}_n=\widehat{\Gamma}_n$. For $x\in I$,
\begin{align*}
    &\widehat{\gamma}_n(x;\mathbf{Z}) -\gamma_0(x;\mathbf{Z}) = \widehat{F}_n(x|\mathbf{A}) - F_0(x|\mathbf{A}) \\
    &\quad + \big[\widehat{S}_n(x|\mathbf{A})-S_0(x|\mathbf{A})\big]\bigg[\frac{\Delta\1(\check{X}\leq x)}{\widehat{S}_n(\check{X}|\mathbf{A})\widehat{G}_n(\check{X}|\mathbf{A})} 
    - \int_0^{\check{X}\land x} \frac{\widehat{\Lambda}_n(du|\mathbf{A})}{\widehat{S}_n(u|\mathbf{A})\widehat{G}_n(u|\mathbf{A})}\bigg] \\
    &\quad + S_0(x|\mathbf{A})\Delta \1(\check{X}\leq x) \Big[\big(\widehat{S}_n(\check{X}|\mathbf{A})\widehat{G}_n(\check{X}|\mathbf{A})\big)^{-1} - \big(S_0(\check{X}|\mathbf{A})G_0(\check{X}|\mathbf{A})\big)^{-1} \Big] \\
    &\quad - S_0(x|\mathbf{A}) \bigg[\int_0^{\check{X}\land x} \frac{\widehat{\Lambda}_n(du|\mathbf{A})}{\widehat{S}_n(u|\mathbf{A})\widehat{G}_n(u|\mathbf{A})} - \int_0^{\check{X}\land x}\frac{\Lambda_0(du|\mathbf{A})}{S_0(u|\mathbf{A})G_0(u|\mathbf{A})} \bigg]
\end{align*}
and
\begin{align*}
    &\left\vert\int_0^{\check{X}\land x} \frac{\widehat{\Lambda}_n(du|\mathbf{A})}{\widehat{S}_n(u|\mathbf{A})\widehat{G}_n(u|\mathbf{A})} - \int_0^{\check{X}\land x}\frac{\Lambda_0(du|\mathbf{A})}{S_0(u|\mathbf{A})G_0(u|\mathbf{A})} \right\vert \\
    &\leq \sup_{x'\in I} \Big\vert \big(\widehat{S}_n(x'|\mathbf{A})\widehat{G}_n(x'|\mathbf{A})\big)^{-1} - \big(S_0(x'|\mathbf{A})G_0(x'|\mathbf{A})\big)^{-1} \Big\vert \widehat{\Lambda}_n(u_0|\mathbf{A}) \\
    & \quad + \left\vert \frac{\widehat{\Lambda}_n(\check{X}\land x|\mathbf{A})-\Lambda_0(\check{X}\land x|\mathbf{A})}{S_0(\check{X}\land x|\mathbf{A})G_0(\check{X}\land x|\mathbf{A})} \right\vert + \left\vert \frac{\widehat{\Lambda}_n(0|\mathbf{A})-\Lambda_0(0|\mathbf{A})}{S_0(0|\mathbf{A})G_0(0|\mathbf{A})} \right\vert \\
    &\quad + \left\vert \int_0^{\check{X}\land x} [\widehat{\Lambda}_n(u|\mathbf{A})-\Lambda_0(u|\mathbf{A})]J_0(du|\mathbf{A})  \right\vert
\end{align*}
using integration by parts, where $J_0(u|\mathbf{a}) = [S_0(u|\mathbf{a})G_0(u|\mathbf{a})]^{-1}$. Thus, there is a fixed $C>0$ such that
\begin{align*}
    |\widehat{\gamma}_n(x;\mathbf{Z}) -\gamma_0(x;\mathbf{Z})|&\leq C \Big[ \sup_{x\in I}|\widehat{S}_n(x|\mathbf{A})-S_0(x|\mathbf{A})| + \sup_{x\in I}|\widehat{G}_n(x|\mathbf{A})-G_0(x|\mathbf{A})|\\
    &\qquad +\sup_{x\in I}|\widehat{\Lambda}_n(x|\mathbf{A})-\Lambda_0(x|\mathbf{A})| \Big],
\end{align*}
From the hypothesis, 
\begin{equation*}
    \frac{1}{n}\sum_{i=1}^n\sup_{x\in I}|\widehat{\gamma}_n(x;\mathbf{Z}_i)-\gamma_0(x;\mathbf{Z}_i)|^2=o_{\P}(1),\ a_n \frac{1}{n}\sum_{i=1}^n\sup_{x\in I_{\x}^{\delta}}\big|\widehat{\gamma}_n(x;\mathbf{Z}_i)  -\gamma_0(x;\mathbf{Z}_i) \big|^2 =o_{\P}(1)
\end{equation*}
follow.

For uniform covering numbers, it suffices to show that each of $\{S(x|\cdot):x\in I\}$, $\{\1(\cdot\leq x):x\in I\}$, and $\{\int_0^{\cdot\land x} \frac{\Lambda(du|\cdot)}{S(u|\cdot)G(u|\cdot)} : x\in I\}$ has an appropriate bound on the uniform covering number by Lemma 5.1 of \cite{vanderVaart-vanderLaan_2006_IJB} (see examples after the lemma). For $\{\int_0^{\cdot\land x} \frac{\Lambda(du|\cdot)}{S(u|\cdot)G(u|\cdot)} : x\in I\}$ with $(S,G,\Lambda)\in \mathfrak{S}_n\times \mathfrak{G}_n\times\mathfrak{L}_n$, the mapping $x\mapsto \int_0^x \frac{\Lambda(du|\cdot)}{S(u|\cdot)G(u|\cdot)}$ is monotone (by the non-decreasing property of $\Lambda$ and $S,G\geq c_1 >0 $) and Lemma 2.6.18 of \cite{vanderVaart-Wellner_1996_Book} implies the desired result.

There is a fixed $C>0$ such that for $x\in I_{\x}^{\eta}$,
\begin{align*}
    |\gamma_0(x;\mathbf{Z}) -\gamma_0(\x;\mathbf{Z})|&\leq C \Big[|S_0(x|\mathbf{A})-S_0(\x|\mathbf{A})| + C\Delta |\1(\check{X}\leq x) - \1(\check{X}\leq \x)| \Big]
\end{align*}
and using $1-S_0(x|\cdot) =\int_0^x f_{X|\mathbf{A}}(u|\cdot) du$ with $f_{X|\mathbf{A}}$ being bounded, we can take
\begin{equation*}
    \bar{D}_{\gamma}^{\eta}(\mathbf{Z}) =  C\Delta \1(\x-\eta \leq \check{X}\leq \x +\eta) + C\eta,
\end{equation*}
which satisfies the desired bound condition.
\newline 
\ref{M-Assumption E: Phi hat} trivially holds since $\widehat{\Phi}_n(x)=\widehat{\Phi}_n^*(x)=x$.

\noindent \ref{M-Assumption E: covariance kernel}:
We have
\begin{equation*}
    \psi_{\x}(v;\mathbf{Z}) = S_0(\x|\mathbf{A})\frac{(\1(\check{X}\leq \x +v)-\1(\check{X}\leq \x))\Delta}{S_0(\check{X}|\mathbf{A})G_0(\check{X}|\mathbf{A})}  + O(|v|)
\end{equation*}
and the first display follows as in the independent censoring case. For the covariance kernel,
\begin{align*}
    \E[\psi_{\x_n}(s\eta_n;\mathbf{Z})\psi_{\x_n}&(t\eta_n;\mathbf{Z})] = \E\bigg[S_0(\x_n|\mathbf{A})^2 \int_{\x_n}^{\x_n +\eta_n (s\land t)} \frac{f_{X|\mathbf{A}}(u|\mathbf{A})}{S_0(u|\mathbf{A})^2G_0(u|\mathbf{A})}du \1(s,t>0) \bigg]\\
    &+ \E\bigg[S_0(\x_n|\mathbf{A})^2 \int_{\x_n +\eta_n (s\lor t)}^{\x_n} \frac{f_{X|\mathbf{A}}(u|\mathbf{A})}{S_0(u|\mathbf{A})^2G_0(u|\mathbf{A})}du \1(s,t<0) \bigg] + O(\eta_n^2)
\end{align*}
and $\eta_n^{-1}\E[\psi_{\x_n}(s\eta_n;\mathbf{Z})\psi_{\x_n}(t\eta_n;\mathbf{Z})]$ converges to $\E[\frac{f_{X|\mathbf{A}}(\x|\mathbf{A})}{G_0(\x|\mathbf{A})}] (|s|\land|t|) \1(\sign(s)=\sign(t))$.

\noindent
\ref{M-Assumption B: uhat}, \ref{M-Assumption B: Phi uhat}, and \ref{M-Assumption B: jump of Phi hat} hold since in this example, $\widehat{u}_n=\widehat{u}_n^*=u_0$ and $\widehat{\Phi}_n,\widehat{\Phi}_n^*$ are the identity map.

\section{Additional example: monotone hazard function}

Let $X$ be a non-negative random variable, $f_0$ be its Lebesgue density, and $S_0(x)=\P[X>x]$ be its survival function. We consider the parameter of estimating the hazard function of $X$, $\theta_0(\x)=f_0(\x)/S_0(\x)$, with possible right-censoring as in the monotone density function example. Observations $\mathbf{Z}_1,\dots,\mathbf{Z}_n$ come from a random sample of $\mathbf{Z}=(\check{X},\Delta)'$ where $\check{X}=\min\{X,C\}$ and $\Delta=\1(X\leq C)$, $C$ being a random censoring time. As pointed out by \cite{Westling-Carone_2020_AoS}, with strictly increasing $\Phi_0$ with $\Phi_0(0)=0$, the function $\Gamma_0$ takes the form $\Gamma_0(x) = \int_0^x \frac{f_0(u)}{S_0(u)}\Phi_0(du)$, and by taking $\Phi_0(x)=\int_0^xS_0(u)du$, $\Gamma_0(x) =F_0(x)=\P[X\leq x]$. Since $\Gamma_0$ is identical to the monotone density case with the choice $\Phi_0=\int_0^xS_0(u)du$, we can leverage the analysis for the monotone density. The interval $I$ equals $[0,u_0^{\mathtt{MD}}]$ where $u_0^{\mathtt{MD}}$ is $u_0$ in the monotone density example. The $u_0$ for the monotone hazard function estimation is $u_0=\Phi_0(u_0^{\mathtt{MD}})$.

Consider the case of completely random censoring i.e.,\ $X\indep C$. As in the setup for Corollary \ref{M-[Corollary] monotone density independent right-censoring}, let $\widehat{S}_n(x)$ be the Kaplan-Meier estimator for $S_0(x)=1-F_0(x)=\P[X> x]$, $\widehat{F}_n=1-\widehat{S}_n$, and $\widehat{G}_n$ be the Kaplan-Meier estimator for $G_0(x)=\P[C>x]$. Also,
\begin{equation*}
    \gamma_0(x;\mathbf{Z}) = F_0(x) + S_0(x)\left[\frac{\Delta\1(\check{X}\leq x)}{S_0(\check{X})G_0(\check{X})} - \int_0^{\check{X}\land x}\frac{\Lambda_0(du)}{S_0(u)G_0(u)} \right]
\end{equation*}
and $\phi_0(x;\mathbf{Z}) =x-\int_0^x \gamma_0(u;\mathbf{Z})du$.
\begin{corollary}\label{Appendix Corollary: Monotone Hazard Function Independent Right-Censoring}
    Suppose that the hypothesis of Corollary \ref{M-[Corollary] monotone density independent right-censoring} and Assumption BW hold. Then, Assumptions \ref{M-Assumption A} and \ref{M-Assumption B} hold with
    \begin{equation*}
        \widehat{\Gamma}_n(x)=1-\widehat{S}_n(x),\quad \widehat{\Gamma}_n^*(x)=\frac{1}{n}\sum_{i=1}^nW_{i,n}\widehat{\gamma}_n(x;\mathbf{Z}_i),
    \end{equation*}
    \begin{equation*}
        \widehat{\gamma}_n(x;\mathbf{Z}) = \widehat{F}_n(x) + \widehat{S}_n(x)\left[\frac{\Delta\1(\check{X}\leq x)}{\widehat{S}_n(\check{X})\widehat{G}_n(\check{X})} - \int_0^{\check{X}\land x} \frac{\widehat{\Lambda}_n(du)}{\widehat{S}_n(u)\widehat{G}_n(u)}\right],
    \end{equation*}
    \begin{equation*}
        \widehat{\Phi}_n(x)= \int_0^x \widehat{F}_n(u) du,\quad \widehat{\Phi}_n^*(x)=\int_0^x [1-\widehat{\Gamma}_n^*(u)]du = \frac{1}{n}\sum_{i=1}^n W_{i,n} \widehat{\phi}_n(x;\mathbf{Z}_i),
    \end{equation*}
    \begin{equation*}
        \widehat{\phi}_n(x;\mathbf{Z}) = x-\int_0^x \widehat{\gamma}_n(u;\mathbf{Z})du,\quad \widehat{u}_n=\widehat{u}_n^*=\widehat{\Phi}_n(u_0^{\mathtt{MD}}),
    \end{equation*}
    \begin{equation*}
        \mathcal{C}_{\x}(s,t) = \frac{f_0(\x)}{G_0(\x)} (|s| \land |t|) \1(\sign(s)=\sign(t)),\quad \mathcal{D}_{\q}(\x) = \frac{S_0(\x)\partial^\q f_0(\x)}{(\q+1)!}.
\end{equation*}
\end{corollary}

\subsection{Proof of Corollary \ref{Appendix Corollary: Monotone Hazard Function Independent Right-Censoring}}

We use the same $\widehat{\gamma}_n$ function and assumptions as in the monotone density setting. Also, the covariance kernels are the same as in the monotone density case. We focus on \ref{M-Assumption E: Phi hat}-\ref{M-Assumption E: covariance kernel} and \ref{M-Assumption B: uhat}-\ref{M-Assumption B: jump of Phi hat}.

\noindent \ref{M-Assumption E: Phi hat}:
Since
\begin{equation*}
    \widehat{\phi}_n(x;\mathbf{Z}) -\phi_0(x;\mathbf{Z})= -\int_0^x [\widehat{\gamma}_n(u;\mathbf{Z})-\gamma_0(u;\mathbf{Z})]du,
\end{equation*}
\begin{equation*}
    \frac{1}{n}\sum_{i=1}^n\sup_{x\in I}|\widehat{\phi}_n(x;\mathbf{Z}_i)-\phi_0(x;\mathbf{Z}_i)|^2=o_{\P}(1),\ a_n\frac{1}{n}\sum_{i=1}^n \sup_{x\in I_{\x}^{\delta}}|\widehat{\phi}_n(x;\mathbf{Z}_i)-\phi_0(x;\mathbf{Z}_i)|^2=o_{\P}(1)
\end{equation*}
follow from $a_n\frac{1}{n}\sum_{i=1}^n\sup_{x\in I}\vert\widehat{\gamma}_n(x;\mathbf{Z}_i) - \gamma_0(x;\mathbf{Z}_i)\vert^2$, which was verified in Section \ref{Appendix subsection: proof of Theorem monotone density independent right-censoring}.
To check $\sup_{x\in I}|\widehat{\Phi}_n(x)-\Phi_0(x)|=o_{\P}(1)$, $\sup_{x\in I}|\frac{1}{n}\sum_{i=1}^n\phi_0(x;\mathbf{Z}_i)-\Phi_0(x)|=o_{\P}(1)$ follows from Glivenko-Cantelli, and
\begin{equation*}
    \sup_{x\in I}\bigg|\frac{1}{n}\sum_{i=1}^n \big[\widehat{\phi}_n(x;\mathbf{Z}_i)-\phi_0(x;\mathbf{Z}_i)\big]\bigg|\leq \frac{1}{n}\sum_{i=1}^n \sup_{x\in I}|\widehat{\gamma}_n(x;\mathbf{Z}_i)-\gamma_0(x;\mathbf{Z}_i)| = o_{\P}(1),
\end{equation*}
where the last equality follows from $\frac{1}{n}\sum_{i=1}^n\sup_{x\in I}|\widehat{\gamma}_n(x;\mathbf{Z}_i)-\gamma_0(x;\mathbf{Z}_i)|^2=o_{\P}(1)$. Now $\sup_{x\in I}|\widehat{\Phi}_n(x)-\Phi_0(x)|=o_{\P}(1)$ follows by the triangle inequality.

For $|v|\leq V$
\begin{align*}
    \big\vert\widehat{\Phi}_n(\x + v) -& \widehat{\Phi}_n(\x) - \bar{\Phi}_n(\x+v) + \bar{\Phi}_n(\x) \big\vert = \left\vert\int_{\x-v}^{\x+v} [\widehat{\Gamma}_n(u)-\bar{\Gamma}_n(u)]du \right\vert  \\
    &\leq 2\vert V\vert \Big( \sup_{\vert v\vert\leq \vert V\vert } \vert \widehat{\Gamma}_n(\x+v)-\widehat{\Gamma}_n(\x) - \bar{\Gamma}_n(\x+v) + \bar{\Gamma}_n(\x)\vert + \vert \widehat{\Gamma}_n(\x)-\bar{\Gamma}_n(\x)\vert \Big)
\end{align*}
Using the argument in Section \ref{Appendix subsection: proof of Theorem monotone density independent right-censoring}, we can bound $\sup_{\vert v\vert\leq \vert V\vert } \vert \widehat{\Gamma}_n(\x+v)-\widehat{\Gamma}_n(\x) - \bar{\Gamma}_n(\x+v) + \bar{\Gamma}_n(\x)\vert$. Then, $\sqrt{n} \vert \widehat{\Gamma}_n(\x)-\bar{\Gamma}_n(\x)\vert=O_{\P}(1)$ implies
\begin{equation*}
    \sqrt{n a_n}\sup_{|v|\leq V} \big\vert  \widehat{\Phi}_n(\x + v) -\widehat{\Phi}_n(\x) - \bar{\Phi}_n(\x+v) + \bar{\Phi}_n(\x)\big\vert \leq o_{\P}(1) + V O_{\P}\big(\sqrt{a_n}\big)
\end{equation*}
uniformly over $V\in (0,2\delta]$.
Theorem 1 of \cite{Lo-Singh_1986_PTRF} implies $\sqrt{n a_n}\sup_{x\in I}\vert \check{\Phi}_n(x ) - \check{\Phi}_n(\x ) - \bar{\Phi}_n(x)+\bar{\Phi}_n(\x )  \vert =o_{\P}(1)$.

The conditions on the uniform covering number hold because $\gamma_0$ and $\widehat{\gamma}_n$ are bounded (for $\widehat{\gamma}_n$, with probability approaching one) and thus $|\phi_0(x_1;\mathbf{Z})-\phi_0(x_2;\mathbf{Z})|\leq C |x_1-x_2|$
and $|\widehat{\phi}_n(x_1;\mathbf{Z})-\widehat{\phi}_n(x_2;\mathbf{Z})|\leq C |x_1-x_2|$ with probability approaching one. By this Lipschitz property, the condition on $\bar{D}_{\phi}^{\eta}(\mathbf{Z})$ also holds.

\noindent \ref{M-Assumption E: covariance kernel}:
Let $\psi_{\x}^{\mathtt{MD}}(v;\mathbf{Z})=\gamma_0(\x +v;\mathbf{Z})-\gamma_0(\x;\mathbf{Z}) - \theta_0(\x)v$ be the $\psi_{\x}$ function for the monotone density. Then, for $x$ sufficiently close to $\x$ and $|v|$ small enough,
\begin{align*}
    \psi_{x}(v;\mathbf{z}) &=\gamma_0(x+v;\mathbf{z})-\gamma_0(x;\mathbf{z}) -\theta_0(x)[\phi_0(x+v;\mathbf{z})-\phi_0(x;\mathbf{z})]\\
    &= \psi_{x}^{\mathtt{MD}}(v;\mathbf{Z}) + \theta_0(x) \int_{x}^{x+v} \gamma_0(u;\mathbf{Z}) du = \psi_{x}^{\mathtt{MD}}(v;\mathbf{Z}) + O(|v|).
\end{align*}
Then, the same argument as in the monotone density case implies the desired result.\smallskip
\newline
\ref{M-Assumption B: uhat} follows from consistency of $\widehat{\Phi}_n$ and $\widehat{\Phi}_n^*$.

\noindent \ref{M-Assumption B: Phi uhat}:
$\widehat{\Phi}_n(x)=\int_0^x \widehat{F}_n(u)du,\widehat{\Phi}_n^*(x)=\int_0^x 1-\widehat{\Gamma}_n^*(u)du$ are non-negative since $\widehat{F}_n\geq 0$ and $1-\widehat{\Gamma}_n^*\geq 0$ with probability approaching one. This property also implies the non-decreasing property as  $\widehat{\Phi}_n,\widehat{\Phi}_n^*$ are integrals. The continuity property also follows from the integral representation.
By definition, $\widehat{\Phi}_n(0)=0=\widehat{\Phi}_n^*(0)$ and $\widehat{\Phi}_n(u_0^{\mathtt{MD}})=\widehat{u}_n=\widehat{u}_n^*=\widehat{\Phi}_n^*(u_0^{\mathtt{MD}})$ with $I=[0,u_0^{\mathtt{MD}}]$.
The closedness of the range follows from continuity and $I$ being a compact interval.\smallskip
\newline
\ref{M-Assumption B: jump of Phi hat} follows from continuity of $\widehat{\Phi}_n$ and $\widehat{\Phi}_n^*$.

\section{Additional example: distribution function estimation with current status data}\label{Appendix Section: Example current status}
We consider the problem of estimating the cdf of $X$ at $\x$, $\theta_0(\x)=F_0(\x)$. Observations $\mathbf{Z}_1,\dots,\mathbf{Z}_n$ come from a random sample of $\mathbf{Z}=(\Delta,C,\mathbf{A}')'$ where $\Delta=\1(X\leq C)$, $C$ is a random censoring time, and $\mathbf{A}$ is a vector of covariates. In this example, we do not observe $\check{X}=X\land C$. Instead, we observe the censoring time and whether the observation was censored. This setup is often referred to as current status data.
Let $H_0(x)=\P[C\leq x]$ be the cdf of $C$. We can use $\Gamma_0(x) = \int_0^x F_0(u) H_0(du)$ and $\Phi_0(x)=H_0(x)$. The interval $I$ is the support of $X$ and $u_0=1$. We also assume $H_0$ admits a Lebesgue density $h_0$.
The structure of the estimation problem turns out to be identical to the one for the monotone regression example, and we can leverage the common structure.

\subsection{Independent right-censoring}
First we consider the case of completely at random censoring $X\indep C$. See \cite{Groeneboom-Wellner_1992_Book} for existing analysis. In this exaple, we do not use covariates $\mathbf{A}$. We set $\gamma_0(x;\mathbf{Z})=\Delta \1(C\leq x)$ and $\phi_0(x;\mathbf{Z})=\1(C\leq x)$. Note that if the notation is mapped by $(\Delta,C)\leftrightarrow (Y,X)$, then these functions are identical to those of the classical monotone regression problem (Corollary \ref{M-[Corollary] classical monotone regression}). Thus, the following result is identical to Corollary \ref{M-[Corollary] classical monotone regression}, up to notation and some changes due to boundedness of $\Delta$.
\begin{corollary}\label{Appendix Corollary: Current Status Data Completely Independent Case}
    Let $\varepsilon=\Delta-\E[\Delta|C]$ and $\x$ be an interior point of $I$. Suppose that Assumption BW holds, $\theta_0=F_0$ satisfies Assumption \ref{M-Assumption A: theta0}, the cdf $\Phi_0=H_0$ satisfies Assumption \ref{M-Assumption A: Phi0}, and $\sigma_0^2(x) = \E[\varepsilon^2|C=x]$ is continuous and positive at $\x$. Then  Assumptions \ref{M-Assumption A} and \ref{M-Assumption B} hold with
    \begin{equation*}
        \widehat{\Gamma}_n(x) = \frac{1}{n}\sum_{i=1}^n \Delta \1(C\leq x),\quad \widehat{\Gamma}_n^*(x) = \frac{1}{n}\sum_{i=1}^n W_{i,n}\Delta \1(C\leq x),
    \end{equation*}
    \begin{equation*}
        \widehat{\Phi}_n(x) = \frac{1}{n}\sum_{i=1}^n  \1(C\leq x),\quad \widehat{\Phi}_n^*(x) = \frac{1}{n}\sum_{i=1}^n W_{i,n}\1(C\leq x),\quad \widehat{u}_n=\widehat{u}_n^*=1,
    \end{equation*}
    \begin{equation*}
        \mathcal{C}_{\x}(s,t) = h_0(\x) \sigma_0^2(\x) (|s| \land |t|)\1(\sign(s)=\sign(t)),\qquad \mathcal{D}_{\q}(\x) = \frac{h_0(\x)\partial^{\q}F_0(\x)}{(\q+1)!}.
\end{equation*}
\end{corollary}

\subsection{Conditionally independent right-censoring}

We consider the case where right-censoring is conditionally independent i.e.,\ $X\indep C|\mathbf{A}$. \cite{vanderVaart-vanderLaan_2006_IJB} analyzed this example as well as settings with time-varying covariates. We are focusing on time-invariant covariates.
Define $F_0(C,\mathbf{A})=\E[\Delta|C,\mathbf{A}]$ and $g_0(C,\mathbf{A})=\frac{h_{C|\mathbf{A}}(C|\mathbf{A})}{h_0(C)}$ where $h_{C|\mathbf{A}}$ is a conditional Lebesgue density of $C$ given $\mathbf{A}$ and $h_0$ is a Lebesgue density of $C$. 
Let $\widehat{F}_n(c,\mathbf{a})$ and $\widehat{g}_n(c,\mathbf{a})$ be preliminary estimators for $F_0(c,\mathbf{a})$ and $g_0(c,\mathbf{a})$, respectively.
 
Identical to the censoring completely at random case, with appropriate changes in the notation (i.e.,\ $(\Delta, C)\leftrightarrow (Y,X)$), the setup is equivalent to that of the monotone regression with covariates. 

\paragraph*{Assumption \thesubsection} 
Let $\varepsilon=\Delta-\E[\Delta|C,\mathbf{A}]$, $\sigma_0^2(C,\mathbf{A})=\E[\varepsilon^2|C,\mathbf{A}]$, and $\delta >0$ be some fixed number.
\begin{enumerate}[label=\normalfont(\roman*),noitemsep,itemindent=*]
	\item $\E[\frac{\sigma_0^2(\x,\mathbf{A})}{g_0(\x,\mathbf{A})}]>0$.

    \item $h_{C|A}$ is bounded and $g_0$ is bounded away from zero.

    \item There exist random variables $A_n=o_{\P}(1)$ and $B_n=O_{\P}(a_n^{1/2})$ such that $$\sqrt{na_n}\sup_{|v|\leq V}|\widehat{\Gamma}_n(\x+v)-\widehat{\Gamma}_n(\x)-\Gamma_0(\x+v)+\Gamma_0(\x)|\leq A_n +B_nV,\qquad V\in (0,2\delta].$$ In addition, $$\frac{a_n}{n}\sum_{i=1}^n|\widehat{F}_n(C_i,\mathbf{A}_i)-F_0(C_i,\mathbf{A}_i)|^2=o_{\P}(1),$$ $$\frac{a_n}{n^2}\sum_{i=1}^n\sum_{j=1}^n|\widehat{F}_n(C_i,\mathbf{A}_j)-F_0(C_i,\mathbf{A}_j)|^2=o_{\P}(1),$$ and $$\frac{a_n}{n}\sum_{i=1}^n|\widehat{g}_n(C_i,\mathbf{A}_i)-g_0(C_i,\mathbf{A}_i)|^2=o_{\P}(1).$$

    \item $\E[\bar{F}(\mathbf{A})^2]<\infty$, where $$\bar{F}(\mathbf{A})=\sup_{|c-c'|\leq\delta}\frac{|F_0(c,\mathbf{A})-F_0(c',\mathbf{A})|}{|c-c'|}.$$ 

    \item $\E[\bar{\sigma}(\mathbf{A})^2]<\infty$, where, for some function $\omega$ with $\lim_{\eta\downarrow0}\omega(\eta)=0$, $$\left\vert\frac{\sigma_0^2(x,\mathbf{A})h_{0}(x)}{g_0(x,\mathbf{A})} -\frac{\sigma_0^2(\x,\mathbf{A})h_{0}(x)}{g_0(\x,\mathbf{A})}\right\vert \leq\omega(|x-\x|) \bar{\sigma}(\mathbf{A}),\qquad x\in I_{\x}^{\delta}.$$
\end{enumerate}
Note that $|\varepsilon|\leq 1$.

\begin{corollary}\label{Appendix Corollary: Current Status Data Conditionally Independence}
    Suppose that $\x$ is in the interior of $I$, $\theta_0$ satisfies \ref{M-Assumption A: theta0}, $\Phi_0$ satisfies \ref{M-Assumption A: Phi0}, and Assumption BW holds. If Assumption Assumption \thesubsection \ holds, Assumptions \ref{M-Assumption A} and \ref{M-Assumption B} hold with
    \begin{equation*}
        \widehat{\Gamma}_n(x)= \frac{1}{n}\sum_{i=1}^n \widehat{\gamma}_n(x;\mathbf{Z}_i),\quad  \widehat{\Gamma}_n^*(x)= \frac{1}{n}\sum_{i=1}^nW_{i,n} \widehat{\gamma}_n(x;\mathbf{Z}_i),
    \end{equation*}
    \begin{equation*}
    \widehat{\gamma}_n(x;\mathbf{Z})=\1(C\leq x)\bigg[\frac{\Delta-\widehat{F}_n(C,\mathbf{A})}{\widehat{g}_n(C,\mathbf{A})} + \frac{1}{n}\sum_{j=1}^n\widehat{F}_n(C,\mathbf{A}_j)\bigg],
\end{equation*}
\begin{equation*}
    \widehat{\Phi}_n(x) = \frac{1}{n}\sum_{i=1}^n \1(C_i\leq x),\quad \widehat{\Phi}_n^*(x) = \frac{1}{n}\sum_{i=1}^nW_{i,n} \1(C_i\leq x),\quad \widehat{u}_n=\widehat{u}_n^*=1,
\end{equation*}
\begin{equation*}
    \mathcal{C}_{\x}(s,t) = h_0(\x) \E\bigg[\frac{\sigma_0^2(\x,\mathbf{A})}{g_0(\x,\mathbf{A})}\bigg] (|s| \land |t|) \1(\sign(s)=\sign(t)),\qquad \mathcal{D}_{\q}(\x) = \frac{h_0(\x)\partial^{\q}F_0(\x)}{(\q+1)!}.
\end{equation*}
\end{corollary}

\subsection{Proof of Corollaries \ref{Appendix Corollary: Current Status Data Completely Independent Case} and \ref{Appendix Corollary: Current Status Data Conditionally Independence}}

As noted above, by mapping the notation $(\Delta, C)\leftrightarrow (Y,X)$, the arguments in Sections \ref{Appendix subsection: proof theorem classical monotone regression} and \ref{Appendix subsection: proof theorem monotone regression with covariates} directly apply to the current status estimators.

\section{Rule-of-thumb step size selection}

Here we develop a rule-of-thumb procedure to choose a step size for the bias-reduced numerical derivative estimator in the context of isotonic regression without covariates. Specifically, we consider the numerical derivative estimator
\begin{equation*}
    \widetilde{\mathcal{D}}_{j,n}^{\mathtt{BR}}(\x) = \epsilon_n^{-(j+1)}\sum_{k=1}^{\underline{\s}+1} \lambda_j^{\mathtt{BR}}(k)[\widehat{\Upsilon}_n(\x+c_k\epsilon_n)-\widehat{\Upsilon}_n(\x)]
\end{equation*}
with $\underline{\s}=3$, $c_1=1,c_2=-1,c_3=2,c_4=-2$. Then,
\begin{equation*}
    \lambda_1^{\mathtt{BR}}(1)=\frac{2}{3}=\lambda_1^{\mathtt{BR}}(2),\quad \lambda_1^{\mathtt{BR}}(3)=-\frac{1}{24}=\lambda_1^{\mathtt{BR}}(4),
\end{equation*}
\begin{equation*}
    \lambda_3^{\mathtt{BR}}(1)=-\frac{1}{6}=\lambda_3^{\mathtt{BR}}(2),\quad \lambda_3^{\mathtt{BR}}(3)=\frac{1}{24}=\lambda_3^{\mathtt{BR}}(4).
\end{equation*}
We use the (asymptotic) MSE-optimal step size discussed in the main paper. See also \ref{Appendix Section: higher-order expansion}. Yet, with the choice of $c_k$'s, part of the bias constant $\sum_{k=1}^{\underline{\s}+1}\lambda_j^{\mathtt{BR}}(k)c_k^{\underline{s}+2}$ equals zero, and we need to turn to the next leading term of the bias, which is
\begin{equation*}
    \epsilon_n^{\underline{\s}+2-j} \frac{\partial^{\underline{\s}+3}\Upsilon_0(\x)}{(\underline{\s}+3)!} \sum_{k=1}^{\underline{\s}+1}\lambda_j^{\mathtt{BR}}(k)c_k^{\underline{s}+3}.
\end{equation*}
Then, letting $\mathsf{B}_j^{\mathtt{BR}}(\x) = \frac{\partial^{6}\Upsilon_0(\x)}{6!} \sum_{k=1}^{4}\lambda_j^{\mathtt{BR}}(k)c_k^{6}$, the MSE-optimal step size is
\begin{equation*}
    \epsilon_{j,n}^{\mathtt{BR}} = \bigg(\frac{(2j+1)\mathsf{V}_j^{\mathtt{BR}}(\x) }{2(5-j)\mathsf{B}_j^{\mathtt{BR}}(\x)^2}\bigg)^{1/11} n^{-1/11}.
\end{equation*}
The bias and variance constants depend on unknown features of the data generating process. Specifically, $\mathsf{B}_j^{\mathtt{BR}}(\x)$ depends on the regression function $\theta_0$, the Lebesgue density of $X$, and their derivatives at $X=\x$ while $\mathsf{V}_j^{\mathtt{BR}}(\x)$ is determined by the density of $X$ and the conditional variance of the regression error $\varepsilon=Y-\theta_0(X)$ at $X=\x$. To implement the construction of the step size, we posit a simple parametric model:
\begin{equation*}
    \E[Y|X] = \gamma_0 + \sum_{k=1}^5 \gamma_k (X-\x_0)^k, \quad X\sim\mathrm{Normal}(\mu,\sigma^2)
\end{equation*}
where $\{\gamma_0,\gamma_1,\gamma_2,\gamma_4,\gamma_4,\gamma_5,\mu,\sigma\}$ are parameters to be estimated. Once we estimate the parameters of this reference model, we can construct a rule-of-thumb step size $\epsilon_{j,n}^{\mathtt{ROT}}$ by replacing $\mathsf{B}_j^{\mathtt{BR}}(\x)$ and $\mathsf{V}_j^{\mathtt{BR}}(\x)$ with their estimates. Note that although the bias and variance constant estimators may not be consistent for the true $\mathsf{B}_j^{\mathtt{BR}}(\x)$ and $\mathsf{V}_j^{\mathtt{BR}}(\x)$, the rate of $\epsilon_{j,n}^{\mathtt{ROT}}$ is MSE-optimal, and the numerical derivative estimator converges to $\mathcal{D}_j(\x)$ sufficiently fast to satisfy Equation \eqref{M-Rate of convergence condition} in the main paper.

\bibliographystyle{imsart-nameyear.bst} 
\bibliography{CJN_2024_AOS--bib}

\makeatletter\@input{xx.tex}\makeatother